МИНИСТЕРСТВО ОБРАЗОВАНИЯ И НАУКИ УКРАИНЫ
ЗАПОРОЖСКИЙ НАЦИОНАЛЬНЫЙ УНИВЕРСИТЕТ

С. В. Курапов
М. В. Давидовский
С. И. Полюга

# АЛГОРИТМИЧЕСКИЕ МЕТОДЫ КОНЕЧНЫХ ДИСКРЕТНЫХ СТРУКТУР

# ГАМИЛЬТОНОВ ЦИКЛ ПОЛНОГО ГРАФА И ЗАДАЧА КОММИВОЯЖЁРА

**(на правах рукописи)**

**Курапов С.В.**

**Давидовский М.В.**

**Полюга С.И.**

Рассматривается задача построения гамильтонового цикла в полном графе. Устанавливается правило построения гамильтонового цикла на базе изометрических циклов графа. Представлен алгоритм построения гамильтонового цикла на основе кольцевого суммирования изометрических циклов графа. На основе матрицы расстояний между вершинами определяется вес каждого цикла как аддитивная сумма весов его ребер. Для построения оптимального маршрута графа используется идея нахождения оптимального маршрута между четырьмя вершинами. Дальнейшие последовательные построения направлены на присоединение соприкасающегося изометрического цила с увеличением количества вершин на одну единицу. Рекурсивный процесс продолжается до подключения всех вершин графа. В работе представлен полиномиальный алгоритм решения симметричной задачи коммивояжёра. Приведены примеры решения задачи.

Для научных работников, студентов и аспирантов, специализирующихся на применении методов прикладной математики.

## Содержание

## ВВЕДЕНИЕ

Задача коммивояжёра (в англоязычной литературе – задача TSP, сокр. от Traveling Salesman Problem – задача странствующего торговца) является одной из классических задач комбинаторной оптимизации. В соотетствии с общепринятой постановкой задачи, необходимо найти оптимальный (с точки зрения некоторого критерия оптимальности) маршрут, проходящий через указанные точки (города) хотя бы по одному разу с возвратом в исходную точку. Критерием оптимальности маршрута в разных формулировках задачи может выступать минимальность расстояния, расходов, их совокупность и т. п., что задаётся соответствующей матрицей. Зачастую в формулировке задачи вводится также следующее условие – маршрут должен проходить через каждую точку один и только один раз. В случае такой постановки задачи выбор маршрутов осуществляется среди гамильтоновых циклов графа, вершинами которого являются указанные точки. Существуют также частные случаи общей постановки задачи коммивояжёра. В *геометрической постановке* (также называемой *планарной* или *евклидовой*) элементы матрицы расстояний представляют собой расстояния между точками на плоскости. В *метрической постановке* для элементов матрицы стоимостей выполняется неравенство треугольника. В зависимости от длительности или длины маршрутов и от направления движения различают *симметричную* и *асимметричную* постановки. В симметричной задаче все пары ребер между одними и теми же вершинами имеют одинаковую длину, в то время как ассиметричная задача представляется с помощью ориентированного графа. Также существует так называемая *обобщённая задача коммивояжёра* – задача комбинаторной оптимизации, заключающаяся в нахождении кратчайшего замкнутого пути, который проходит через одну вершину в каждом кластере, где кластеры определяются специальным разбиением множества вершин, а стоимость перехода из одной вершины в другую – соответствующим значением матрицы стоимостей.

Большинство постановок являются оптимизационными и относятся к классу NP-трудных задач. Кроме того, задача коммивояжёра относится к числу трансвычислительных, что означает что на практике она не может быть решена методом перебора уже при относительно небольшом (более 66) числе точек.

С задачей коммивояжёра связана другая классическая задача теории графов – задача о поиске гамильтонового цикла графа. Примером является известная математическая задача «О ходе коня». Эту задачу исследовал выдающийся математик Леонард Эйлер в работе «Решение одного любопытного вопроса, который, кажется, не подчиняется никакому исследованию». Другим примером может служить икосиан – описанная в 1859 г. Уильямом Гамильтоном математическая головомка «Путешествие вокруг света», в которой надо пройти по додекаэдру (графу с 20 вершинами), пройдя при этом по каждой вершине ровно один раз.

Впервые формулировка задачи, подарившая ей классическое название, по-видимому, встречается в изданном в 1832 г. справочнике для коммивояжёров под названием «Коммивояжёр – как он должен вести себя и что должен делать для того, чтобы доставлять товар и иметь успех в своих делах – советы старого курьера» (нем. Der Handlungsreisende – wie er sein soll und was er zu tun hat, um Aufträge zu erhalten und eines glücklichen Erfolgs in seinen Geschäften gewiß zu sein – von einem alten Commis-Voyageur). В этой книге задача рассматривается с практической точки зрения и приводятся примеры туров по Германии и Швейцарии, однако книга не содержит математического обоснования задачи и авторами не были предложены математические методы её решения.

Впервые задача коммивояжёра была представлена в терминах математической оптимизации Карлом Менгером в 1930 году в следующей формулировке: *мы называем задачей посыльного (поскольку этот вопрос возникает у каждого почтальона, в частности, её решают многие путешественники) задачу найти кратчайший путь между конечным множеством мест, расстояние между которыми известно*. Современная широко известная в англоязычной литературе постановка задачи под названием «*задача странствующего торговца*» (англ. traveling salesman problem) принадлежит профессору Принстонского университета Хасслеру Уитни (англ. Hassler Whitney).

Несмотря на простоту постановки задачи и возможность достаточно простого нахождения хороших решений, задача коммивояжёра интересна тем, что нахождение оптимального маршрута является сложной задачей. Поэтому исследование задачи коммивояжёра имеет не только практический смысл, но и важное теоретическое значение, как модель для создания и апробации новых математических алгоритмов. Так многие методы дискретной оптимизации (i.e., метод ветвей и границ, метод отсечений, различные эвристические алгоритмы и т. д.) были разработаны на примере задачи коммивояжёра. Задача исследовалась многими учеными как теоретически, так и с точки зрения её приложений в информатике, экономике, химии и биологии.

Важный вклад в исследование задачи внесли такие учёные как Джордж Данциг, Делберт Рей Фалкерсон и Селмер Джонсон, которые в 1954 г. сформулировали её в виде задачи дискретной оптимизации. Используя метод отсечений, они построили маршрут для одной частной постановки задачи с 49 городами, а также доказали его оптимальность. NP-трудность задачи была показана в 1972 году Ричардом Карпом как следствие из доказанной им NP-полноты задачи поиска гамильтоновых путей в графе. В свою очередь, Мартину Грётчелу, Манфреду Падбергу и Джованни Ринальди в 1970 – 1980 годах удалось получить решение для отдельного случая задачи с 2393 городами на основе новых методов деления плоскостью (ветвей и границ) [18].

В 1990-е годы для сравнения результатов работы различных групп исследователей был создан набор стандартизованных экземпляров задачи коммивояжёра различной степени сложности (TSPLIB). В марте 2005 года программой Конкорд была решена задача с 33 810 вершинами. Решение представляло собой оптимальный маршрут длиной 66 048 945. Кроме того, было доказано отсутствие более коротких маршрутов. В 2006 году было найдено решение для задачи с 85 900 вершинами.

Метрическая, симметричная и асимметричная постановки задачи коммивояжёра NP-эквивалентны. При условии, что $P \neq NP$, не существует алгоритма, который для некоторого полинома вычислял бы решения задачи коммивояжёра, отличающегося от оптимального решения на некоторый коэффициент. Все известные на сегодняшний день эффективные методы решения задачи коммивояжёра являются эвристическими, большинство из которых дают некоторое приближённое решение, а не оптимальное. На практике наиболее часто применяются метод ветвей и границ [18], эластичной сети [31], алгоритм муравьиной колонии [30], а также алгоритмы динамического программирования [28, 29] и генетические алгоритмы [25, 33].

## Глава 1. ОСНОВНЫЕ ПОНЯТИЯ И ОПРЕДЕЛЕНИЯ

Приведём основные понятия и определения, использующиеся в данной работе. Для удобства восприятия читателем, сопроводим их краткими пояснениями и иллюстрациями.

**Определение 1.1.** Связный неориентированный граф G, не имеющий мостов и точек сочленения, без петель и кратных рёбер и без вершин с локальной степенью меньшей или равной двум, называется *несепарабельным графом* G.

**Определение 1.2.** *Простой цикл графа* – связный замкнутый маршрут графа, в котором валентность всех вершин равна двум [5-7].

**Определение 1.3.** *Квазицикл графа* – замкнутый маршрут в графе, в котором валентность вершин кратна двум [5-7,21,23,27].

**Определение 1.4.** *Гамильтонов цикл графа* – простой цикл, проходящий через все вершины графа [27].

Количество вершин в графе G(E,V) обозначается латинской буквой $n$.

**Определение 1.5.** *Изометрическим циклом* в графе называется простой цикл, для которого кратчайший путь между любыми двумя его вершинами состоит из рёбер этого цикла. Изометрический цикл – частный случай изометрического подграфа [32].

Другими словами, изометрическим циклом в графе называется подграф G′ в виде простого цикла, если между двумя любыми несмежными вершинами данного подграфа в соответствующем графе G не существует маршрутов меньшей длины, чем маршруты, принадлежащие данному циклу.

Граф $G^* = (V^*,E^*)$ называется частью графа G = (V,E), если $V^* \subseteq V$ и $U^* \subseteq U$ т.е. часть графа образуется из исходного графа удалением некоторых вершин и рёбер.

Особо важную роль играют два типа частей графа: подграф и суграф [27].

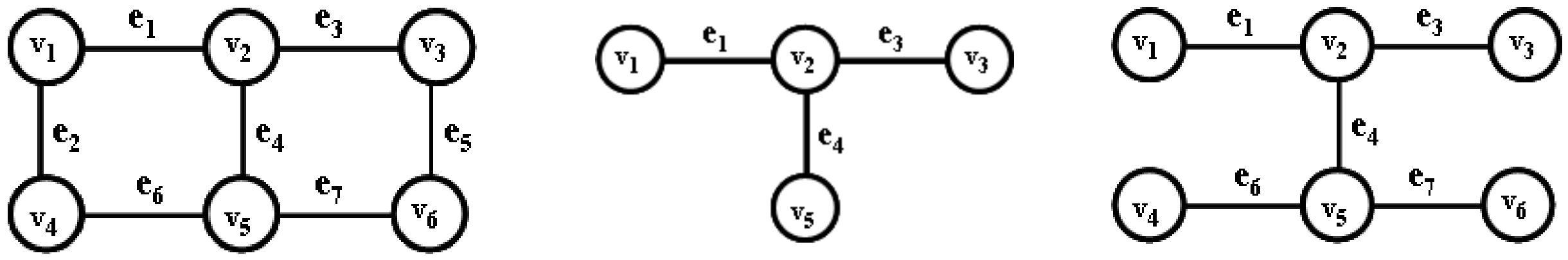

а) исходный граф б) подграф в) суграф

Рис. 1.1. Получение частей графа из исходного графа.

**Определение 1.6.** Часть $G^* = (V^*,E^*)$ называют *подграфом* графа G = (V,E), если $E^* = \{xy \in E / x,y \in V^*\}$. Подграф образуется из исходного графа некоторым количеством выделенных вершин и некоторым количеством инцидентных им рёбер (рис. 1.1(б)).

**Определение 1.7.** Часть графа $G^* = (V^*,E^*)$ называют *суграфом* графа G = (X,U) если $X^* = X$, т.е. суграф образуется из исходного графа удалением только рёбер, без удаления вершин (рис. 1.1(в)).

Введём несколько операций над суграфами графа G (рис. 1.2(а)). Рассмотрим суграфы $G_1 = (V,E_1)$ и $G_2 = (V, E_2)$ (рис. 1.2(б) и рис. 1.2(в)).

**Определение 1.8.** *Объединение* суграфов $G_1$ и $G_2$, обозначаемое как $G_1 \cup G_2$, представляет собой такой суграф $G_3= (V, E_1 \cup E_2)$, что множество его рёбер является объединением $E_1$ и $E_2$. Например, суграфы $G_1$ и $G_2$ и их объединение представлено на рис. 1.2(г).

**Определение 1.9.** *Пересечение* суграфов $G_1$ и $G_2$, обозначаемое как $G_1 \cap G_2$, представляет собой граф $G_3= (V,E_1 \cap E_2)$. Таким образом, множество рёбер $G_3$ состоит только из рёбер, присутствующих одновременно в $G_1$ и $G_2$. Пересечение суграфов $G_1$ и $G_2$ показано на рис. 1.2(д).

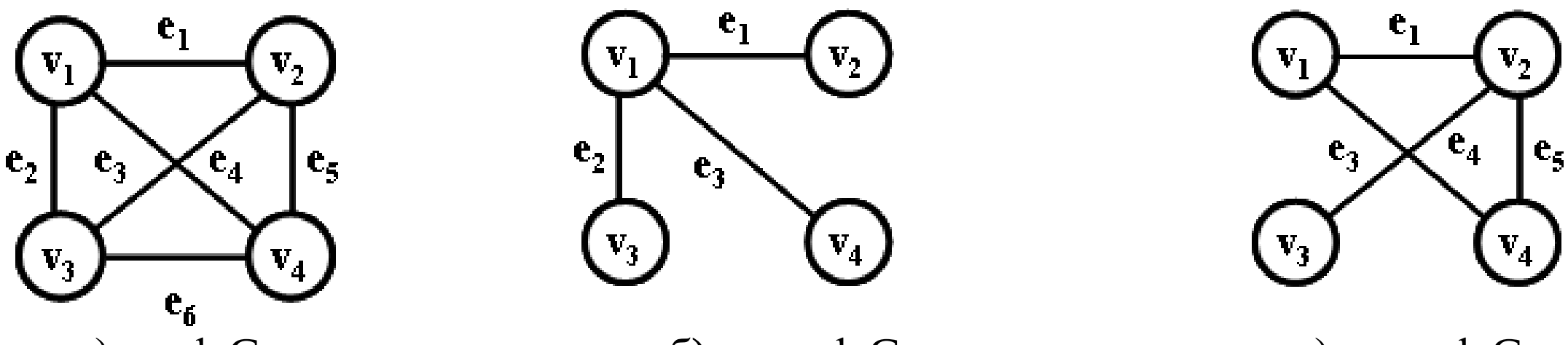

а) граф G. б) суграф $G_1$. в) суграф $G_2$.

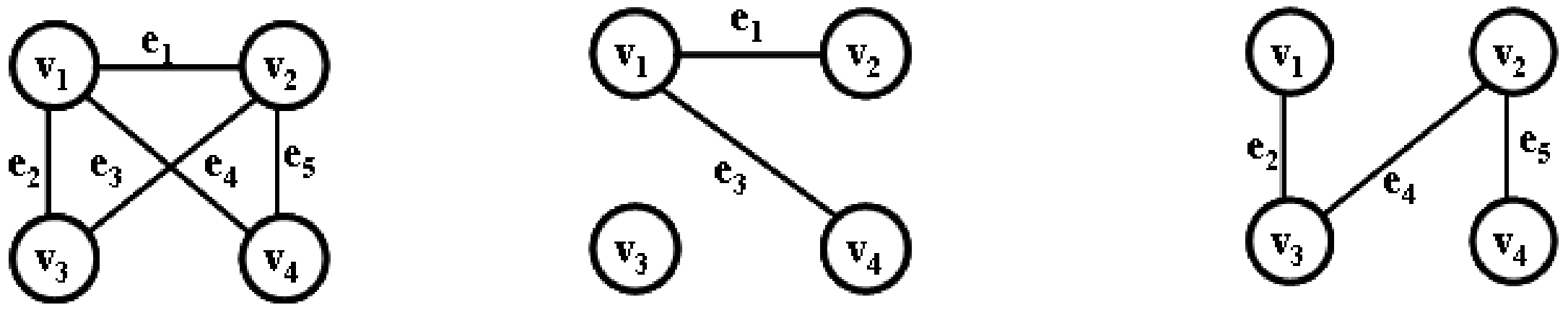

г) объединение суграфов $G_1$ и $G_2$. д) пересечение суграфов $G_1$ и $G_2$. е) кольцевая сумма суграфов $G_1$ и $G_2$.

Рис. 1.2. Суграфы.

**Определение 1.10.** *Кольцевая сумма* двух суграфов $G_1$ и $G_2$, обозначаемая как $G_1 \oplus G_2$, представляет собой суграф $G_3$, порождённый на множестве рёбер $E_1 \oplus E_2$. Другими словами, суграф $G_3$ состоит только из рёбер, присутствующих либо в $G_1$, либо в $G_2$, но не в обоих суграфах одновременно [23]. Кольцевая сумма суграфов показана на рис. 1.2(е).

Легко убедиться в том, что рассмотренные операции коммутативны, т.е. $G_1 \cup G_2 = G_2 \cup G_1$, $G_1 \cap G_2 = G_2 \cap G_1$, $G_1 \oplus G_2 = G_2 \oplus G_1$.

Заметим также, что эти операции бинарны, т. е. определены по отношению к двум суграфам.

Таким образом, появляется операция сложения суграфов (будем называть её *кольцевой суммой* [24,27]). Эта операция отлична от известной арифметической операции сложения:

$$(V,E_1;P) \oplus (V,E_2;P) = (V,(E_1 \cup E_2)\backslash(E_1 \cap E_2);P). \quad (1.1)$$

Пусть $G = (V,E;P)$ – граф с пронумерованным множеством рёбер $E = \{e_1,e_2,...,e_m\}$, тогда суграф графа можно представлять в виде характеристического вектора с коэффициентами

$$\alpha_i = \begin{cases} 1, \text{если } e_i \in E; \\ 0, \text{если } e_i \notin E. \end{cases}$$

Следующее множество элементов <1,0,....0>, <0,1,.....,0>,....,<0,0,....,1> представляет собой характеристические вектора для однорёберных суграфов, записанные в виде кортежа. Размерность характеристического вектора определяется количеством рёбер графа.

В качестве примера рассмотрим граф G представленный на рис. 1.2(а). Характеристический вектор графа G = <1,1,1,1,1,1> можно записать в виде множества элементов $\{e_1,e_2,e_3,e_4,e_5,e_6\}$, рассматривая каждый элемент множества как результат произведения на соответствующий коэффициент из характеристического вектора. Последовательность записи коэффициентов в характеристическом векторе осуществляется справа налево.

Суграф $G_1$ представленный на рис. 1.2(б) запишется в виде характеристического вектора <0,0,0,1,1,1> или в виде подмножества элементов $\{e_1,e_2,e_3\}$. Суграф $G_2$ представленный на рис. 1.2(в) запишется в виде характеристического вектора <0,1,1,1,0,1> или в виде подмножества элементов $\{e_1,e_3,e_4,e_5\}$. Суграф $G_1 \cup G_2$ представленный на рис. 1.2(г) запишется в виде характеристического вектора <0,1,1,1,1,1> или в виде подмножества элементов $\{e_1,e_2,e_3,e_4,e_5\}$. Суграф $G_1 \cap G_2$ представленный на рис. 1.2(д) запишется в виде характеристического вектора <0,0,0,1,0,1> или в виде подмножества элементов $\{e_1,e_3\}$. Суграф $G_1 \oplus G_2$ представленный на рис. 1.2(е) запишется в виде характеристического вектора <0,1,1,0,1,0> или в виде подмножества элементов $\{e_2,e_4,e_5\}$.

$$G_1 \oplus G_2 = <0,0,0,1,1,1> + <0,1,1,1,0,1> = <0,1,1,0,1,0>$$

или

$$\{e_1,e_2,e_3) \oplus \{e_1,e_3,e_4,e_5\} = \{e_2,e_4,e_5\};$$

$$G_1 \oplus G_2 = (G_1 \cup G_2) \backslash (G_1 \cap G_2) = <0,1,1,1,1,1> + <0,0,0,1,0,1> = <0,1,1,0,1,0>$$

или

$$\{e_1,e_2,e_3,e_4,e_5) \oplus \{e_1,e_3\} = \{e_2,e_4,e_5\}.$$

**Определение 1.11.** *Множество изометрических циклов* графа является инвариантом графа, обозначается как $C_\tau$ и представляет собой кортеж:

$C_\tau$ = $<c_1,c_2,c_3,\ldots,c_k>$, где $k$ – количество циклов в $C_\tau$.

Рассмотрим несепарабельный граф $G_1$ (рис. 1.3).

Количество вершин графа = 10.
Количество рёбер графа = 20.
Количество изометрических циклов = 16.

Инцидентность графа:

ребро $e_1$: $(v_1,v_2)$ или $(v_2,v_1)$;
ребро $e_2$: $(v_1,v_3)$ или $(v_3,v_1)$;
ребро $e_3$: $(v_1,v_4)$ или $(v_4,v_1)$;
ребро $e_4$: $(v_1,v_7)$ или $(v_7,v_1)$;
ребро $e_5$: $(v_2,v_3)$ или $(v_3,v_2)$;
ребро $e_6$: $(v_2,v_5)$ или $(v_5,v_2)$;
ребро $e_7$: $(v_2,v_7)$ или $(v_7,v_2)$;
ребро $e_8$: $(v_3,v_7)$ или $(v_7,v_3)$;
ребро $e_9$: $(v_3,v_9)$ или $(v_9,v_3)$;
ребро $e_{10}$: $(v_4,v_5)$ или $(v_5,v_4)$;
ребро $e_{11}$: $(v_4,v_8)$ или $(v_8,v_4)$;
ребро $e_{12}$: $(v_4,v_{10})$ или $(v_{10},v_4)$;
ребро $e_{13}$: $(v_5,v_6)$ или $(v_6,v_5)$;
ребро $e_{14}$: $(v_5,v_8)$ или $(v_8,v_5)$;
ребро $e_{15}$: $(v_6,v_8)$ или $(v_8,v_6)$;
ребро $e_{16}$: $(v_6,v_9)$ или $(v_9,v_6)$;
ребро $e_{17}$: $(v_6,v_{10})$ или $(v_{10},v_6)$;
ребро $e_{18}$: $(v_7,v_9)$ или $(v_9,v_7)$;
ребро $e_{19}$: $(v_8,v_{10})$ или $(v_{10},v_8)$;
ребро $e_{20}$: $(v_9,v_{10})$ или $(v_{10},v_9)$.

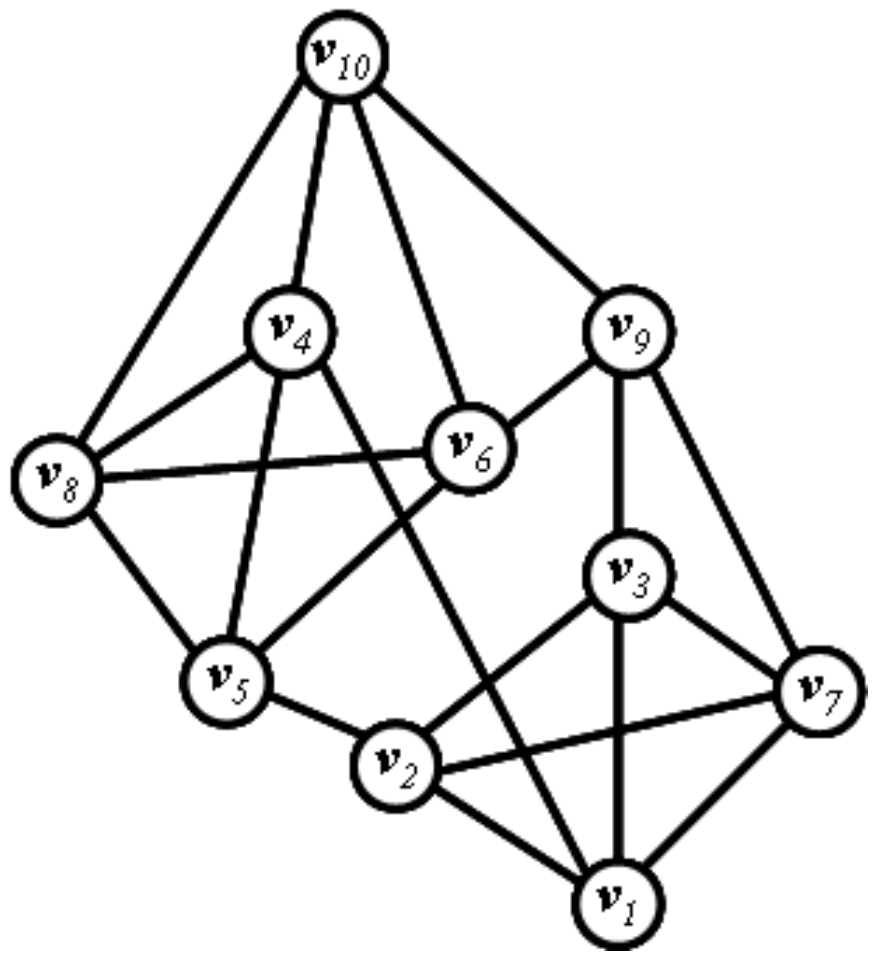

Рис 1.3. Несепарабельный граф $G_1$.

Множество изометрических циклов графа $G_1$:

цикл $c_1 = \{e_1,e_2,e_5\} \leftrightarrow \{v_1,v_2,v_3\}$;
цикл $c_2 = \{e_1,e_3,e_6,e_{10}\} \leftrightarrow \{v_1,v_2,v_4,v_5\}$;
цикл $c_3 = \{e_1,e_4,e_7\} \leftrightarrow \{v_1,v_2,v_7\}$;
цикл $c_4 = \{e_2,e_4,e_8\} \leftrightarrow \{v_1,v_3,v_7\}$;
цикл $c_5 = \{e_2,e_3,e_9,e_{12},e_{20}\} \leftrightarrow \{v_1,v_3,v_4,v_9,v_{10}\}$;
цикл $c_6 = \{e_3,e_4,e_{12},e_{18},e_{20}\} \leftrightarrow \{v_1,v_4,v_7,v_{10},v_9\}$;
цикл $c_7 = \{e_5,e_7,e_8\} \leftrightarrow \{v_2,v_3,v_7\}$;
цикл $c_8 = \{e_5,e_6,e_9,e_{13},e_{16}\} \leftrightarrow \{v_2,v_3,v_5,v_9,v_6\}$;
цикл $c_9 = \{e_6,e_7,e_{13},e_{16},e_{18}\} \leftrightarrow \{v_2,v_5,v_7,v_6,v_9\}$;
цикл $c_{10} = \{e_8,e_9,e_{18}\} \leftrightarrow \{v_3,v_7,v_9\}$;
цикл $c_{11} = \{e_{10},e_{11},e_{14}\} \leftrightarrow \{v_4,v_5,v_8\}$;
цикл $c_{12} = \{e_{10},e_{12},e_{13},e_{17}\} \leftrightarrow \{v_4,v_5,v_{10},v_6\}$;
цикл $c_{13} = \{e_{11},e_{12},e_{19}\} \leftrightarrow \{v_4,v_8,v_{10}\}$;
цикл $c_{14} = \{e_{13},e_{14},e_{15}\} \leftrightarrow \{v_5,v_6,v_8\}$;
цикл $c_{15} = \{e_{15},e_{17},e_{19}\} \leftrightarrow \{v_6,v_8,v_{10}\}$;
цикл $c_{16} = \{e_{16},e_{17},e_{20}\} \leftrightarrow \{v_6,v_9,v_{10}\}$. (1.2)

**Определение 1.12.** *Вектором циклов по ребру* называется кортеж количества циклов, проходящих по ребру, для данного подмножества циклов $Y(G)$.

Например, множество изометрических циклов графа $G_1$ характеризуется следующим вектором $P_e$:

$P_e = <p_1.p_2.p_3,p_4,p_5,p_6,p_7,p_8,p_9,p_{10},p_{11},p_{12},p_{13},p_{14},p_{15},p_{16},p_{17},p_{18},p_{19}.p_{20}> =$

$=<5,5,6,3,3,3,3,3,3,3,2,4,4,2,2,3,3,3,2,3>$

Подмножество циклов Y(G) может характеризоваться значением квадратичного функционала Маклейна [20] без учёта нулевых элементов в векторе $P_e$

$$F_1(Y) = \sum_{i=1}^{m} p_i^2 - 3\sum_{i=1}^{m} p_i + 2m, \quad (1.3)$$

или значением кубического функционала Маклейна, учитывающего наличие нулевых элементов в векторе $P_e$

$$F_2(Y) = \sum_{i=1}^{m} p_i^3 - 3\sum_{i=1}^{m} p_i^2 + 2\sum_{i=1}^{m} p_i, \quad (1.4)$$

где $p_i$ – количество циклов, проходящих по ребру $e_i$.

Для множества изометрических циклов графа $G_1$ значение квадратичного функционала Маклейна равно 78, а значение кубического функционала Маклейна равно 354.

**Определение 1.13.** *Вращение вершины* – циклический порядок обхода инцидентных рёбер (представляется диаграммой вращения) [24].

**Определение 1.14.** *Топологический рисунок* плоского несепарабельного графа – это система вращений вершин, индуцирующая (порождающая) систему простых циклов, или система простых циклов, индуцирующая вращение вершин.

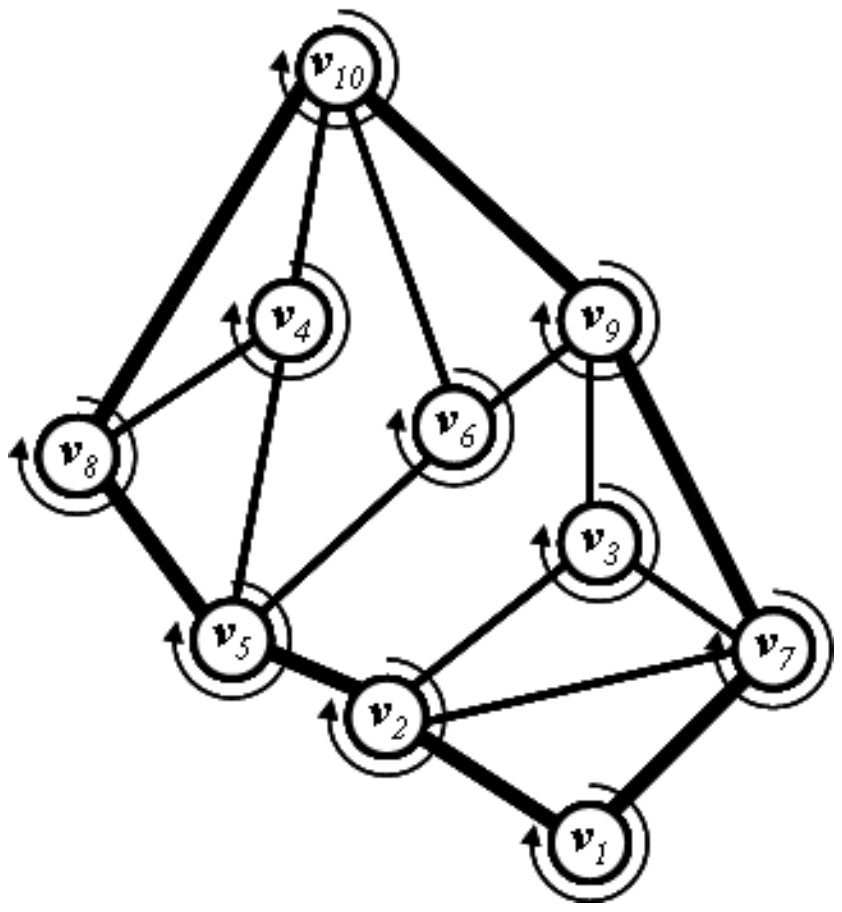

Рис. 1.4. Топологический рисунок плоской части графа и обод.

Вращение вершин в топологическом рисунке:

$\sigma(v_1)$: $v_2$ $v_7$

$\sigma(v_2)$: $v_1$ $v_5$ $v_3$ $v_7$

$\sigma(v_3)$: $v_2$ $v_9$ $v_7$

$\sigma(v_4)$: $v_8$ $v_{10}$ $v_5$

$\sigma(v_5)$: $v_8$ $v_4$ $v_6$ $v_2$

$\sigma(v_6)$: $v_9$ $v_5$ $v_{10}$

$\sigma(v_7)$: $v_1$ $v_2$ $v_3$ $v_9$

$\sigma(v_8)$: $v_4$ $v_5$ $v_{10}$

$\sigma(v_9)$: $v_7$ $v_3$ $v_6$ $v_{10}$

$\sigma(v_{10})$: $v_9$ $v_6$ $v_4$ $v_8$

Индуцируемая вращением вершин система циклов имеет следующий вид:

$c_3 = \{e_1,e_4,e_7\} \leftrightarrow <v_1,v_7,v_2> \leftrightarrow (v_1,v_7)+(v_7,v_2)+(v_2,v_1)$;
$c_7 = \{e_5,e_7,e_8\} \leftrightarrow <v_2,v_7,v_3> \leftrightarrow (v_2,v_7)+(v_7,v_3)+(v_3,v_2)$;
$c_8 = \{e_5,e_6,e_9,e_{13},e_{16}\} \leftrightarrow <v_2,v_3,v_9,v_6,v_5> \leftrightarrow$
$\leftrightarrow (v_2,v_3)+(v_3,v_9)+(v_9,v_6) +(v_6,v_5)+(v_5,v_2)$;
$c_{10} = \{e_8,e_9,e_{18}\} \leftrightarrow <v_3,v_7,v_9> \leftrightarrow (v_3,v_7)+(v_7,v_9)+(v_9,v_3)$;
$c_{11} = \{e_{10},e_{11},e_{14}\} \leftrightarrow <v_4,v_8,v_5> \leftrightarrow (v_4,v_8)+(v_8,v_5)+(v_5,v_4)$;
$c_{12} = \{e_{10},e_{12},e_{13},e_{17}\} \leftrightarrow <v_4,v_5,v_6,v_{10}> \leftrightarrow (v_4,v_5)+(v_5,v_6)+(v_6,v_{10}) +(v_{10},v_4)$;
$c_{13} = \{e_{11},e_{12},e_{19}\} \leftrightarrow \{v_4,v_{10},v_8\} \leftrightarrow (v_4,v_{10})+(v_{10},v_8)+(v_8,v_4)$;

$c_{16} = \{e_{16},e_{17},e_{20}\} \leftrightarrow <v_6,v_9,v_{10}> \leftrightarrow (v_6,v_9)+(v_9,v_{10})+(v_{10},v_6)$.

В топологическом рисунке цикл может быть записан не только в виде множества рёбер и вершин, но и в виде циклического порядка ориентированных рёбер, называемого векторной записью цикла.

$c_{12} = (v_4,v_5)+(v_5,v_6)+(v_6,v_{10}) +(v_{10},v_4)$.

**Определение 1.15.** *Ободом* подмножества изометрических циклов графа называется кольцевая сумма всех изометрических циклов суграфа.

Например, обод плоской части графа (выделен жирной линией), представлен на рис. 1.4:

$c_0 = c_3 \oplus c_7 \oplus c_8 \oplus c_{10} \oplus c_{11} \oplus c_{12} \oplus c_{13} \oplus c_{16} = \{e_1,e_4,e_6,e_{14},e_{18},e_{19},e_{20}\} \leftrightarrow$

$\leftrightarrow <v_1,v_2,v_5,v_8,v_{10},v_9,v_7,v_1> \leftrightarrow (v_1,v_2)+(v_2,v_5)+(v_5,v_8)+(v_8,v_{10})+(v_{10},v_{90})+(v_9,v_7)+(v_7,v_1)$.

Кольцевая сумма независимых циклов и обода в топологическом рисунке графа есть пустое множество:

$$\sum_{i=1}^{k} c_i \oplus c_0 = \varnothing \qquad (2.18)$$

Перейдём к рассмотрению полных графов $K_n$.

Количество рёбер в полном графе определим по формуле:

$$m = \frac{n(n-1)}{2}. \qquad (2.19)$$

Количество изометрических (треугольных) циклов:

$$k_c = \frac{n(n-1)(n-2)}{6} \qquad (2.20)$$

Введём следующее определение гамильтонового цикла.

**Определение 1.16**. *Гамильтоновым циклом в полном графе*, будем называть простой цикл обода плоского суграфа, образованный линейной комбинацией изометрических циклов графа, значение кубического функционала Маклейна которого равно нулю.

# Глава 2. ГАМИЛЬТОНОВ ЦИКЛ ПОЛНОГО ГРАФА

## 2.1. Принцип построения гамильтонового цикла в полном графе

Рассмотрим процесс образования гамильтонового цикла полного графа. Построение будем проводить на примере полного графа $K_5$.

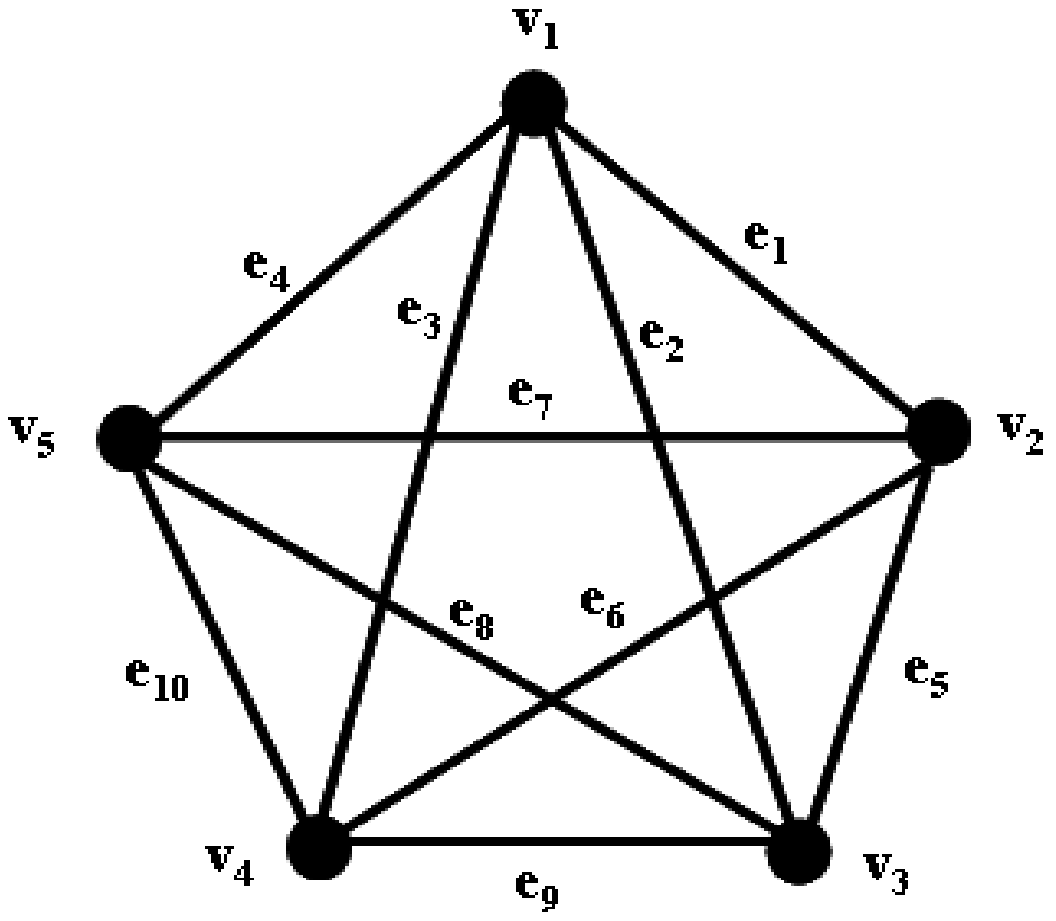

Рис. 2.1. Полный граф $K_5$.

Множество изометрических (треугольных) циклов графа $K_5$.

цикл $c_1 = \{e_1,e_2,e_5\} \leftrightarrow \{v_1,v_2,v_3\}$;
цикл $c_2 = \{e_1,e_3,e_6\} \leftrightarrow \{v_1,v_2,v_4\}$;
цикл $c_3 = \{e_1,e_4,e_7\} \leftrightarrow \{v_1,v_2,v_5\}$;
цикл $c_4 = \{e_2,e_3,e_9\} \leftrightarrow \{v_1,v_3,v_4\}$;
цикл $c_5 = \{e_2,e_4,e_8\} \leftrightarrow \{v_1,v_3,v_5\}$;
цикл $c_6 = \{e_3,e_4,e_{10}\} \leftrightarrow \{v_1,v_4,v_5\}$;
цикл $c_7 = \{e_5,e_6,e_9\} \leftrightarrow \{v_2,v_3,v_4\}$;
цикл $c_8 = \{e_5,e_7,e_8\} \leftrightarrow \{v_2,v_3,v_5\}$;
цикл $c_9 = \{e_6,e_7,e_{10}\} \leftrightarrow \{v_2,v_4,v_5\}$;
цикл $c_{10} = \{e_8,e_9,e_{10}\} \leftrightarrow \{v_3,v_4,v_5\}$.

Вектор количества проходящих по рёбрам циклов для множества циклов $C_\tau$ имеет вид:

$P_e = <3,3,3,3,3,3,3,3,3,3>$.

Вектор количества проходящих по вершинам циклов для множества циклов $C_\tau$ имеет вид:

$P_v = <6,6,6,6,6>$.

Упрощенно вектор $P_e$ будем называть вектором циклов по рёбрам, а вектор $P_v$ – вектором циклов по вершинам.

Значение кубического функционала Маклейна для множества изометрических циклов графа $C_\tau$ равно 60.

В соответствии с принципом построения гамильтонового цикла в полном графе исключаем из множества изометрических циклов $C_\tau$ цикл $c_1$.

После исключения цикла $c_1$, вектор циклов по ребрам для подмножества циклов $\frac{\delta C_\tau}{\delta c_1}$

имеет вид: $P_e$ = <2,2,3,3,2,3,3,3,3,3>.

После исключения цикла $c_1$, вектор циклов по вершинам для подмножества циклов $\frac{\delta C_\tau}{\delta c_1}$ имеет вид: $P_v$ = <5,5,5,6,6>.

Значение кубического функционала Маклейна для подмножества изометрических циклов графа $F(\frac{\delta C_\tau}{\delta c_1}) = 42$.

$\frac{\delta C_\tau}{\delta c_1}$ – обозначение процесса исключения цикла $c_1$ из множества изометрических циклов $C_\tau$ согласно [2].

Исключаем из множества $C_\tau$ цикл $c_6$.

После исключения цикла $c_6$ вектор циклов по ребрам для подмножества циклов $\frac{\delta^2 C_\tau}{\delta c_1 \delta c_6}$ имеет вид:

$P_e$ = <2,2,2,2,2,3,3,3,3,2>.

После исключения цикла $c_6$ вектор циклов по вершинам для подмножества циклов $\frac{\delta^2 C_\tau}{\delta c_1 \delta c_6}$ имеет вид:

$P_v$ = <4,5,5,5,5>.

Значение кубического функционала Маклейна для подмножества изометрических циклов графа $F(\frac{\delta^2 C_\tau}{\delta c_1 \delta c_6}) = 24$.

Исключаем из множества цикл $c_8$.

После исключения цикла $c_8$ вектор циклов по ребрам для подмножества циклов $\frac{\delta^3 C_\tau}{\delta c_1 \delta c_6 \delta c_8}$ имеет вид:

$P_e$ = <2,2,2,2,1,3,2,2,3,2>.

После исключения цикла $c_8$ вектор циклов по вершинам для подмножества циклов $\frac{\delta^3 C_\tau}{\delta c_1 \delta c_6 \delta c_8}$ имеет вид:

$P_v$ = <4,4,4,5,4>.

Значение кубического функционала Маклейна для подмножества изометрических циклов графа F( $\frac{\delta^3 C_\tau}{\delta c_1 \delta c_6 \delta c_8}$ ) = 12.

Исключаем из множества цикл $c_2$.

После исключения цикла $c_2$ вектор циклов по ребрам для подмножества циклов $\frac{\delta^4 C_\tau}{\delta c_1 \delta c_6 \delta c_8 \delta c_2}$ имеет вид:

$P_e$ = <1,2,1,2,1,2,2,2,3,2>.

После исключения цикла $c_2$ вектор циклов по вершинам для подмножества циклов $\frac{\delta^4 C_\tau}{\delta c_1 \delta c_6 \delta c_8 \delta c_2}$ имеет вид:

$P_v$ = <3,3,4,4,4>.

Значение кубического функционала Маклейна для подмножества изометрических циклов графа F( $\frac{\delta^4 C_\tau}{\delta c_1 \delta c_6 \delta c_8 \delta c_2}$ ) = 6.

Количество циклов в оставшейся системе равно цикломатическому числу графа $\vartheta(G) = m - n + 1$. Данная система циклов линейно независима и характеризует базис линейного подпространства циклов C пространства суграфов графа G.

Продолжим процесс исключения циклов. Обозначим базис Y(G) = $\frac{\delta^4 C_\tau}{\delta c_1 \delta c_6 \delta c_8 \delta c_2}$

Исключаем из множества цикл $c_4$.

После исключения цикла $c_4$ вектор циклов по ребрам для подмножества циклов $\frac{\partial Y(G)}{\partial c_4}$ имеет вид:

$P_e$ = <1,1,0,2,1,2,2,2,2,2>.

После исключения цикла $c_4$ вектор циклов по вершинам для подмножества циклов $\frac{\partial Y(G)}{\partial c_4}$ имеет вид:

$P_v$ = <2,3,3,3,4>.

Значение кубического функционала Маклейна для подмножества изометрических циклов графа F( $\frac{\partial Y(G)}{\partial c_4}$ ) = 0.

Так как кубический функционал Маклейна равен нулю, то подмножество циклов характеризует плоский суграф с удалённым ребром $e_3$ (рис. 2.2).

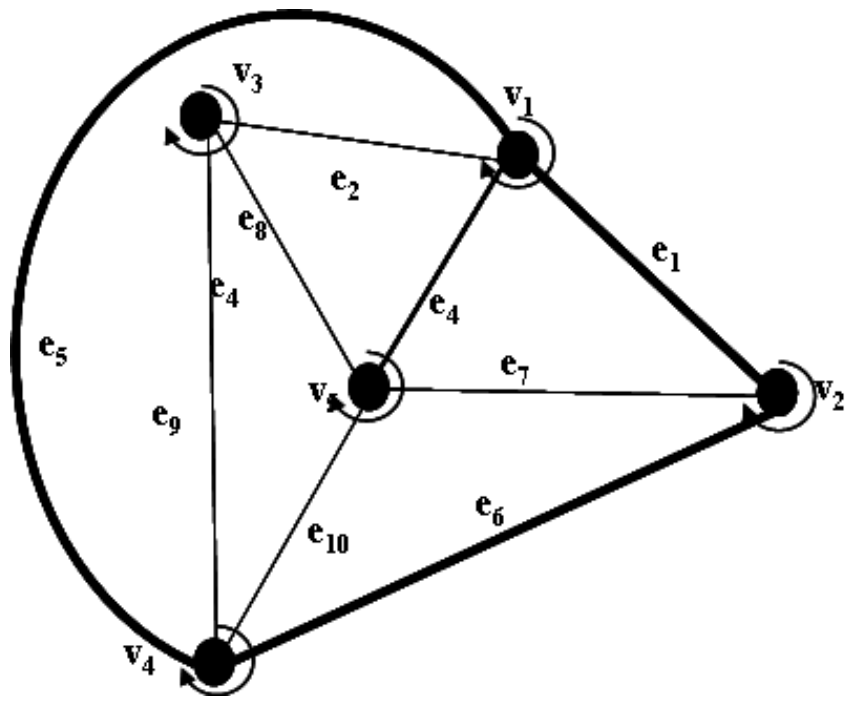

Рис. 2.2. Топологический рисунок плоского суграфа без ребра $e_3$.

Исключаем из подмножества $\frac{\partial Y(G)}{\partial c_4}$ цикл $c_3$.

После исключения цикла $c_3$ вектор циклов по ребрам для подмножества циклов $\frac{\partial^2 Y(G)}{\partial c_4 \partial c_3}$ имеет вид:

$P_e$ = <0,1,0,1,1,2,1,2,2,2>.

После исключения цикла $c_3$ вектор циклов по вершинам для подмножества циклов $\frac{\partial^2 Y(G)}{\partial c_4 \partial c_3}$ имеет вид:

$P_v$ = <1,2,3,3,3>.

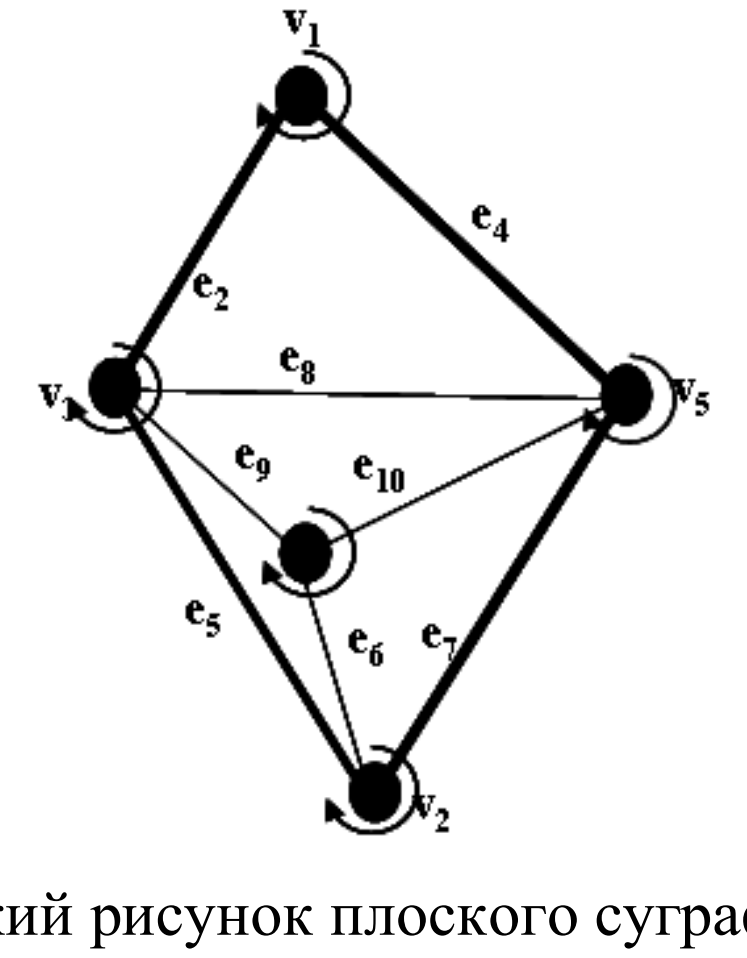

Рис. 2.3. Топологический рисунок плоского суграфа без рёбер $c_3$ и $c_4$.

Значение кубического функционала Маклейна для подмножества изометрических циклов графа $F(\frac{\partial^2 Y(G)}{\partial c_4 \partial c_3}) = 0$.

Так как кубический функционал Маклейна равен нулю, то подмножество циклов характеризует плоский суграф с удаленными ребрами $e_3$ и $c_4$ (рис. 2.3).

Исключаем из подмножества $\frac{\partial^2 Y(G)}{\partial c_4 \partial c_3}$ цикл $c_9$.

После исключения цикла $c_9$ вектор циклов по ребрам для подмножества циклов $\frac{\partial^3 Y(G)}{\partial c_4 \partial c_3 \partial c_9}$ имеет вид:

$P_e$ = <0,1,0,1,1,1,0,2,2,1>.

После исключения цикла $c_9$ вектор циклов по вершинам для подмножества циклов $\frac{\partial^3 Y(G)}{\partial c_4 \partial c_3 \partial c_9}$ имеет вид:

$P_v$ = <1,1,3,2,2>.

Значение кубического функционала Маклейна для подмножества изометрических циклов графа $F(\frac{\partial^3 Y(G)}{\partial c_4 \partial c_3 \partial c_9}) = 0$.

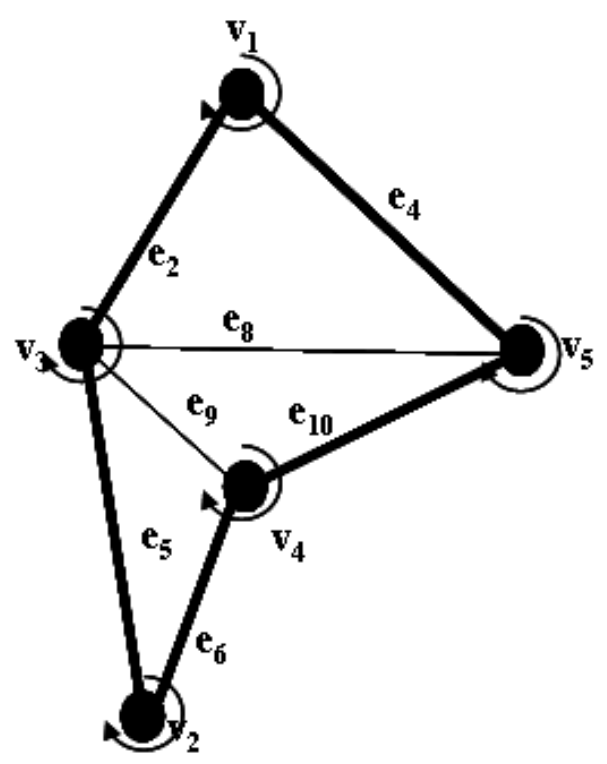

Рис. 2.4. Топологический рисунок плоского суграфа и обод.

В результате исключения образована следующая система циклов:

цикл $c_5 = \{e_2,e_4,e_8\} \leftrightarrow \{v_1,v_3,v_5\}$;
цикл $c_7 = \{e_5,e_6,e_9\} \leftrightarrow \{v_2,v_3,v_4\}$;
цикл $c_{10} = \{e_8,e_9,e_{10}\} \leftrightarrow \{v_3,v_4,v_5\}$.

Данное подмножество циклов описывает топологический рисунок плоского суграфа с тремя удаленными ребрами. Обод плоского суграфа является простым циклом, проходящим по всем вершинам графа (рис. 2.4).

Заметим, что количество единиц в векторе $P_e$ равно количеству вершин в гамильтоновом цикле. С другой стороны, в векторе $P_e$ единицам соответствуют ребра, входящие в обод суграфа. Количество двоек в $P_e$ определяет количество пересекающихся рёбер, а количество нулей определяет количество исключенных ребер в процессе планаризации графа.

Данное построение гамильтонового цикла можно осуществить и другим способом. Для этого достаточно осуществлять построение гамильтонового цикла в соответствии с правилом:

к предыдущему простому циклу можно последовательно присоединять только один изометрический цикл, имеющий одно пересекающееся ребро и одну отличную от предыдущих вершину.

**Определение 2.1**. Два цикла, имеющие в своем составе одно пересекающееся ребро и отличающиеся между собой одной вершиной, называются *соприкасающимися*.

Например,

$c_5 \oplus c_{10} = \{e_2,e_4,e_8\} \oplus \{e_8,e_9,e_{10}\} = \{e_2,e_4,e_9,e_{10}\} \leftrightarrow \{v_1,v_3,v_4,v_5\}$;

$c_5 \oplus c_{10} \oplus c_7 = \{e_2,e_4,e_9,e_{10}\} \oplus \{e_5,e_6,e_9\} = \{e_2,e_4,e_5,e_6,e_{10}\} \leftrightarrow \{v_1,v_2.v_3,v_4,v_5\}$.

Циклы $c_{10}$ и $c_5$ – соприкасающиеся. Цикл $c_7$ – соприкасающийся с простым циклом, образованным кольцевой суммой $c_5 \oplus c_{10}$. Данный процесс построения гамильтонового цикла для полного графа, использующий правило подключения изометрического цикла, представлен на рис. 2.5.

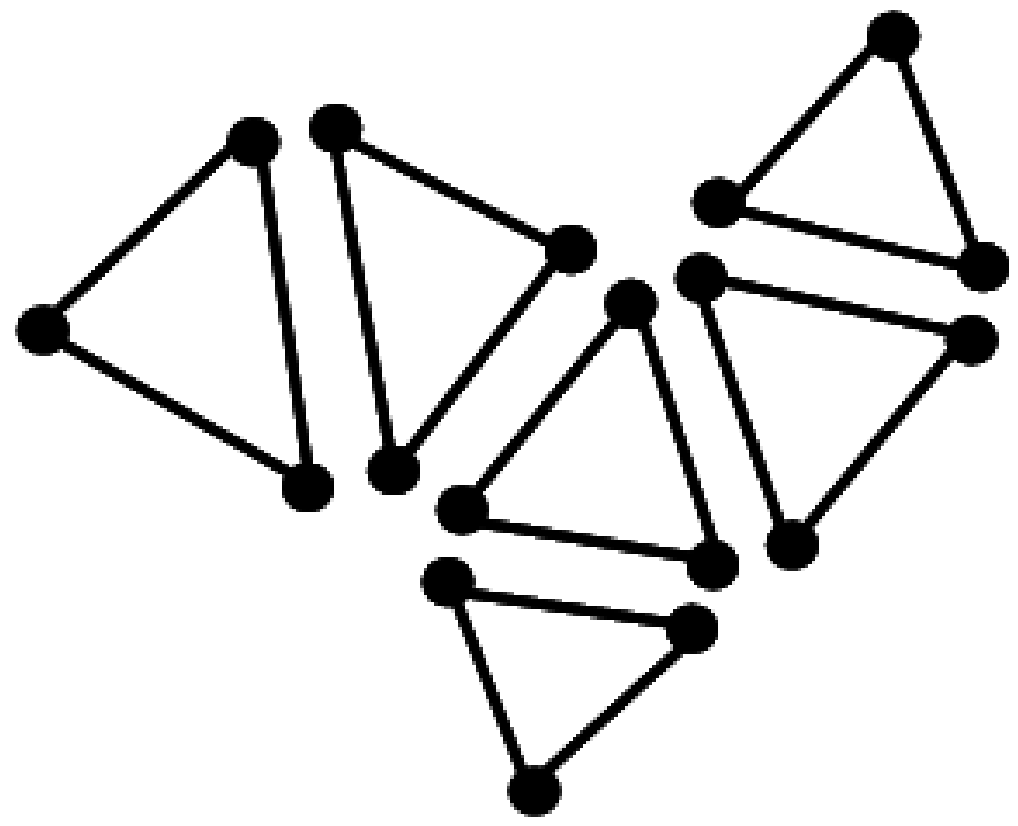

Рис. 2.5. Процесс построения гамильтонового цикла.

Если мы на последнем шаге процесса исключения циклов в качестве исключаемого цикла выберем не цикл $c_9$, а цикл $c_5$, то после исключения цикла $c_5$ вектор циклов по ребрам для подмножества циклов $\frac{\partial^3 Y(G)}{\partial c_4 \partial c_3 \partial c_5}$ будет иметь вид:

$P_e$ = <0,0,0,0,1,2,1,1,2,2>.

После исключения цикла $c_5$ вектор циклов по вершинам для подмножества циклов $\frac{\partial^2 Y(G)}{\partial c_4 \partial c_3 \partial c_5}$ имеет вид:

$P_v$ = <0,2,2,3,2>.

Значение кубического функционала Маклейна для множества изометрических циклов графа $F(\frac{\partial^3 Y(G)}{\partial c_4 \partial c_3 \partial c_5}) = 0$.

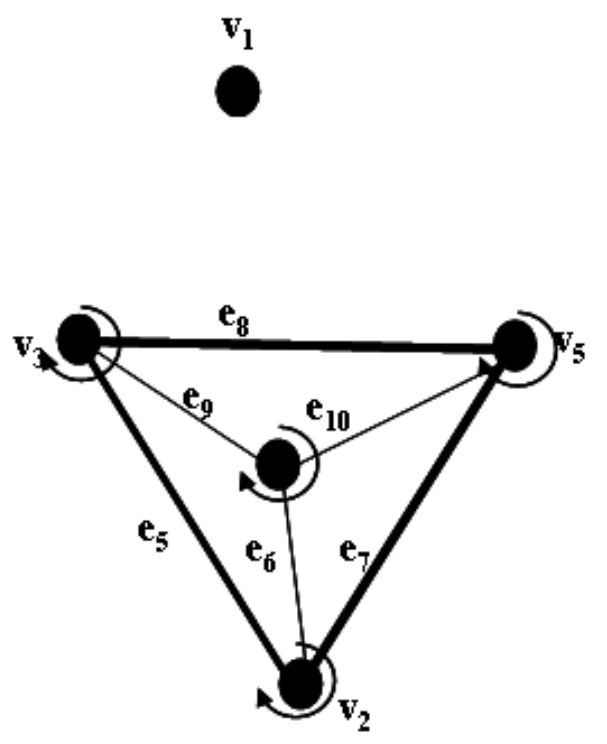

Рис. 2.6. Топологический рисунок плоского суграфа с несоприкасающимися циклами.

В данном случае имеется два пересекающихся ребра $e_9$ и $e_{10}$. Правило присоединения циклов нарушено. Данное построение не является гамильтоновым циклом, так как исключена вершина $v_1$.

## 2.2. Алгоритм построения гамильтонового цикла полного графа

**Инициализация.** Выделяем в полном графе G(E,V) систему изометрических циклов $C_\tau$. Записываем каждый цикл $c_i$ в виде подмножества рёбер $E_i$ и подмножества вершин $V_i$.

**Шаг 1. [Выбор].** Выбираем произвольный изометрический цикл.

**Шаг 2. [Поиск].** Осуществляем поиск соприкасающегося цикла.

**Шаг 3. [Суммирование].** Осуществляем кольцевое суммирование предыдущего цикла и соприкасающегося. Результат суммирования помещаем в подмножество ребер цикла. Добавляем новую вершину в подмножество вершин цикла.

**Шаг 4. [Проверка].** Если в подмножество вершин включены все вершины графа, то конец работы алгоритма. Иначе идти на шаг 2.

В качестве примера рассмотрим процесс построения гамильтононового цикла в графе $K_6$ (рис. 2.7).

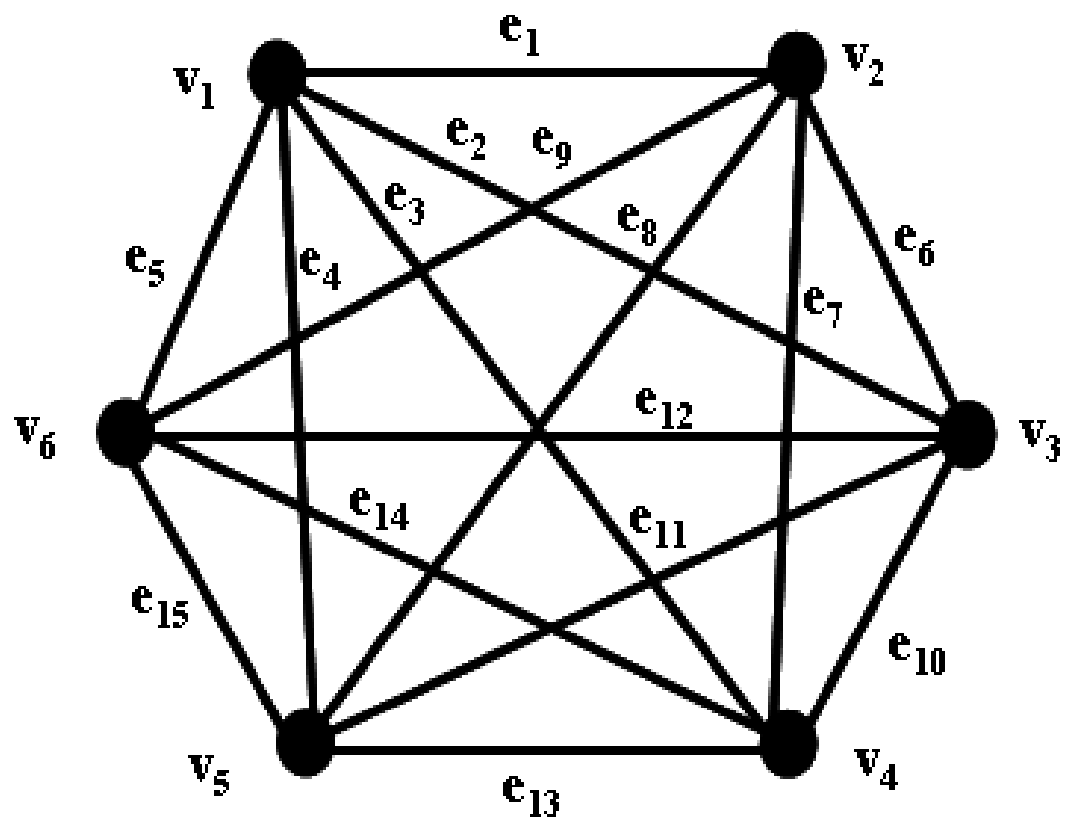

Рис. 2.7. Полный граф $K_6$.

Изометрические (треугольные) циклы графа:

$c_1 = \{e_1,e_2,e_6\} \leftrightarrow \{v_1,v_2,v_3\}$; $c_2 = \{e_1,e_3,e_7\} \leftrightarrow \{v_1,v_2,v_4\}$;
$c_3 = \{e_1,e_4,e_8\} \leftrightarrow \{v_1,v_2,v_5\}$; $c_4 = \{e_1,e_5,e_9\} \leftrightarrow \{v_1,v_2,v_6\}$;

$c_5 = \{e_2,e_3,e_{10}\} \leftrightarrow \{v_1,v_3,v_4\}$; $c_6 = \{e_2,e_4,e_{11}\} \leftrightarrow \{v_1,v_3,v_5\}$;
$c_7 = \{e_2,e_5,e_{12}\} \leftrightarrow \{v_1,v_3,v_6\}$; $c_8 = \{e_3,e_4,e_{13}\} \leftrightarrow \{v_1,v_4,v_5\}$;
$c_9 = \{e_3,e_5,e_{14}\} \leftrightarrow \{v_1,v_4,v_6\}$; $c_{10} = \{e_4,e_5,e_{15}\} \leftrightarrow \{v_1,v_5,v_6\}$;
$c_{11} = \{e_6,e_7,e_{10}\} \leftrightarrow \{v_2,v_3,v_4\}$; $c_{12} = \{e_6,e_8,e_{11}\} \leftrightarrow \{v_2,v_3,v_5\}$;
$c_{13} = \{e_6,e_9,e_{12}\} \leftrightarrow \{v_2,v_3,v_6\}$; $c_{14} = \{e_7,e_8,e_{13}\} \leftrightarrow \{v_2,v_4,v_5\}$;
$c_{15} = \{e_7,e_9,e_{14}\} \leftrightarrow \{v_2,v_4,v_6\}$; $c_{16} = \{e_8,e_9,e_{15}\} \leftrightarrow \{v_2,v_5,v_6\}$;
$c_{17} = \{e_{10},e_{11},e_{13}\} \leftrightarrow \{v_3,v_4,v_5\}$; $c_{18} = \{e_{10},e_{12},e_{14}\} \leftrightarrow \{v_3,v_4,v_6\}$;
$c_{19} = \{e_{11},e_{12},e_{15}\} \leftrightarrow \{v_3,v_5,v_6\}$; $c_{20} = \{e_{13},e_{14},e_{15}\} \leftrightarrow \{v_4,v_5,v_6\}$.

Выбираем цикл $c_1$. Ищем первый попавшийся соприкасающийся цикл. Это цикл $c_2$. Суммируем эти циклы:

$z_1 = c_1 \oplus c_2 = \{e_1,e_2,e_6\} \oplus \{e_1,e_3,e_7\} = \{e_2,e_3,e_6,e_7\} \leftrightarrow \{v_1,v_2,v_3,v_4\}$.

Ищем цикл соприкасающийся с $z_1$. Это цикл $c_8$. Суммируем циклы:

$z_2 = z_1 \oplus c_8 = \{e_2,e_3,e_6,e_7\} \oplus \{e_3,e_4,e_{13}\} = \{e_2,e_4,e_6,e_7,c_{13}\} \leftrightarrow \{v_1,v_2,v_3,v_4,v_5\}$.

Ищем цикл соприкасающийся с $z_2$. Это цикл $c_{15}$. Суммируем циклы:

$z_3 = z_2 \oplus c_{13} = \{e_2,e_4,e_6,e_7,c_{13}\} \oplus \{e_6,e_9,e_{12}\} = \{e_2,e_4,e_7,e_9,c_{12},c_{13}\} \leftrightarrow \{v_1,v_2,v_3,v_4,v_5,v_6\}$.

В подмножестве вершин цикла $z_3$ присутствуют все вершины. Конец работы алгоритма. Гамильтонов цикл состоит из следующего подмножества рёбер: $\{e_2,e_4,e_7,e_9,c_{12},c_{13}\}$.

## Глава 3. ЗАДАЧА КОММИВОЯЖЁРА

### 3.1. Гамильтонов цикл в задаче коммивояжёра

В задаче о коммивояжёре существенную роль играет вес ребра. Поэтому требуется найти не просто гамильтонов цикл, а гамильтонов цикл с минимальным суммарным весом. И здесь существенную роль играет способ построения такого цикла.

В традиционным подходе гамильтонов цикл рассматривается как циклический маршрут, проходящий по всем вершинам графа. Отсюда вытекает понятие оптимального маршрута – это гамильтонов цикл с минимальной аддитивной суммой весов рёбер.

$$\sum_{i=1}^{n} w_i \to \min, \tag{3.1}$$

здесь $w_i$ – вес ребра $e_i$.

Однако в такой конструкции не определён способ последовательности подключения ребер для формирования гамильтонового цикла в зависимости от веса ребра.

Для построения плоского топологического рисунка гамильтонового цикла воспользуемся принципом построения гамильтонового цикла как последовательности объединения соприкасающихся изометрических циклов.

*Гамильтоновым циклом* в графе называется простой цикл с количеством ребер равным количеству вершин графа, который характеризует обод плоского суграфа с нулевым значением кубического функционала Маклейна.

Известно, что кольцевая сумма независимых циклов [3,13] и обода в топологическом рисунке плоского суграфа есть пустое множество [15,16].

$$\sum_{i=1}^{n-2} c_i \oplus c_0 = \varnothing. \tag{3.2}$$

Отсюда

$$\sum_{i=1}^{n-2} c_i = c_0 \tag{3.3}$$

С другой стороны

$$\sum_{j=1}^{n} e_j = c_0 \tag{3.4}$$

В итоге

$$\sum_{i=1}^{n-2} c_i = \sum_{j=1}^{n} e_j = c_0 \tag{3.5}$$

**Определение 3.1.** Вес цикла $w_{c_i}$ есть аддитивная сумма рёбер цикла $c_i$.

С учетом веса ребер и веса циклов можно записать общее выражение

$$\sum_{i=1}^{n-2} c_{\min} = \min \qquad (3.6)$$

Данное выражение можно трактовать как последовательное соединение треугольных соприкасающихся циклов

$$(...(((w_{c_1} \oplus w_{c_2})_{\min} \oplus w_{c_3})_{\min} \oplus w_{c_4})_{\min} \oplus ... \oplus w_{c_{n-2}})_{\min} = \min \qquad (3.7)$$

Таким образом, построение осуществляется последовательным поиском и кольцевым сложением треугольных циклов с минимальным суммарным весом. Сначала образуется простой цикл с минимальным весом (стоимостью) длиной четыре, путём кольцевого суммирования соприкасающихся треугольных циклов. Далее, путём присоединения соприкасающегося цикла, образуется простой цикл минимальной стоимости длиной пять и так далее, до присоединения последнего соприкасающегося цикла. Всего количество треугольных циклов для построения гамильтонового цикла равно $n$-2.

Перейдём к рассмотрению данного метода на конкретных примерах.

## 3.2. Задача коммивояжёра для полного графа $K_4$

Для решения задачи коммивояжёра воспользуемся алгоритмом описанным в главе 2 и динамическим принципом Беллмана [22]. Рассмотрим особый случай выделения минимального гамильтонового цикла в взвешенном графе $K_4$, описанный в книге [4].

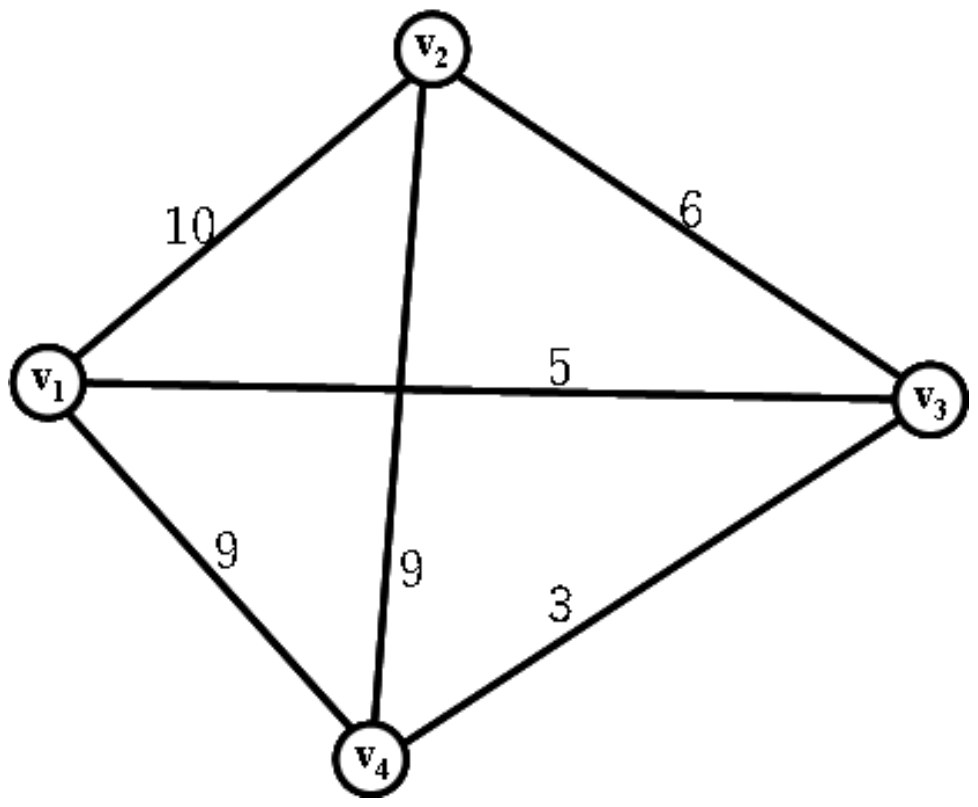

Рис. 3.1. Граф с взвешенными ребрами.

Зададим веса ребер:

ребро $e_1 \leftrightarrow (v_1,v_2) = 10$;
ребро $e_2 \leftrightarrow (v_1,v_3) = 5$;
ребро $e_3 \leftrightarrow (v_1,v_4) = 9$;
ребро $e_4 \leftrightarrow (v_2,v_3) = 6$;
ребро $e_5 \leftrightarrow (v_2,v_4) = 9$;
ребро $e_6 \leftrightarrow (v_3,v_4) = 3$.

Определим веса изометрических циклов:

цикл $c_1 = \{e_1,e_2,e_4\} \leftrightarrow \{v_1,v_2,v_3\} = 10+5+6 = 21$;
цикл $c_2 = \{e_1,e_3,e_5\} \leftrightarrow \{v_1,v_2,v_4\} = 10+9+9 = 28$;
цикл $c_3 = \{e_2,e_3,e_6\} \leftrightarrow \{v_1,v_3,v_4\} = 9+5+3 = 17$;
цикл $c_4 = \{e_4,e_5,e_6\} \leftrightarrow \{v_2,v_3,v_4\} = 9+3+6 = 18$.

$c_1 \oplus c_2 = \{e_2,e_3,e_4,e_5\} \leftrightarrow \{v_1,v_2,v_3,v_4\} = 5+6+9+9 = 29;$
$c_1 \oplus c_3 = \{e_1,e_3,e_4,e_6\} \leftrightarrow \{v_1,v_2,v_3,v_4\} = 10+9+6+3 = 28;$
$c_1 \oplus c_4 = \{e_1,e_2,e_5,e_6\} \leftrightarrow \{v_1,v_2,v_3,v_4\} = 10+5+9+3 = 27;$
$c_2 \oplus c_3 = \{e_1,e_2,e_5,e_6\} \leftrightarrow \{v_1,v_2,v_3,v_4\} = 10+5+9+3 = 27;$
$c_2 \oplus c_4 = \{e_1,e_3,e_4,e_6\} \leftrightarrow \{v_1,v_2,v_3,v_4\} = 10+9+6+3 = 28;$
$c_3 \oplus c_4 = \{e_2,e_3,e_4,e_5\} \leftrightarrow \{v_1,v_2,v_3,v_4\} = 5+6+9+9 = 29,$

Кольцевая сумма двух циклов определяет все возможные ориентированные циклы длиной четыре (рис. 3.2).

a) $<v_1,v_2,v_3.v_4> \leftrightarrow c_1 \oplus c_3 = \{e_1,e_3,e_4,e_6\} = 10+9+6+3 = 28;$
b) $<v_1,v_2,v_4.v_3> \leftrightarrow c_2 \oplus c_3 = \{e_1,e_2,e_5,e_6\} = 10+5+9+3 = 27;$
c) $<v_1,v_3,v_2.v_4> \leftrightarrow c_1 \oplus c_2 = \{e_2,e_3,e_4,e_5\} = 5+6+9+9 = 29;$
d) $<v_1,v_4,v_3.v_2> \leftrightarrow c_2 \oplus c_4 = \{e_1,e_3,e_4,e_6\} = 10+9+6+3 = 28;$
e) $<v_1,v_4,v_2.v_3> \leftrightarrow c_3 \oplus c_4 = \{e_2,e_3,e_4,e_5\} = 5+6+9+9 = 29;$
f) $<v_1,v_3,v_4.v_2> \leftrightarrow c_1 \oplus c_4 = \{e_1,e_2,e_5,e_6\} = 10+5+9+3 = 27.$

**Определение 2.2.** *Дубль-цикл* – это простой цикл, который допускает по крайней мере два различных нетривиальных разложения в сумму изометрических циклов.

В нашем случае может быть образован следующий дубль-цикл:

$\{a,b,c\} \oplus \{c,d,e\} = \{b,g,e\} \oplus \{a,d,g\} = \{a,b,d,e\},$

или можно организовать следующие дубль-циклы:

$\{a,b,c\} \oplus \{b,g,e\} = \{c,d,e\} \oplus \{a,d,g\} = \{a,c,e,g\},$
$\{a,b,c\} \oplus \{c,d,e\} \oplus \{b,g,e\} = \{a,d,g\}.$

По сути, множество дубль-циклов – это подмножество простых циклов. Поэтому дубль-циклы обладают всеми свойствами простых циклов. Изометрические циклы, простые циклы и дубль-циклы являются частными случаями квазициклов.

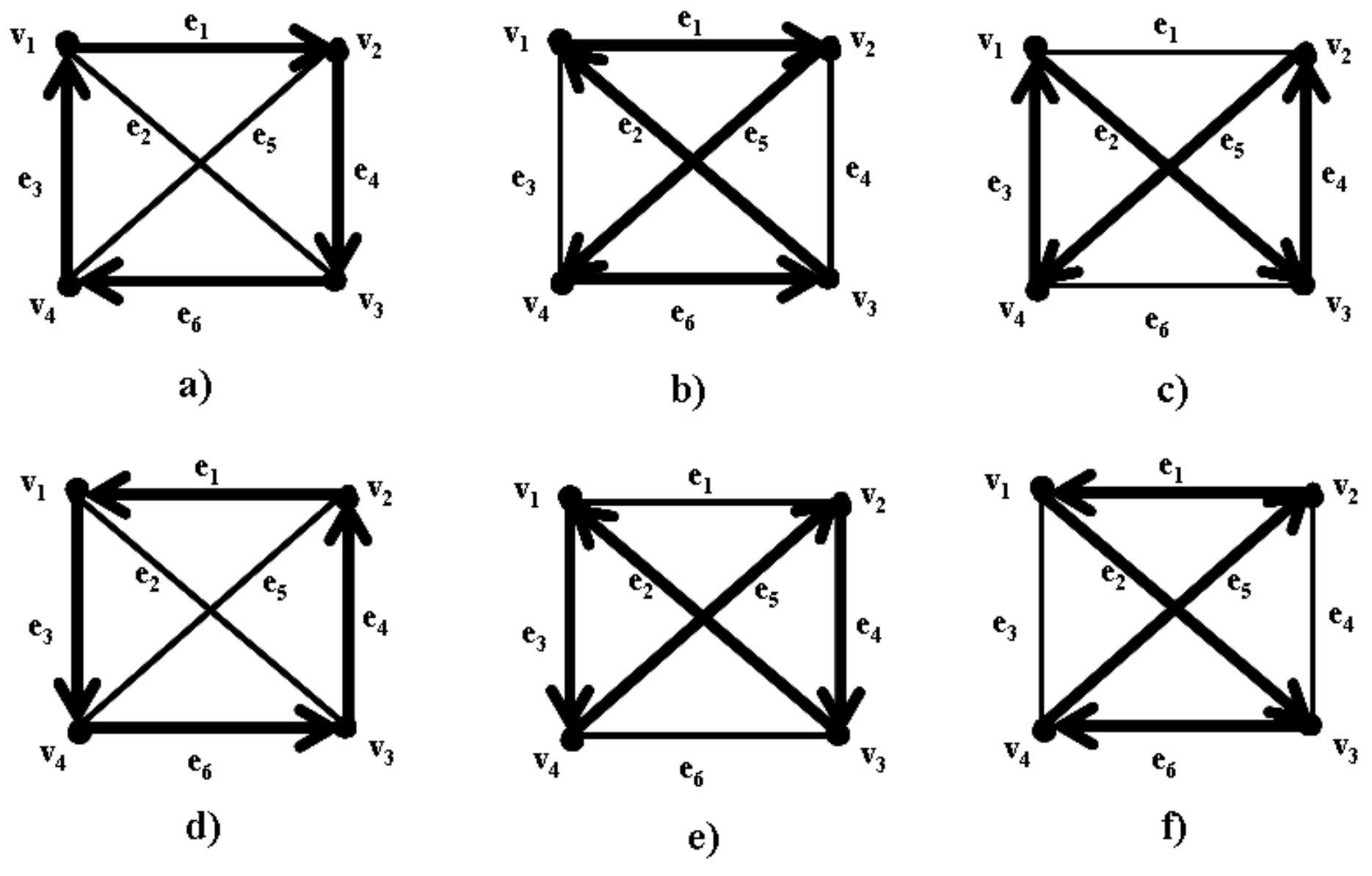

Рис. 3.2. Ориентированные циклы.

Следующие циклы являются дубль-циклами:
$c_1 \oplus c_3 = c_2 \oplus c_4$;   $c_2 \oplus c_3 = c_1 \oplus c_4$;   $c_1 \oplus c_2 = c_3 \oplus c_4$.

Таким образом, рассматриваемое количество циклов длиной четыре можно сократить до трех.

Если рассматривать гамильтоновы циклы с минимальным весом, то можно взять только один цикл с минимальным весом. В нашем случае – это цикл b) или f).

$c_1 \oplus c_4 = \{e_1,e_2,e_5,e_6\} \leftrightarrow \{v_1,v_2,v_3,v_4\} = 10+5+9+3 = 27$.

### 3.3. Задача коммивояжёра для полного графа $K_5$

Рассмотрим следующую задачу коммивояжёра на графе $K_5$ (рис. 3.3).

Зададим веса ребер:

ребро $e_1 \leftrightarrow (v_1,v_2) = 6$;
ребро $e_2 \leftrightarrow (v_1,v_3) = 10$;
ребро $e_3 \leftrightarrow (v_1,v_4) = 5$;
ребро $e_4 \leftrightarrow (v_1,v_5) = 11$;
ребро $e_5 \leftrightarrow (v_2,v_3) = 10$;
ребро $e_6 \leftrightarrow (v_2,v_4) = 9$;
ребро $e_7 \leftrightarrow (v_2,v_5) = 7$;
ребро $e_8 \leftrightarrow (v_3,v_4) = 9$;
ребро $e_9 \leftrightarrow (v_3,v_5) = 8$;
ребро $e_{10} \leftrightarrow (v_4,v_5) = 11$.

Определим веса изометрических циклов:

цикл $c_1 = \{e_1,e_2,e_5\} \leftrightarrow \{v_1,v_2,v_3\} = 10+5+6 = 21$;
цикл $c_2 = \{e_1,e_3,e_6\} \leftrightarrow \{v_1,v_2,v_4\} = 9+5+6 = 20$;
цикл $c_3 = \{e_1,e_4,e_7\} \leftrightarrow \{v_1,v_2,v_5\} = 11+7+6 = 24$;
цикл $c_4 = \{e_2,e_3,e_8\} \leftrightarrow \{v_1,v_3,v_4\} = 10+5+9 = 24$;
цикл $c_5 = \{e_2,e_4,e_9\} \leftrightarrow \{v_1,v_3,v_5\} = 10+11+8 = 29$;
цикл $c_6 = \{e_3,e_4,e_{10}\} \leftrightarrow \{v_1,v_4,v_5\} = 11+5+11 = 27$;
цикл $c_7 = \{e_5,e_6,e_8\} \leftrightarrow \{v_2,v_3,v_4\} = 10+9+9 = 28$;
цикл $c_8 = \{e_5,e_7,e_9\} \leftrightarrow \{v_2,v_3,v_5\} = 10+7+9 = 26$;
цикл $c_9 = \{e_6,e_7,e_{10}\} \leftrightarrow \{v_2,v_4,v_5\} = 9+7+11 = 27$;
цикл $c_{10} = \{e_8,e_9,e_{10}\} \leftrightarrow \{v_3,v_4,v_5\} = 11+9+8 = 28$.

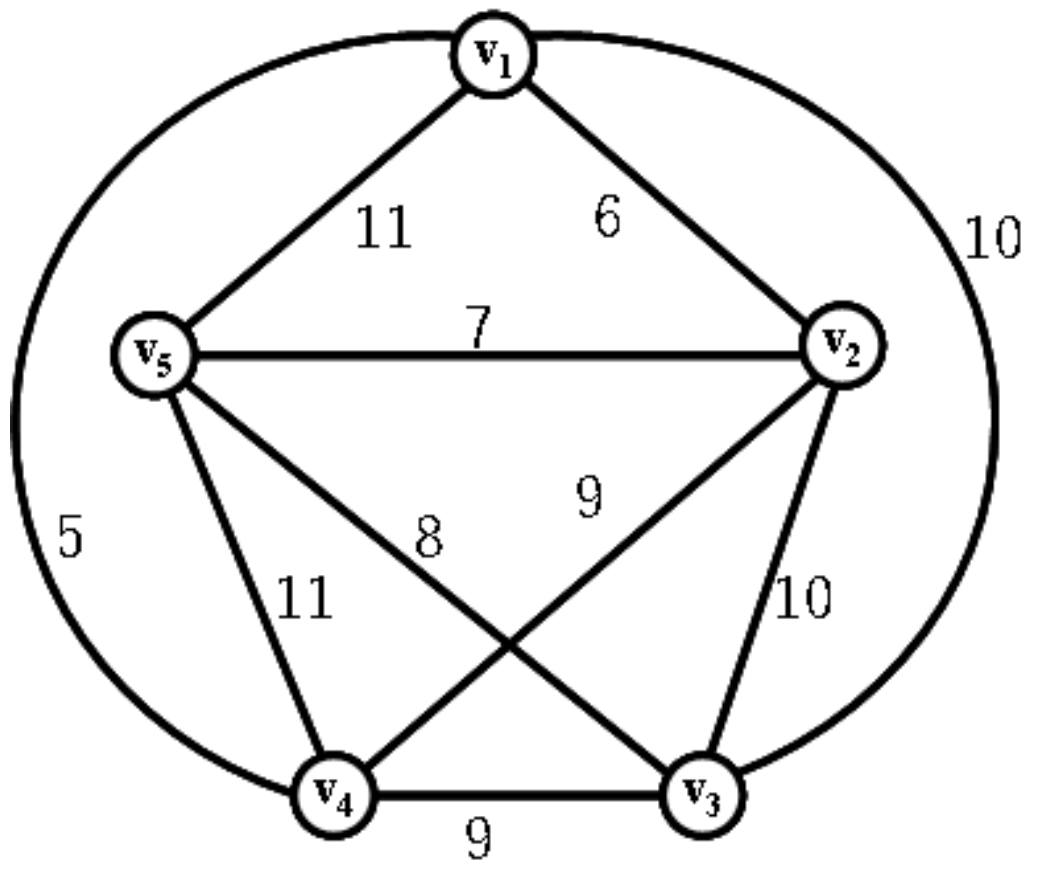

Рис. 3.3. Граф $K_5$ с весами ребер.

Формируем сумму двух циклов для построения множества соприкасающихся циклов длиной четыре:

$c_1 \oplus c_2 = \{e_1,e_2,e_5\} \oplus \{e_1,e_3,e_6\} \leftrightarrow \{e_2,e_3,e_5,e_6\} \leftrightarrow \{v_1,v_2,v_3,v_4\} = 10+5+10+9 = 34;$
$c_1 \oplus c_3 = \{e_1,e_2,e_5\} \oplus \{e_1,e_4,e_7\} \leftrightarrow \{e_2,e_4,e_5,e_7\} \leftrightarrow \{v_1,v_2,v_3,v_5\} = 10+11+10+7 = 38;$
$c_1 \oplus c_4 = \{e_1,e_2,e_5\} \oplus \{e_2,e_3,e_8\} \leftrightarrow \{e_1,e_3,e_5,e_8\} \leftrightarrow \{v_1,v_2,v_3,v_4\} = 6+5+10+9 = 30;$
$c_1 \oplus c_5 = \{e_1,e_2,e_5\} \oplus \{e_2,e_4,e_9\} \leftrightarrow \{e_1,e_4,e_5,e_9\} \leftrightarrow \{v_1,v_2,v_3,v_5\} = 6+11+10+8 = 35;$
$c_1 \oplus c_6 = \{e_1,e_2,e_5\} \cap \{e_3,e_4,e_{10}\} \leftrightarrow \varnothing;$
$c_1 \oplus c_7 = \{e_1,e_2,e_5\} \oplus \{e_5,e_6,e_8\} \leftrightarrow \{e_1,e_2,e_6,e_8\} \leftrightarrow \{v_1,v_2,v_3,v_4\} = 6+10+9+9 = 34;$
$c_1 \oplus c_8 = \{e_1,e_2,e_5\} \oplus \{e_5,e_7,e_9\} \leftrightarrow \{e_1,e_2,e_7,e_9\} \leftrightarrow \{v_1,v_2,v_3,v_5\} = 6+10+7+8 = 31;$
$c_1 \oplus c_9 = \{e_1,e_2,e_5\} \cap \{e_6,e_7,e_{10}\} \leftrightarrow \varnothing;$
$c_1 \oplus c_{10} = \{e_1,e_2,e_5\} \cap \{e_8,e_9,e_{10}\} \leftrightarrow \varnothing;$
$c_2 \oplus c_3 = \{e_1,e_3,e_6\} \oplus \{e_1,e_4,e_7\} \leftrightarrow \{e_3,e_4,e_6,e_7\} \leftrightarrow \{v_1,v_2,v_4,v_5\} = 5+11+9+7 = 32;$
$c_2 \oplus c_4 = \{e_1,e_3,e_6\} \oplus \{e_2,e_3,e_8\} \leftrightarrow \{e_1,e_2,e_6,e_8\} \leftrightarrow \{v_1,v_2,v_3,v_4\} = 6+10+9+9 = 34;$
$c_2 \oplus c_5 = \{e_1,e_3,e_6\} \cap \{e_2,e_4,e_9\} \leftrightarrow \varnothing,$
$c_2 \oplus c_6 = \{e_1,e_3,e_6\} \oplus \{e_3,e_4,e_{10}\} \leftrightarrow \{e_1,e_4,e_6,e_{10}\} \leftrightarrow \{v_1,v_2,v_4,v_5\} = 6+11+9+11 = 37;$
$c_2 \oplus c_7 = \{e_1,e_3,e_6\} \oplus \{e_5,e_6,e_8\} \leftrightarrow \{e_1,e_3,e_5,e_8\} \leftrightarrow \{v_1,v_2,v_3,v_4\} = 6+5+10+9 = 30;$
$c_2 \oplus c_8 = \{e_1,e_3,e_6\} \cap \{e_5,e_7,e_9\} \leftrightarrow \varnothing;$
$c_2 \oplus c_9 = \{e_1,e_3,e_6\} \oplus \{e_6,e_7,e_{10}\} \leftrightarrow \{e_1,e_3,e_7,e_{10}\} \leftrightarrow \{v_1,v_2,v_4,v_5\} = 6+5+7+11 = 29;$
$c_2 \oplus c_{10} = \{e_1,e_3,e_6\} \cap \{e_8,e_9,e_{10}\} \leftrightarrow \varnothing;$
$c_3 \oplus c_4 = \{e_1,e_4,e_7\} \cap \{e_2,e_3,e_8\} \leftrightarrow \varnothing;$
$c_3 \oplus c_5 = \{e_1,e_4,e_7\} \oplus \{e_2,e_4,e_9\} \leftrightarrow \{e_1,e_2,e_7,e_9\} \leftrightarrow \{v_1,v_2,v_3,v_5\} = 6+10+7+8 = 31;$
$c_3 \oplus c_6 = \{e_1,e_4,e_7\} \oplus \{e_3,e_4,e_{10}\} \leftrightarrow \{e_1,e_3,e_7,e_{10}\} \leftrightarrow \{v_1,v_2,v_4,v_5\} = 6+5+7+11 = 29;$
$c_3 \oplus c_7 = \{e_1,e_4,e_7\} \cap \{e_5,e_6,e_8\} \leftrightarrow \varnothing;$
$c_3 \oplus c_8 = \{e_1,e_4,e_7\} \oplus \{e_5,e_7,e_9\} \leftrightarrow \{e_1,e_4,e_5,e_9\} \leftrightarrow \{v_1,v_2,v_3,v_5\} = 6+11+10+8 = 35;$
$c_3 \oplus c_9 = \{e_1,e_4,e_7\} \oplus \{e_6,e_7,e_{10}\} \leftrightarrow \{e_1,e_4,e_6,e_{10}\} \leftrightarrow \{v_1,v_2,v_4,v_5\} = 6+11+9+11 = 37;$
$c_3 \oplus c_{10} = \{e_1,e_4,e_7\} \cap \{e_8,e_9,e_{10}\} \leftrightarrow \varnothing;$
$c_4 \oplus c_5 = \{e_2,e_3,e_8\} \oplus \{e_2,e_4,e_9\} \leftrightarrow \{e_3,e_4,e_8,e_9\} \leftrightarrow \{v_1,v_3,v_4,v_5\} = 5+11+9+8 = 33;$
$c_4 \oplus c_6 = \{e_2,e_3,e_8\} \oplus \{e_3,e_4,e_{10}\} \leftrightarrow \{e_2,e_4,e_8,e_{10}\} \leftrightarrow \{v_1,v_3,v_4,v_5\} = 10+11+9+11 = 41;$
$c_4 \oplus c_7 = \{e_2,e_3,e_8\} \oplus \{e_5,e_6,e_8\} \leftrightarrow \{e_2,e_3,e_5,e_6\} \leftrightarrow \{v_1,v_2,v_3,v_4\} = 10+5+10+9 = 34;$
$c_4 \oplus c_8 = \{e_2,e_3,e_8\} \cap \{e_5,e_7,e_9\} \leftrightarrow \varnothing;$
$c_4 \oplus c_9 = \{e_2,e_3,e_8\} \cap \{e_6,e_7,e_{10}\} \leftrightarrow \varnothing;$
$c_4 \oplus c_{10} = \{e_2,e_3,e_8\} \oplus \{e_8,e_9,e_{10}\} \leftrightarrow \{e_2,e_3,e_9,e_{10}\} \leftrightarrow \{v_1,v_3,v_4,v_5\} = 10+5+8+11 = 34;$
$c_5 \oplus c_6 = \{e_2,e_4,e_9\} \oplus \{e_3,e_4,e_{10}\} \leftrightarrow \{e_2,e_3,e_9,e_{10}\} \leftrightarrow \{v_1,v_3,v_4,v_5\} = 10+5+8+11 = 34;$
$c_5 \oplus c_7 = \{e_2,e_4,e_9\} \cap \{e_5,e_6,e_8\} \leftrightarrow \varnothing;$
$c_5 \oplus c_8 = \{e_2,e_4,e_9\} \oplus \{e_5,e_7,e_9\} \leftrightarrow \{e_2,e_4,e_5,e_7\} \leftrightarrow \{v_1,v_2,v_3,v_5\} = 10+11+10+7 = 38;$
$c_5 \oplus c_9 = \{e_2,e_4,e_9\} \cap \{e_6,e_7,e_{10}\} \leftrightarrow \varnothing;$
$c_5 \oplus c_{10} = \{e_2,e_4,e_9\} \oplus \{e_8,e_9,e_{10}\} \leftrightarrow \{e_2,e_4,e_8,e_{10}\} \leftrightarrow \{v_1,v_3,v_4,v_5\} = 10+11+9+11 = 41;$
$c_6 \oplus c_7 = \{e_3,e_4,e_{10}\} \cap \{e_5,e_6,e_8\} \leftrightarrow \varnothing;$
$c_6 \oplus c_8 = \{e_3,e_4,e_{10}\} \cap \{e_5,e_7,e_9\} \leftrightarrow \varnothing;$
$c_6 \oplus c_9 = \{e_3,e_4,e_{10}\} \oplus \{e_6,e_7,e_{10}\} \leftrightarrow \{e_3,e_4,e_6,e_7\} \leftrightarrow \{v_1,v_2,v_4,v_5\} = 5+11+9+7 = 32;$
$c_6 \oplus c_{10} = \{e_3,e_4,e_{10}\} \oplus \{e_8,e_9,e_{10}\} \leftrightarrow \{e_3,e_4,e_8,e_9\} \leftrightarrow \{v_1,v_3,v_4,v_5\} = 5+11+9+8 = 33;$
$c_7 \oplus c_8 = \{e_5,e_6,e_8\} \oplus \{e_5,e_7,e_9\} \leftrightarrow \{e_6,e_7,e_8,e_9\} \leftrightarrow \{v_2,v_3,v_4,v_5\} = 9+7+9+8 = 33;$
$c_7 \oplus c_9 = \{e_5,e_6,e_8\} \oplus \{e_6,e_7,e_{10}\} \leftrightarrow \{e_5,e_7,e_8,e_{10}\} \leftrightarrow \{v_2,v_3,v_4,v_5\} = 10+7+9+11 = 37;$
$c_7 \oplus c_{10} = \{e_5,e_6,e_8\} \oplus \{e_8,e_9,e_{10}\} \leftrightarrow \{e_5,e_6,e_9,e_{10}\} \leftrightarrow \{v_2,v_3,v_4,v_5\} = 10+9+8+11 = 38;$
$c_8 \oplus c_9 = \{e_5,e_7,e_9\} \oplus \{e_6,e_7,e_{10}\} \leftrightarrow \{e_5,e_6,e_9,e_{10}\} \leftrightarrow \{v_2,v_3,v_4,v_5\} = 10+9+8+11 = 38;$
$c_8 \oplus c_{10} = \{e_5,e_7,e_9\} \oplus \{e_8,e_9,e_{10}\} \leftrightarrow \{e_5,e_7,e_8,e_{10}\} \leftrightarrow \{v_2,v_3,v_4,v_5\} = 10+7+9+11 = 37;$
$c_9 \oplus c_{10} = \{e_6,e_7,e_{10}\} \oplus \{e_8,e_9,e_{10}\} \leftrightarrow \{e_6,e_7,e_8,e_9\} \leftrightarrow \{v_2,v_3,v_4,v_5\} = 9+7+9+8 = 33.$

Количество кольцевых сумм двух циклов можно рассчитать по формуле:

$$k_{2c} = \left( \frac{n(n-1)(n-2)}{6} \right)^2 = \frac{n^2(n-1)^2(n-2)^2}{36} \approx \frac{n^6}{36}. \qquad (3.8)$$

Однако количество соприкасающихся циклов явно меньше, и может быть определено как утроенное сочетание вершин по четыре, согласно положению о дубль-циклах:

$$k_4 = 3\frac{n(n-1)(n-2)(n-3)}{4!} \approx \frac{n^4}{8}. \quad (3.9)$$

Выделим циклы длиной четыре относительно подмножества четвёрок вершин:

$z_1 = \{e_2,e_3,e_5,e_6\} \leftrightarrow \{v_1,v_2,v_3,v_4\} = 10+5+10+9 = 34$;
$z_2 = \{e_1,e_3,e_5,e_8\} \leftrightarrow \{v_1,v_2,v_3,v_4\} = 6+5+10+9 = 30$;
$z_3 = \{e_1,e_2,e_6,e_8\} \leftrightarrow \{v_1,v_2,v_3,v_4\} = 6+10+9+9 = 34$;

$z_4 = \{e_3,e_4,e_8,e_9\} \leftrightarrow \{v_1,v_3,v_4,v_5\} = 5+11+9+8 = 33$;
$z_5 = \{e_2,e_4,e_8,e_{10}\} \leftrightarrow \{v_1,v_3,v_4,v_5\} = 10+11+9+11 = 41$;
$z_6 = \{e_2,e_3,e_9,e_{10}\} \leftrightarrow \{v_1,v_3,v_4,v_5\} = 10+5+8+11 = 34$;

$z_7 = \{e_2,e_4,e_5,e_7\} \leftrightarrow \{v_1,v_2,v_3,v_5\} = 10+11+10+7 = 38$;
$z_8 = \{e_1,e_4,e_5,e_9\} \leftrightarrow \{v_1,v_2,v_3,v_5\} = 6+11+10+8 = 35$;
$z_9 = \{e_1,e_2,e_7,e_9\} \leftrightarrow \{v_1,v_2,v_3,v_5\} = 6+10+7+8 = 31$;

$z_{10} = \{e_3,e_4,e_6,e_7\} \leftrightarrow \{v_1,v_2,v_4,v_5\} = 5+11+9+7 = 32$;
$z_{11} = \{e_1,e_4,e_6,e_{10}\} \leftrightarrow \{v_1,v_2,v_4,v_5\} = 6+11+9+11 = 37$;
$z_{12} = \{e_1,e_3,e_7,e_{10}\} \leftrightarrow \{v_1,v_2,v_4,v_5\} = 6+5+7+11 = 29$;

$z_{13} = \{e_6,e_7,e_8,e_9\} \leftrightarrow \{v_2,v_3,v_4,v_5\} = 9+7+9+8 = 33$;
$z_{14} = \{e_5,e_7,e_8,e_{10}\} \leftrightarrow \{v_2,v_3,v_4,v_5\} = 10+7+9+11 = 37$;
$z_{15} = \{e_5,e_6,e_9,e_{10}\} \leftrightarrow \{v_2,v_3,v_4,v_5\} = 10+9+8+11 = 38$;

Из каждой четвертки вершин выбираем циклы с минимальным весом:

$z_2 = \{e_1,e_3,e_5,e_8\} \leftrightarrow \{v_1,v_2,v_3,v_4\} = 6+5+10+9 = 30$;
$z_4 = \{e_3,e_4,e_8,e_9\} \leftrightarrow \{v_1,v_3,v_4,v_5\} = 5+11+9+8 = 33$;
$z_8 = \{e_1,e_2,e_7,e_9\} \leftrightarrow \{v_1,v_2,v_3,v_5\} = 6+10+7+8 = 31$;
$z_{12} = \{e_1,e_3,e_7,e_{10}\} \leftrightarrow \{v_1,v_2,v_4,v_5\} = 6+5+7+11 = 29$;
$z_{13} = \{e_6,e_7,e_8,e_9\} \leftrightarrow \{v_2,v_3,v_4,v_5\} = 9+7+9+8 = 33$;

Выбираем цикл длиной четыре с минимальным весом:

$z_{12} = \{e_1,e_3,e_7,e_{10}\} \leftrightarrow \{v_1,v_2,v_4,v_5\} = 6+5+7+11 = 29$;

Для построения гамильтонового цикла, состоящего из пяти ребер, необходимо к циклу $z_{12}$ прибавить соприкасающийся изометрический цикл, имеющий в своём составе вершину $v_3$.

цикл $c_4 = \{e_2,e_3,e_8\} \leftrightarrow \{v_1,v_3,v_4\} = 10+5+9 = 24$;
цикл $c_5 = \{e_2,e_4,e_9\} \leftrightarrow \{v_1,v_3,v_5\} = 10+11+8 = 29$;
цикл $c_7 = \{e_5,e_6,e_8\} \leftrightarrow \{v_2,v_3,v_4\} = 10+9+9 = 28$;
цикл $c_8 = \{e_5,e_7,e_9\} \leftrightarrow \{v_2,v_3,v_5\} = 10+7+9 = 26$;
цикл $c_{10} = \{e_8,e_9,e_{10}\} \leftrightarrow \{v_3,v_4,v_5\} = 11+9+8 = 28$.

Удалим из подмножества те циклы, которые не имеют совпадающих рёбер с циклом $z_9$

цикл $c_4 = \{e_2,e_3,e_8\} \leftrightarrow \{v_1,v_3,v_4\} = 10+5+9 = 24$;

цикл $c_8 = \{e_5,e_7,e_9\} \leftrightarrow \{v_2,v_3,v_5\} = 10+7+9 = 26$;
цикл $c_{10} = \{e_8,e_9,e_{10}\} \leftrightarrow \{v_3,v_4,v_5\} = 11+9+8 = 28$.

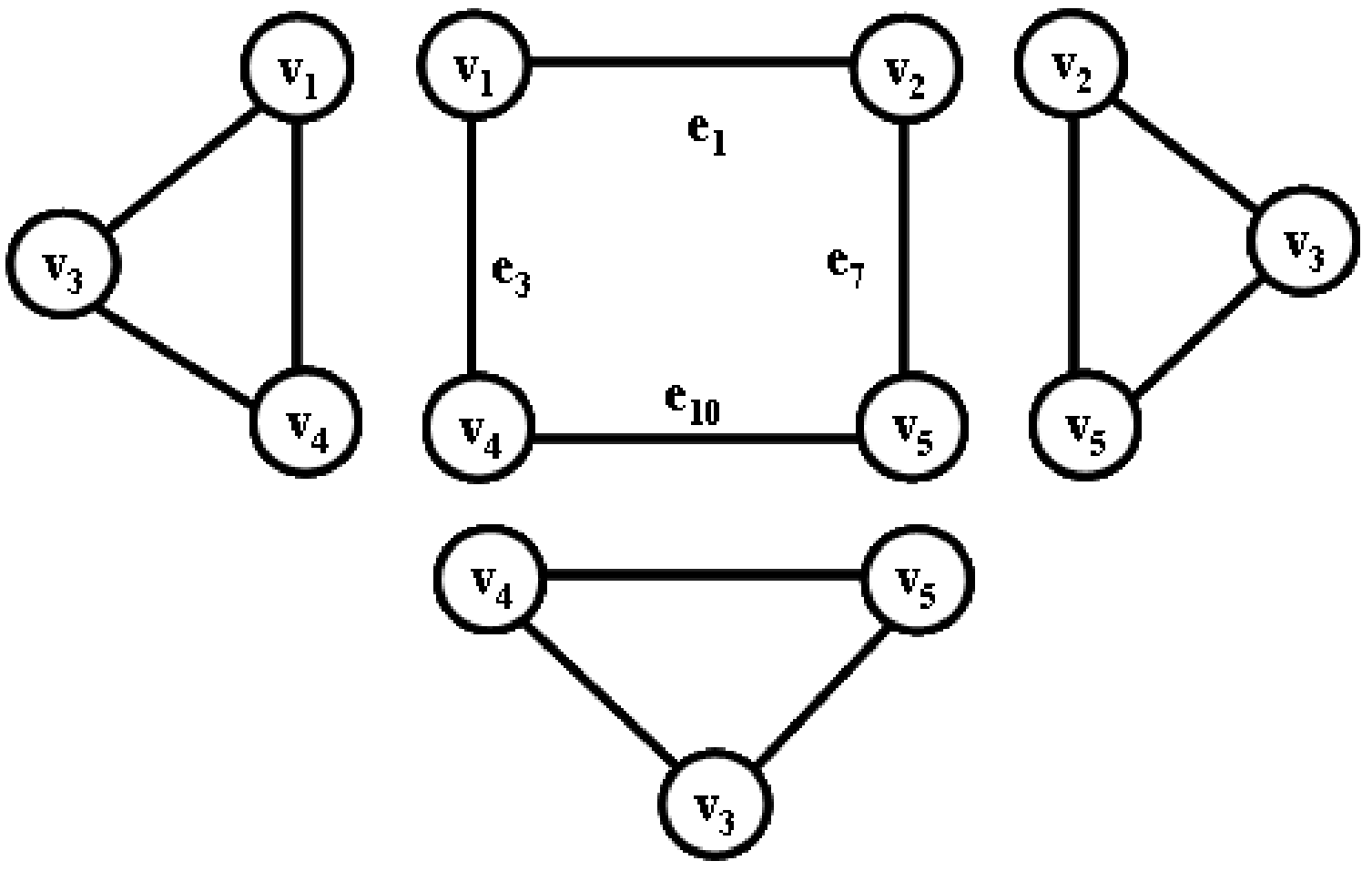

Рис. 3.4. Процесс подключения изометрических циклов.

цикл $c_4 \oplus z_9 = \{e_2,e_3,e_8\} \oplus \{e_1,e_3,e_7,e_{10}\} \leftrightarrow \{e_1,e_2,e_7,e_8,e_{10}\} = 10+9+6+7+11 = 43$;
цикл $c_8 \oplus z_9 = \{e_5,e_7,e_9\} \oplus \{e_1,e_3,e_7,e_{10}\} \leftrightarrow \{e_1,e_3,e_5,e_9,e_{10}\} = 10+9+6+5+11 = 41$;
цикл $c_{10} \oplus z_9 = \{e_8,e_9,e_{10}\} \oplus \{e_1,e_3,e_7,e_{10}\} \leftrightarrow \{e_1,e_3,e_7,e_9,e_{10}\} = 9+8+6+5+7 = 35$.

Определяем минимальный гамильтонов цикл:

$c_{min} = \{e_1,e_3,e_7,e_9,e_{10}\} \leftrightarrow \{v_1,v_2,v_3,v_4,v_5\} = 9+8+6+5+7 = 35$.

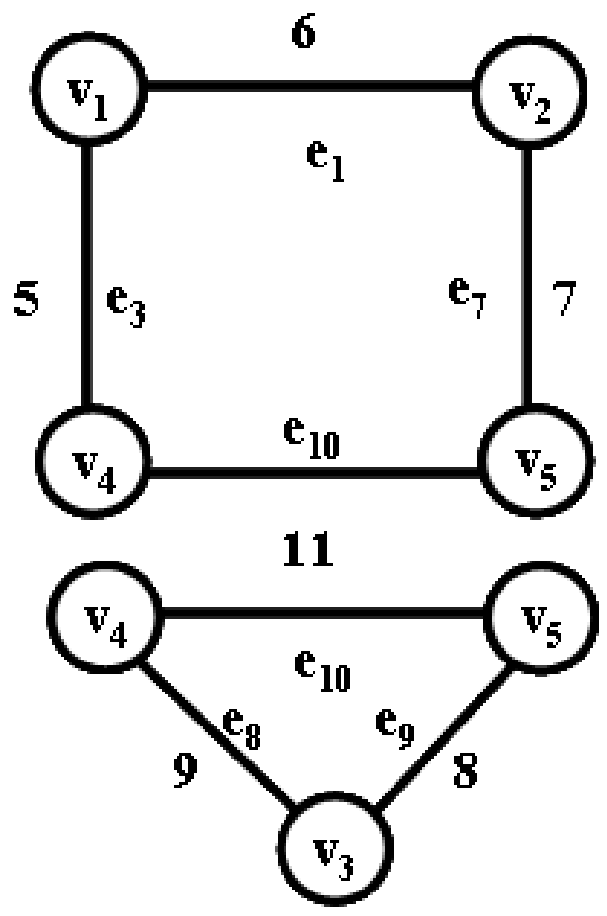

Рис. 3.5. Построение минимального гамильтонового цикла.

### 3.4. Задача коммивояжёра для полного графа $K_5$ и различные веса ребер

Проведем эксперимент. Увеличим стоимость ребра $e_{10}$. Предположим, что увеличенная стоимость не позволит построить гамильтонов цикл минимальной стоимости.

ребро $e_1 \leftrightarrow (v_1,v_2) = 6$;
ребро $e_2 \leftrightarrow (v_1,v_3) = 10$;
ребро $e_3 \leftrightarrow (v_1,v_4) = 5$;
ребро $e_4 \leftrightarrow (v_1,v_5) = 11$;
ребро $e_5 \leftrightarrow (v_2,v_3) = 10$;
ребро $e_6 \leftrightarrow (v_2,v_4) = 9$;
ребро $e_7 \leftrightarrow (v_2,v_5) = 7$;
ребро $e_8 \leftrightarrow (v_3,v_4) = 9$;

ребро $e_9 \leftrightarrow (v_3,v_5) = 8$;
ребро $e_{10} \leftrightarrow (v_4,v_5) = 40$.

Определим стоимость каждого цикла:

цикл $c_1 = \{e_1,e_2,e_5\} \leftrightarrow \{v_1,v_2,v_3\} = 10+5+6 = 21$;
цикл $c_2 = \{e_1,e_3,e_6\} \leftrightarrow \{v_1,v_2,v_4\} = 9+5+6 = 20$;
цикл $c_3 = \{e_1,e_4,e_7\} \leftrightarrow \{v_1,v_2,v_5\} = 11+7+6 = 24$;
цикл $c_4 = \{e_2,e_3,e_8\} \leftrightarrow \{v_1,v_3,v_4\} = 10+5+9 = 24$;
цикл $c_5 = \{e_2,e_4,e_9\} \leftrightarrow \{v_1,v_3,v_5\} = 10+11+8 = 29$;
цикл $c_6 = \{e_3,e_4,e_{10}\} \leftrightarrow \{v_1,v_4,v_5\} = 11+5+40 = 56$;
цикл $c_7 = \{e_5,e_6,e_8\} \leftrightarrow \{v_2,v_3,v_4\} = 10+9+9 = 28$;
цикл $c_8 = \{e_5,e_7,e_9\} \leftrightarrow \{v_2,v_3,v_5\} = 10+7+8 = 25$;
цикл $c_9 = \{e_6,e_7,e_{10}\} \leftrightarrow \{v_2,v_4,v_5\} = 9+7+40 = 56$;
цикл $c_{10} = \{e_8,e_9,e_{10}\} \leftrightarrow \{v_3,v_4,v_5\} = 40+9+8 = 57$.

Построим циклы длиной четыре:

$z_1 = \{e_2,e_3,e_5,e_6\} \leftrightarrow \{v_1,v_2,v_3,v_4\} = 10+5+10+9 = 34$;
$z_2 = \{e_2,e_4,e_5,e_7\} \leftrightarrow \{v_1,v_2,v_3,v_5\} = 10+11+10+7 = 38$;
$z_3 = \{e_1,e_3,e_5,e_8\} \leftrightarrow \{v_1,v_2,v_3,v_4\} = 6+5+10+9 = 30$;
$z_4 = \{e_1,e_4,e_5,e_9\} \leftrightarrow \{v_1,v_2,v_3,v_5\} = 6+11+10+8 = 35$;
$z_5 = \{e_1,e_2,e_6,e_8\} \leftrightarrow \{v_1,v_2,v_3,v_4\} = 6+10+9+9 = 34$;
$z_6 = \{e_1,e_2,e_7,e_9\} \leftrightarrow \{v_1,v_2,v_3,v_5\} = 6+10+7+8 = 31$;
$z_7 = \{e_3,e_4,e_6,e_7\} \leftrightarrow \{v_1,v_2,v_4,v_5\} = 5+11+9+7 = 32$;
$z_8 = \{e_1,e_4,e_6,e_{10}\} \leftrightarrow \{v_1,v_2,v_4,v_5\} = 6+11+9+40 = 66$;
$z_9 = \{e_1,e_3,e_7,e_{10}\} \leftrightarrow \{v_1,v_2,v_4,v_5\} = 6+5+7+40 = 58$;
$z_{10} = \{e_3,e_4,e_8,e_9\} \leftrightarrow \{v_1,v_3,v_4,v_5\} = 5+11+9+8 = 33$;
$z_{11} = \{e_2,e_4,e_8,e_{10}\} \leftrightarrow \{v_1,v_3,v_4,v_5\} = 10+11+9+40 = 70$;
$z_{12} = \{e_2,e_3,e_9,e_{10}\} \leftrightarrow \{v_1,v_3,v_4,v_5\} = 10+5+8+40 = 63$;
$z_{13} = \{e_6,e_7,e_8,e_9\} \leftrightarrow \{v_2,v_3,v_4,v_5\} = 9+7+9+8 = 33$;
$z_{14} = \{e_5,e_7,e_8,e_{10}\} \leftrightarrow \{v_2,v_3,v_4,v_5\} = 10+7+9+40 = 66$;
$z_{15} = \{e_5,e_6,e_9,e_{10}\} \leftrightarrow \{v_2,v_3,v_4,v_5\} = 10+9+8+40 = 67$;

Разобьём подмножество циклов, сгруппировав их относительно четвёрок вершин:

$z_1 = \{e_2,e_3,e_5,e_6\} \leftrightarrow \{v_1,v_2,v_3,v_4\} = 10+5+10+9 = 34$;
$z_5 = \{e_1,e_2,e_6,e_8\} \leftrightarrow \{v_1,v_2,v_3,v_4\} = 6+10+9+9 = 34$;
$z_3 = \{e_1,e_3,e_5,e_8\} \leftrightarrow \{v_1,v_2,v_3,v_4\} = 6+5+10+9 = 30$;

$z_2 = \{e_2,e_4,e_5,e_7\} \leftrightarrow \{v_1,v_2,v_3,v_5\} = 10+11+10+7 = 38$;
$z_4 = \{e_1,e_4,e_5,e_9\} \leftrightarrow \{v_1,v_2,v_3,v_5\} = 6+11+10+8 = 35$;
$z_6 = \{e_1,e_2,e_7,e_9\} \leftrightarrow \{v_1,v_2,v_3,v_5\} = 6+10+7+8 = 31$;

$z_7 = \{e_3,e_4,e_6,e_7\} \leftrightarrow \{v_1,v_2,v_4,v_5\} = 5+11+9+7 = 32$;
$z_8 = \{e_1,e_4,e_6,e_{10}\} \leftrightarrow \{v_1,v_2,v_4,v_5\} = 6+11+9+40 = 66$;
$z_9 = \{e_1,e_3,e_7,e_{10}\} \leftrightarrow \{v_1,v_2,v_4,v_5\} = 6+5+7+40 = 58$;

$z_{10} = \{e_3,e_4,e_8,e_9\} \leftrightarrow \{v_1,v_3,v_4,v_5\} = 5+11+9+8 = 33$;
$z_{11} = \{e_2,e_4,e_8,e_{10}\} \leftrightarrow \{v_1,v_3,v_4,v_5\} = 10+11+9+40 = 70$;
$z_{12} = \{e_2,e_3,e_9,e_{10}\} \leftrightarrow \{v_1,v_3,v_4,v_5\} = 10+5+8+40 = 63$;

$z_{13} = \{e_6,e_7,e_8,e_9\} \leftrightarrow \{v_2,v_3,v_4,v_5\} = 9+7+9+8 = 33$;

$z_{14} = \{e_5,e_7,e_8,e_{10}\} \leftrightarrow \{v_2,v_3,v_4,v_5\} = 10+7+9+40 = 66;$
$z_{15} = \{e_5,e_6,e_9,e_{10}\} \leftrightarrow \{v_2,v_3,v_4,v_5\} = 10+9+8+40 = 67;$

Выделим из троек минимальные циклы длиной четыре:

$z_3 = \{e_1,e_3,e_5,e_8\} \leftrightarrow \{v_1,v_2,v_3,v_4\} = 6+5+10+9 = 30;$
$z_6 = \{e_1,e_2,e_7,e_9\} \leftrightarrow \{v_1,v_2,v_3,v_5\} = 6+10+7+8 = 31;$
$z_7 = \{e_3,e_4,e_6,e_7\} \leftrightarrow \{v_1,v_2,v_4,v_5\} = 5+11+9+7 = 32;$
$z_{10} = \{e_3,e_4,e_8,e_9\} \leftrightarrow \{v_1,v_3,v_4,v_5\} = 5+11+9+8 = 33;$
$z_{13} = \{e_6,e_7,e_8,e_9\} \leftrightarrow \{v_2,v_3,v_4,v_5\} = 9+7+9+8 = 33;$

Для выделенного цикла длиной четыре
$z_3 = \{e_1,e_3,e_5,e_8\} \leftrightarrow \{v_1,v_2,v_3,v_4\} = 6+5+10+9 = 30;$
выделим соприкасающиеся циклы:
цикл $c_3 = \{e_1,e_4,e_7\} \leftrightarrow \{v_1,v_2,v_5\} = 11+7+6 = 24;$
цикл $c_6 = \{e_3,e_4,e_{10}\} \leftrightarrow \{v_1,v_4,v_5\} = 11+5+40 = 56;$
цикл $c_8 = \{e_5,e_7,e_9\} \leftrightarrow \{v_2,v_3,v_5\} = 10+7+8 = 26;$
цикл $c_{10} = \{e_8,e_9,e_{10}\} \leftrightarrow \{v_3,v_4,v_5\} = 40+9+8 = 57.$

Строим гамильтонов цикл длиной пять с минимальным весом:

$z_3 \oplus c_3 = \{e_1,e_3,e_5,e_8\} \oplus \{e_1,e_4,e_7\}=\{e_3,e_4,e_5,e_7,e_8\} \leftrightarrow \{v_1,v_2,v_3,v_4,v_5\} = 5+10+9+11+7 = 42;$
$z_3 \oplus c_6 = \{e_1,e_3,e_5,e_8\} \oplus \{e_3,e_4,e_{10}\}=\{e_1,e_4,e_5,e_8,e_{10}\} \leftrightarrow \{v_1,v_2,v_3,v_4,v_5\} = 6+10+9+11+40 = 76;$
$z_3 \oplus c_8 = \{e_1,e_3,e_5,e_8\} \oplus \{e_5,e_7,e_9\}=\{e_1,e_3,e_7,e_8,e_9\} \leftrightarrow \{v_1,v_2,v_3,v_4,v_5\} =6+5+9+8+7 = 35;$
$z_3 \oplus c_{10} = \{e_1,e_3,e_5,e_8\} \oplus \{e_8,e_9,e_{10}\}=\{e_1,e_3,e_5,e_9,e_{10}\} \leftrightarrow \{v_1,v_2,v_3,v_4,v_5\} =6+5+10+8+40 = 69;$

Выбираем оптимальный маршрут:
$z_3 \oplus c_8 = \{e_1,e_3,e_5,e_8\} \oplus \{e_5,e_7,e_9\}=\{e_1,e_3,e_7,e_8,e_9\} \leftrightarrow \{v_1,v_2,v_3,v_4,v_5\} =6+5+9+8+7 = 35.$

Продолжим эксперимент. Увеличим стоимость ребра $e_5$, так как существует предположение, что увеличение стоимости данного ребра не приведёт к выбору оптимального маршрута.

ребро $e_1 \leftrightarrow (v_1,v_2) = 6;$
ребро $e_2 \leftrightarrow (v_1,v_3) = 10;$
ребро $e_3 \leftrightarrow (v_1,v_4) = 5;$
ребро $e_4 \leftrightarrow (v_1,v_5) = 11;$
ребро $e_5 \leftrightarrow (v_2,v_3) = 50;$
ребро $e_6 \leftrightarrow (v_2,v_4) = 9;$
ребро $e_7 \leftrightarrow (v_2,v_5) = 7;$
ребро $e_8 \leftrightarrow (v_3,v_4) = 9;$
ребро $e_9 \leftrightarrow (v_3,v_5) = 8;$
ребро $e_{10} \leftrightarrow (v_4,v_5) = 40.$

Определим стоимость изометрического цикла:
цикл $c_1 = \{e_1,e_2,e_5\} \leftrightarrow \{v_1,v_2,v_3\} = 50+5+6 = 71;$
цикл $c_2 = \{e_1,e_3,e_6\} \leftrightarrow \{v_1,v_2,v_4\} = 9+5+6 = 20;$
цикл $c_3 = \{e_1,e_4,e_7\} \leftrightarrow \{v_1,v_2,v_5\} = 11+7+6 = 24;$
цикл $c_4 = \{e_2,e_3,e_8\} \leftrightarrow \{v_1,v_3,v_4\} = 10+5+9 = 24;$
цикл $c_5 = \{e_2,e_4,e_9\} \leftrightarrow \{v_1,v_3,v_5\} = 10+11+8 = 29;$
цикл $c_6 = \{e_3,e_4,e_{10}\} \leftrightarrow \{v_1,v_4,v_5\} = 11+5+40 = 56;$
цикл $c_7 = \{e_5,e_6,e_8\} \leftrightarrow \{v_2,v_3,v_4\} = 50+9+9 = 68;$
цикл $c_8 = \{e_5,e_7,e_9\} \leftrightarrow \{v_2,v_3,v_5\} = 50+7+8 = 65;$
цикл $c_9 = \{e_6,e_7,e_{10}\} \leftrightarrow \{v_2,v_4,v_5\} = 9+7+40 = 56;$

цикл $c_{10} = \{e_8,e_9,e_{10}\} \leftrightarrow \{v_3,v_4,v_5\}$ = 40+9+8 = 57.

Построим циклы длиной четыре:

$z_1 = \{e_2,e_3,e_5,e_6\} \leftrightarrow \{v_1,v_2,v_3,v_4\}$ = 10+5+50+9 = 74;
$z_2 = \{e_2,e_4,e_5,e_7\} \leftrightarrow \{v_1,v_2,v_3,v_5\}$ = 10+11+50+7 = 78;
$z_3 = \{e_1,e_3,e_5,e_8\} \leftrightarrow \{v_1,v_2,v_3,v_4\}$ = 6+5+50+9 = 70;
$z_4 = \{e_1,e_4,e_5,e_9\} \leftrightarrow \{v_1,v_2,v_3,v_5\}$ = 6+11+50+8 = 75;
$z_5 = \{e_1,e_2,e_6,e_8\} \leftrightarrow \{v_1,v_2,v_3,v_4\}$ = 6+10+9+9 = 34;
$z_6 = \{e_1,e_2,e_7,e_9\} \leftrightarrow \{v_1,v_2,v_3,v_5\}$ = 6+10+7+8 = 31;
$z_7 = \{e_3,e_4,e_6,e_7\} \leftrightarrow \{v_1,v_2,v_4,v_5\}$ = 5+11+9+7 = 32;
$z_8 = \{e_1,e_4,e_6,e_{10}\} \leftrightarrow \{v_1,v_2,v_4,v_5\}$ = 6+11+9+40 = 66;
$z_9 = \{e_1,e_3,e_7,e_{10}\} \leftrightarrow \{v_1,v_2,v_4,v_5\}$ = 6+5+7+40 = 58;
$z_{10} = \{e_3,e_4,e_8,e_9\} \leftrightarrow \{v_1,v_3,v_4,v_5\}$ = 5+11+9+8 = 33;
$z_{11} = \{e_2,e_4,e_8,e_{10}\} \leftrightarrow \{v_1,v_3,v_4,v_5\}$ = 10+11+9+40 = 70;
$z_{12} = \{e_2,e_3,e_9,e_{10}\} \leftrightarrow \{v_1,v_3,v_4,v_5\}$ = 10+5+8+40 = 63;
$z_{13} = \{e_6,e_7,e_8,e_9\} \leftrightarrow \{v_2,v_3,v_4,v_5\}$ = 9+7+9+8 = 33;
$z_{14} = \{e_5,e_7,e_8,e_{10}\} \leftrightarrow \{v_2,v_3,v_4,v_5\}$ = 50+7+9+40 = 106;
$z_{15} = \{e_5,e_6,e_9,e_{10}\} \leftrightarrow \{v_2,v_3,v_4,v_5\}$ = 50+9+8+40 = 107;

Циклы длиной четыре (с группировкой по вершинам):

$z_1 = \{e_2,e_3,e_5,e_6\} \leftrightarrow \{v_1,v_2,v_3,v_4\}$ = 10+5+50+9 = 74;
$z_5 = \{e_1,e_2,e_6,e_8\} \leftrightarrow \{v_1,v_2,v_3,v_4\}$ = 6+10+9+9 = 34;
$z_3 = \{e_1,e_3,e_5,e_8\} \leftrightarrow \{v_1,v_2,v_3,v_4\}$ = 6+5+50+9 = 70;

$z_2 = \{e_2,e_4,e_5,e_7\} \leftrightarrow \{v_1,v_2,v_3,v_5\}$ = 10+11+50+7 = 78;
$z_4 = \{e_1,e_4,e_5,e_9\} \leftrightarrow \{v_1,v_2,v_3,v_5\}$ = 6+11+50+8 = 75;
$z_6 = \{e_1,e_2,e_7,e_9\} \leftrightarrow \{v_1,v_2,v_3,v_5\}$ = 6+10+7+8 = 31;

$z_7 = \{e_3,e_4,e_6,e_7\} \leftrightarrow \{v_1,v_2,v_4,v_5\}$ = 5+11+9+7 = 32;
$z_8 = \{e_1,e_4,e_6,e_{10}\} \leftrightarrow \{v_1,v_2,v_4,v_5\}$ = 6+11+9+40 = 66;
$z_9 = \{e_1,e_3,e_7,e_{10}\} \leftrightarrow \{v_1,v_2,v_4,v_5\}$ = 6+5+7+40 = 58;

$z_{10} = \{e_3,e_4,e_8,e_9\} \leftrightarrow \{v_1,v_3,v_4,v_5\}$ = 5+11+9+8 = 33;
$z_{11} = \{e_2,e_4,e_8,e_{10}\} \leftrightarrow \{v_1,v_3,v_4,v_5\}$ = 10+11+9+40 = 70;
$z_{12} = \{e_2,e_3,e_9,e_{10}\} \leftrightarrow \{v_1,v_3,v_4,v_5\}$ = 10+5+8+40 = 63;

$z_{13} = \{e_6,e_7,e_8,e_9\} \leftrightarrow \{v_2,v_3,v_4,v_5\}$ = 9+7+9+8 = 33;
$z_{14} = \{e_5,e_7,e_8,e_{10}\} \leftrightarrow \{v_2,v_3,v_4,v_5\}$ = 50+7+9+40 = 106;
$z_{15} = \{e_5,e_6,e_9,e_{10}\} \leftrightarrow \{v_2,v_3,v_4,v_5\}$ = 50+9+8+40 = 107;

Минимальные циклы длиной четыре:

$z_3 = \{e_1,e_3,e_5,e_8\} \leftrightarrow \{v_1,v_2,v_3,v_4\}$ = 6+5+50+9 = 70;
$z_6 = \{e_1,e_2,e_7,e_9\} \leftrightarrow \{v_1,v_2,v_3,v_5\}$ = 6+10+7+8 = 31;
$z_7 = \{e_3,e_4,e_6,e_7\} \leftrightarrow \{v_1,v_2,v_4,v_5\}$ = 5+11+9+7 = 32;
$z_{10} = \{e_3,e_4,e_8,e_9\} \leftrightarrow \{v_1,v_3,v_4,v_5\}$ = 5+11+9+8 = 33;
$z_{13} = \{e_6,e_7,e_8,e_9\} \leftrightarrow \{v_2,v_3,v_4,v_5\}$ = 9+7+9+8 = 33;

Строим гамильтонов цикл с минимальным весом для цикла

$z_6 = \{e_1,e_2,e_7,e_9\} \leftrightarrow \{v_1,v_2,v_3,v_5\}$ = 6+10+7+8 = 31;

выбираем сопряжённые циклы из множества изометрических циклов:

цикл $c_2 = \{e_1,e_3,e_6\} \leftrightarrow \{v_1,v_2,v_4\} = 9+5+6 = 20$;
цикл $c_4 = \{e_2,e_3,e_8\} \leftrightarrow \{v_1,v_3,v_4\} = 10+5+9 = 24$;
цикл $c_9 = \{e_6,e_7,e_{10}\} \leftrightarrow \{v_2,v_4,v_5\} = 9+7+40 = 56$;
цикл $c_{10} = \{e_8,e_9,e_{10}\} \leftrightarrow \{v_3,v_4,v_5\} = 40+9+8 = 57$.

Определяем гамильтонов цикл с минимальным весом:

$z_6 \oplus c_4 = \{e_1,e_2,e_7,e_9\} \oplus \{e_2,e_3,e_8\}=\{e_1,e_3,e_7,e_8,e_9\} \leftrightarrow \{v_1,v_2,v_3,v_4,v_5\} = 6+7+8+5+9 = 35$;

Таким образом, изменение в сторону увеличения стоимости рёбер, не принадлежащих оптимальному маршруту, не влияет на результат решения задачи.

### 3.5. Задача коммивояжёра для полного графа $K_6$

Рассмотрим задачу о коммивояжёре, представленную следующей матрицей стоимости:

| | $v_1$ | $v_2$ | $v_3$ | $v_4$ | $v_5$ | $v_6$ |
|---|---|---|---|---|---|---|
| $v_1$ | | 6 | 4 | 8 | 7 | 14 |
| $v_2$ | 6 | | 7 | 11 | 7 | 10 |
| $v_3$ | 4 | 7 | | 4 | 3 | 10 |
| $v_4$ | 8 | 11 | 4 | | 5 | 11 |
| $v_5$ | 7 | 7 | 3 | 5 | | 7 |
| $v_6$ | 14 | 10 | 10 | 11 | 7 | |

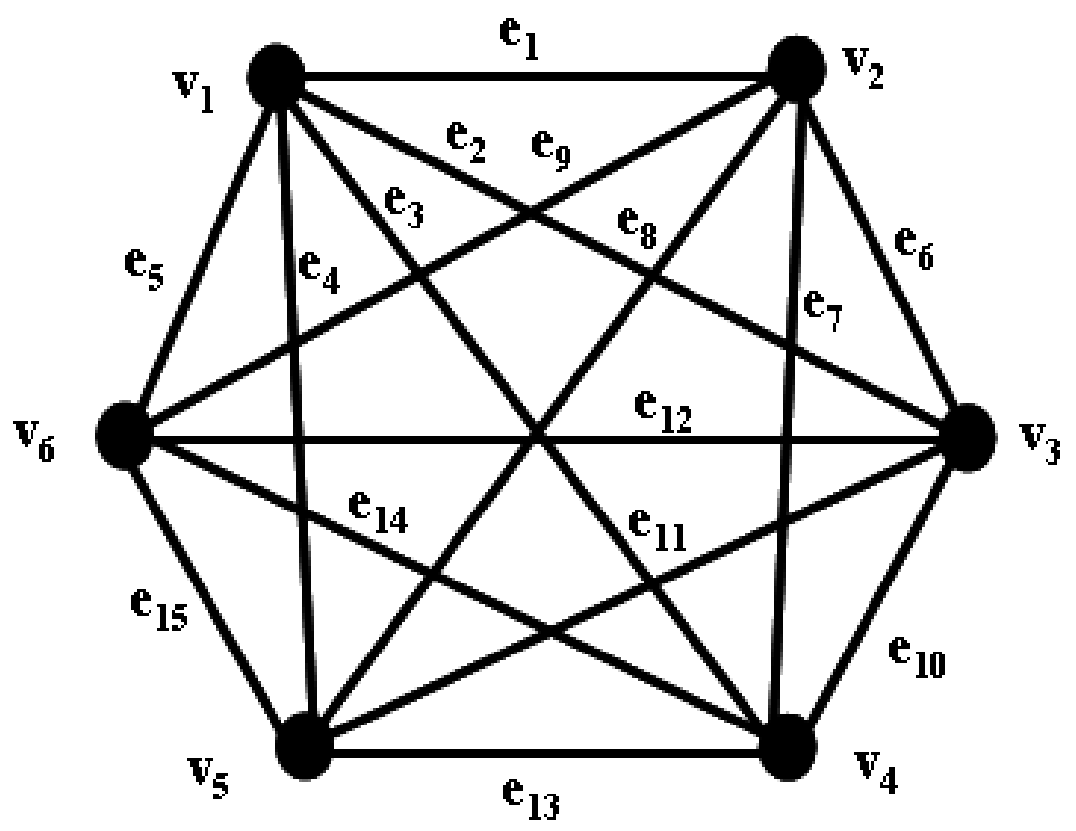

Рис. 3.6. Матрица стоимости и граф с правильной нумерацией ребер

Для полного графа $K_6$ $m$=15.

Веса рёбер графа:

ребро $e_1 = (v_1,v_2) = 6$;
ребро $e_2 = (v_1,v_3) = 4$;
ребро $e_3 = (v_1,v_4) = 8$;
ребро $e_4 = (v_1,v_5) = 7$;
ребро $e_5 = (v_1,v_6) = 14$;
ребро $e_6 = (v_2,v_3) = 7$;
ребро $e_7 = (v_2,v_4) = 11$;
ребро $e_8 = (v_2,v_5) = 7$;
ребро $e_9 = (v_2,v_6) = 10$;
ребро $e_{10} = (v_3,v_4) = 4$;
ребро $e_{11} = (v_3,v_5) = 3$;
ребро $e_{12} = (v_3,v_6) = 10$;
ребро $e_{13} = (v_4,v_5) = 5$;
ребро $e_{14} = (v_4,v_6) = 11$;
ребро $e_{15} = (v_5,v_6) = 7$.

Изометрические циклы графа:

$c_1 = \{e_1,e_2,e_6\} \leftrightarrow \{v_1,v_2,v_3\}=6+4+7=17$;
$c_2 = \{e_1,e_3,e_7\} \leftrightarrow \{v_1,v_2,v_4\}=6+8+11=25$;

$c_3 = \{e_1,e_4,e_8\} \leftrightarrow \{v_1,v_2,v_5\}=6+7+7=20$; $c_4 = \{e_1,e_5,e_9\} \leftrightarrow \{v_1,v_2,v_6\}=6+14+10=30$;
$c_5 = \{e_2,e_3,e_{10}\} \leftrightarrow \{v_1,v_3,v_4\}=4+8+4=16$; $c_6 = \{e_2,e_4,e_{11}\} \leftrightarrow \{v_1,v_3,v_5\}=4+7+3=14$;
$c_7 = \{e_2,e_5,e_{12}\} \leftrightarrow \{v_1,v_3,v_6\}=4+14+10=28$; $c_8 = \{e_3,e_4,e_{13}\} \leftrightarrow \{v_1,v_4,v_5\}=8+7+5=20$;
$c_9 = \{e_3,e_5,e_{14}\} \leftrightarrow \{v_1,v_4,v_6\}=8+14+11=32$; $c_{10} = \{e_4,e_5,e_{15}\} \leftrightarrow \{v_1,v_5,v_6\}=7+14+5=28$;
$c_{11} = \{e_6,e_7,e_{10}\} \leftrightarrow \{v_2,v_3,v_4\}=7+11+4=22$; $c_{12} = \{e_6,e_8,e_{11}\} \leftrightarrow \{v_2,v_3,v_5\}=7+7+3=15$;
$c_{13} = \{e_6,e_9,e_{12}\} \leftrightarrow \{v_2,v_3,v_6\}=7+10+10=27$; $c_{14} = \{e_7,e_8,e_{13}\} \leftrightarrow \{v_2,v_4,v_5\}=11+7+5=23$;
$c_{15} = \{e_7,e_9,e_{14}\} \leftrightarrow \{v_2,v_4,v_6\}=11+10+11=32$; $c_{16} = \{e_8,e_9,e_{15}\} \leftrightarrow \{v_2,v_5,v_6\}=7+10+7=24$;
$c_{17} = \{e_{10},e_{11},e_{13}\} \leftrightarrow \{v_3,v_4,v_5\}=4+3+5=12$; $c_{18} = \{e_{10},e_{12},e_{14}\} \leftrightarrow \{v_3,v_4,v_6\}=4+10+11=24$;
$c_{19} = \{e_{11},e_{12},e_{15}\} \leftrightarrow \{v_3,v_5,v_6\}=3+10+7=20$; $c_{20} = \{e_{13},e_{14},e_{15}\} \leftrightarrow \{v_4,v_5,v_6\}=5+11+7=23$.

Для каждого набора из четырёх вершин построим по три цикла длиной четыре.

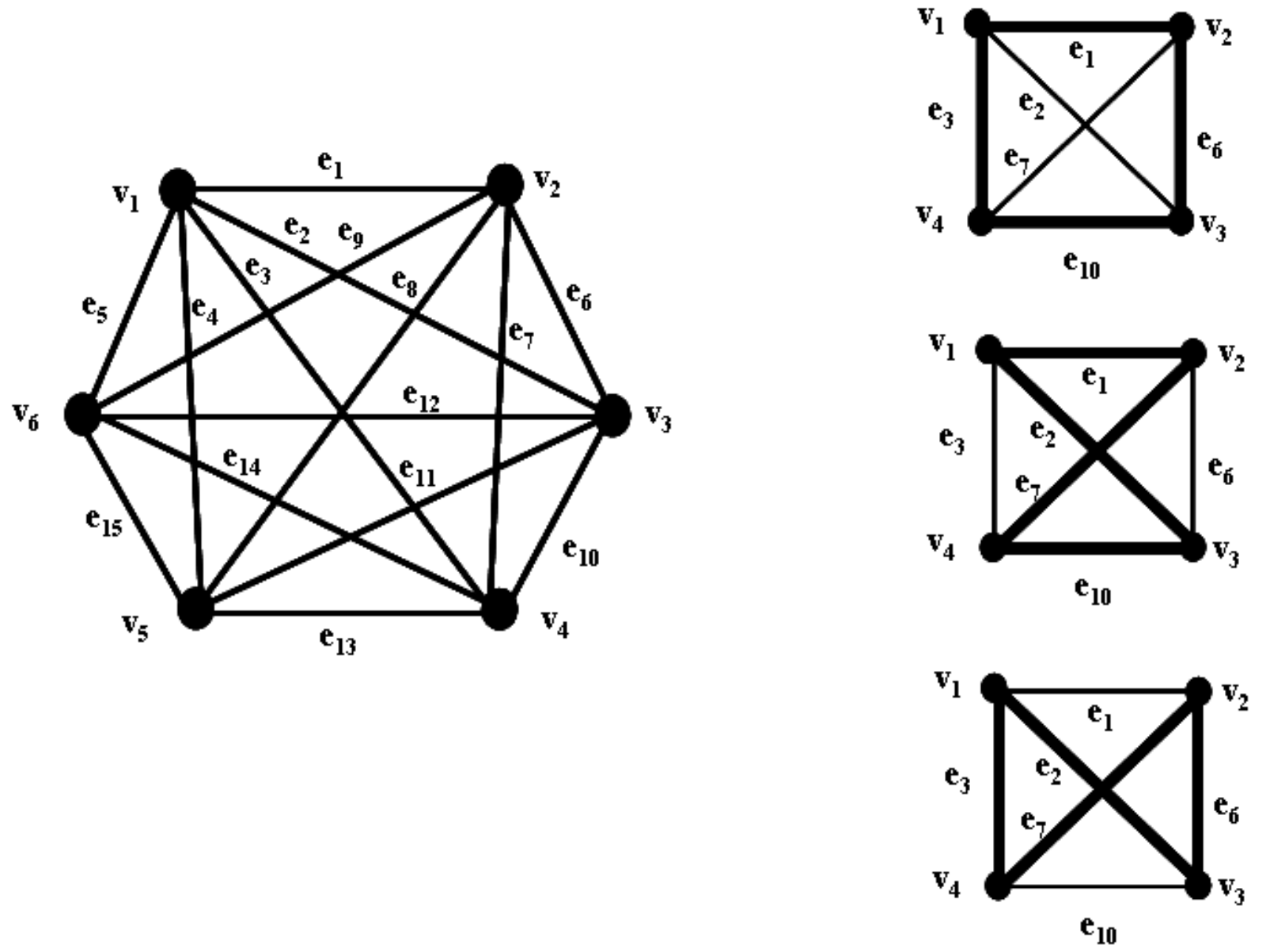

Рис. 3.7. Циклы для подмножества вершин $\{v_1,v_2,v_3,v_4\}$.

$z_1 = \{e_1,e_4,e_6,e_{10}\} = \{v_1,v_2,v_3,v_4\} = 6+7+7+4 = 24$;
$z_2 = \{e_1,e_2,e_7,e_{10}\} = \{v_1,v_2,v_3,v_4\} = 6+4+11+4 = 25$;
$z_3 = \{e_2,e_3,e_6,e_7\} = \{v_1,v_2,v_3,v_4\} = 48+7+11 = 30$;

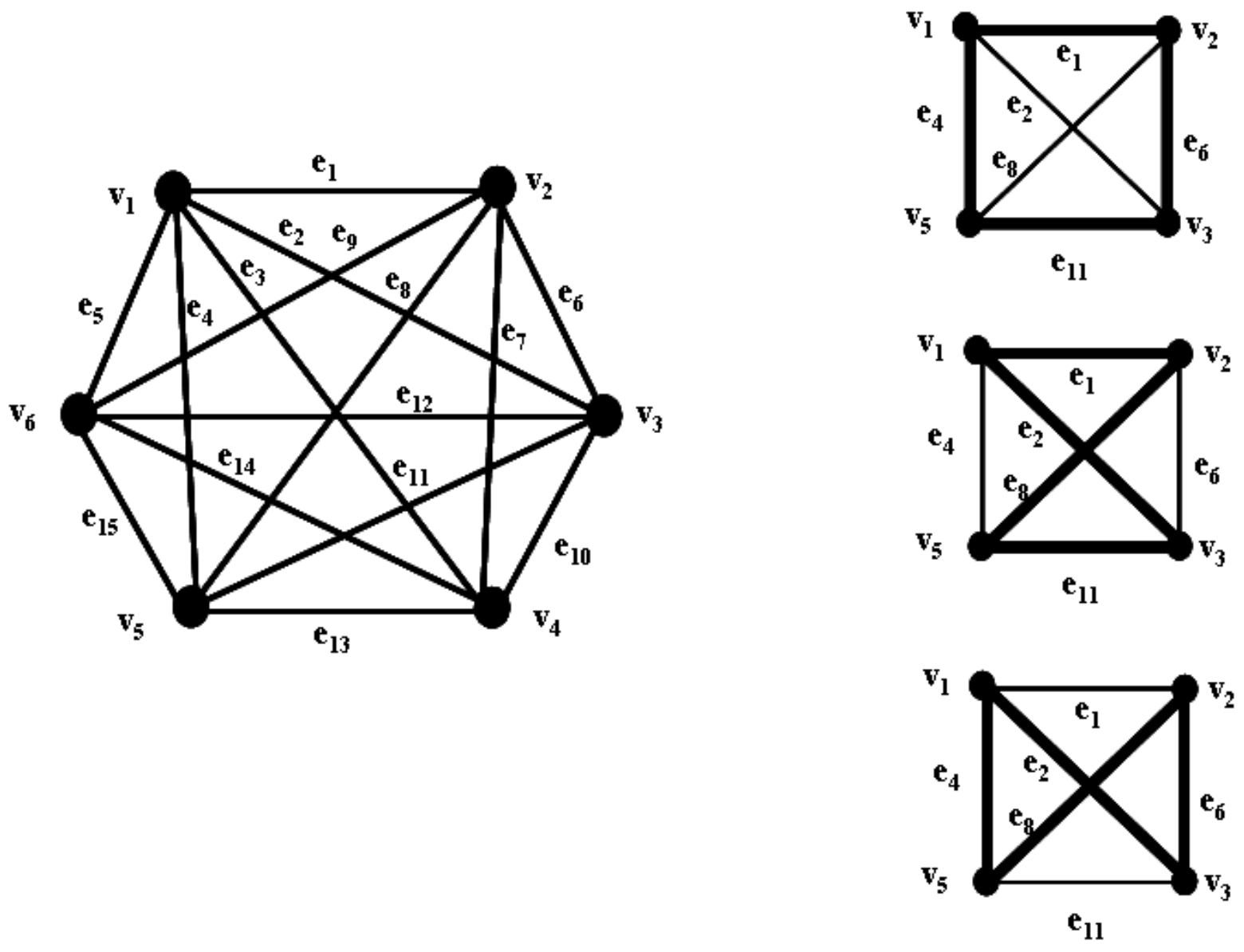

Рис. 3.8. Циклы для подмножества вершин $\{v_1,v_2,v_3,v_5\}$.

$z_4 = \{e_1,e_4,e_6,e_{11}\} = \{v_1,v_2,v_3,v_5\} = 6+7+7+3 = 23;$
$z_5 = \{e_1,e_2,e_8,e_{11}\} = \{v_1,v_2,v_3,v_5\} = 6+4+7+3 = 20;$
$z_6 = \{e_2,e_4,e_6,e_8\} = \{v_1,v_2,v_3,v_5\} = 4+7+7+7 = 25;$

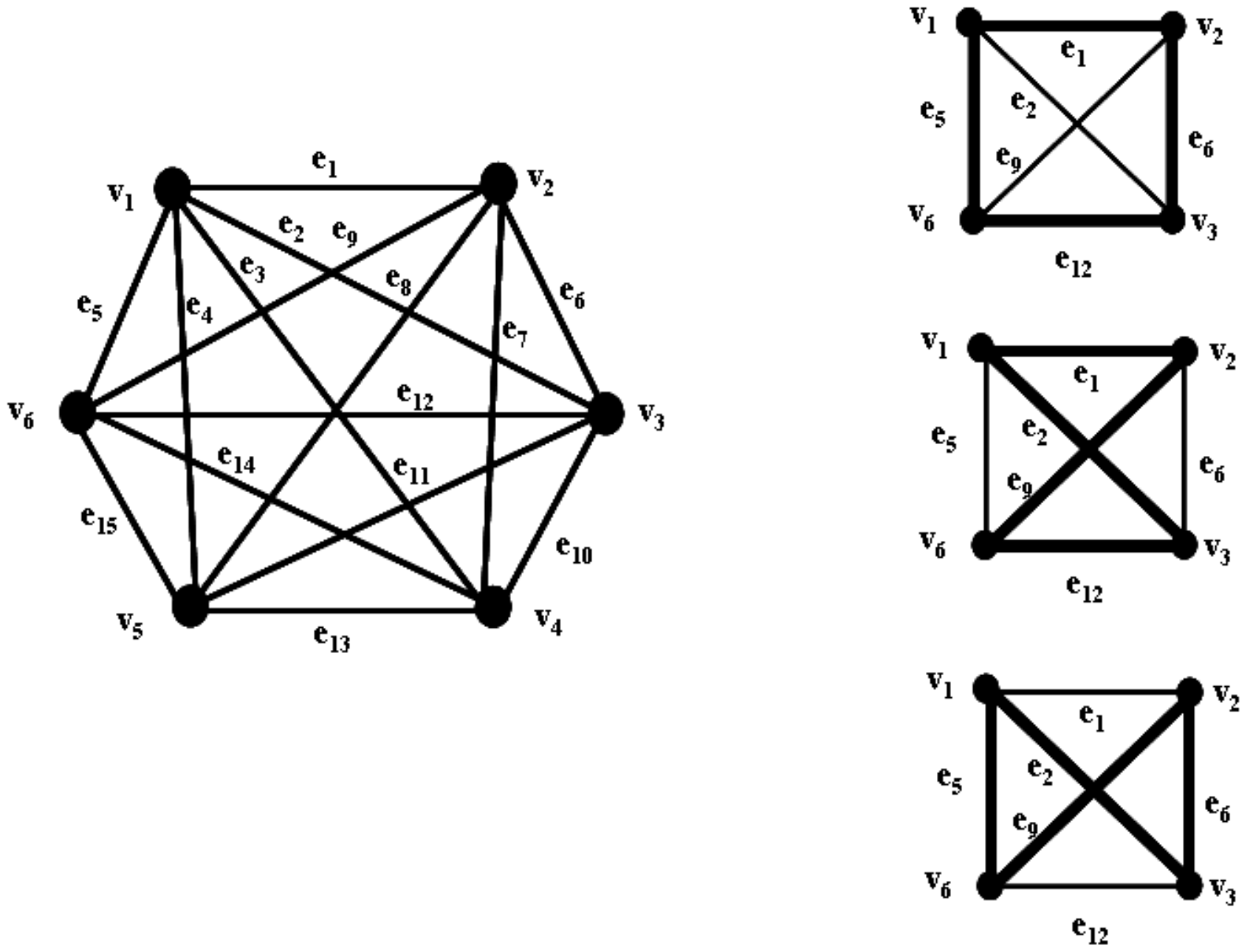

Рис. 3.9. Циклы для подмножества вершин $\{v_1,v_2,v_3,v_6\}$.

$z_7 = \{e_1,e_5,e_6,e_{12}\} = \{v_1,v_2,v_3,v_6\} = 6+4+7+10 = 27;$
$z_8 = \{e_1,e_2,e_9,e_{12}\} = \{v_1,v_2,v_3,v_6\} = 6+4+10+10 = 30;$
$z_9 = \{e_2,e_5,e_6,e_{12}\} = \{v_1,v_2,v_3,v_6\} = 4+14+7+10 = 35;$

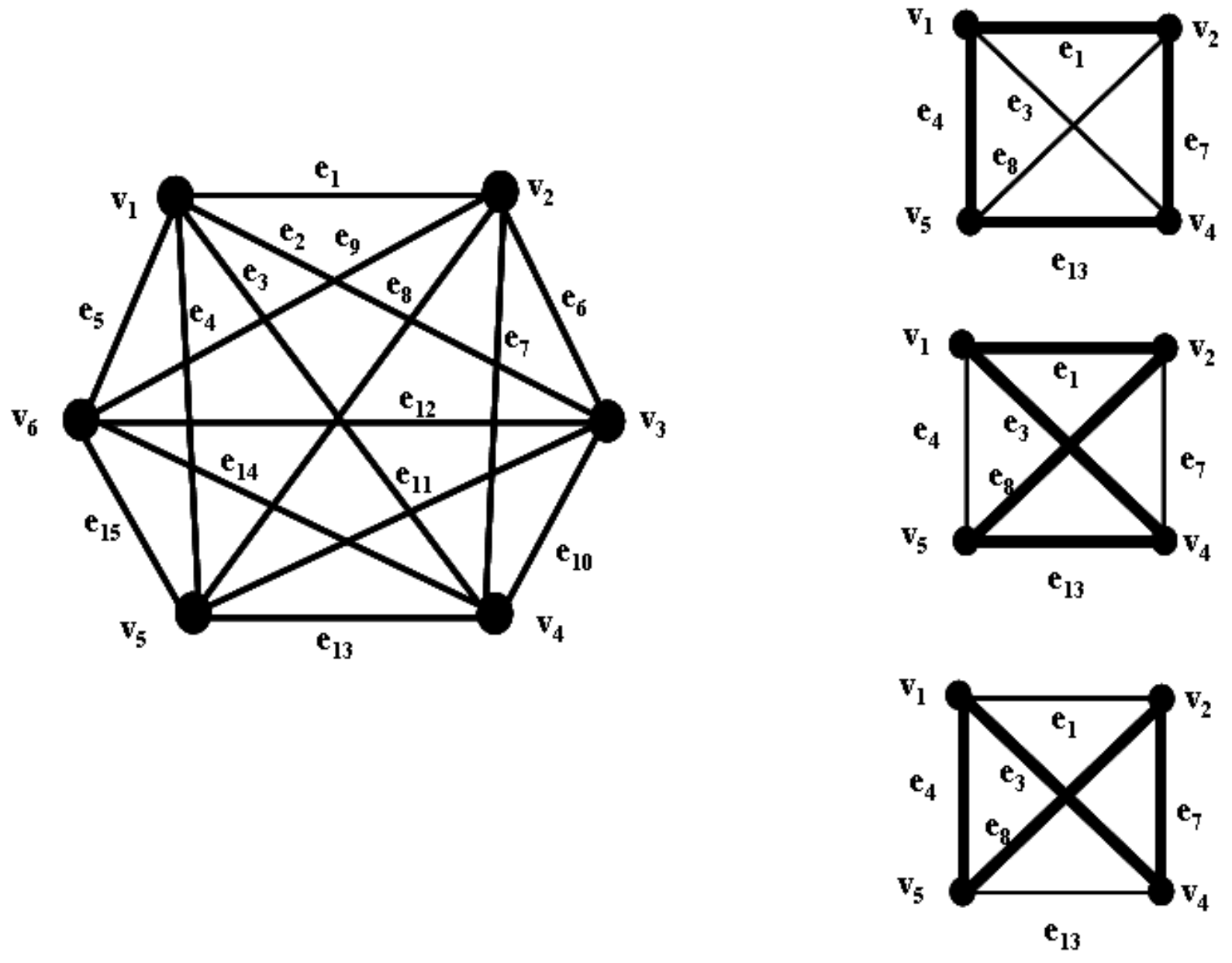

Рис. 3.10. Циклы для подмножества вершин $\{v_1,v_2,v_4,v_5\}$.

$z_{10} = \{e_1,e_4,e_7,e_{13}\} = \{v_1,v_2,v_4,v_5\} = 6+7+11+5 = 29;$
$z_{11} = \{e_1,e_3,e_8,e_{13}\} = \{v_1,v_2,v_4,v_5\} = 6+7+8+5 = 26;$
$z_{12} = \{e_3,e_4,e_7,e_8\} = \{v_1,v_2,v_4,v_5\} = 8+7+11+7 = 33;$

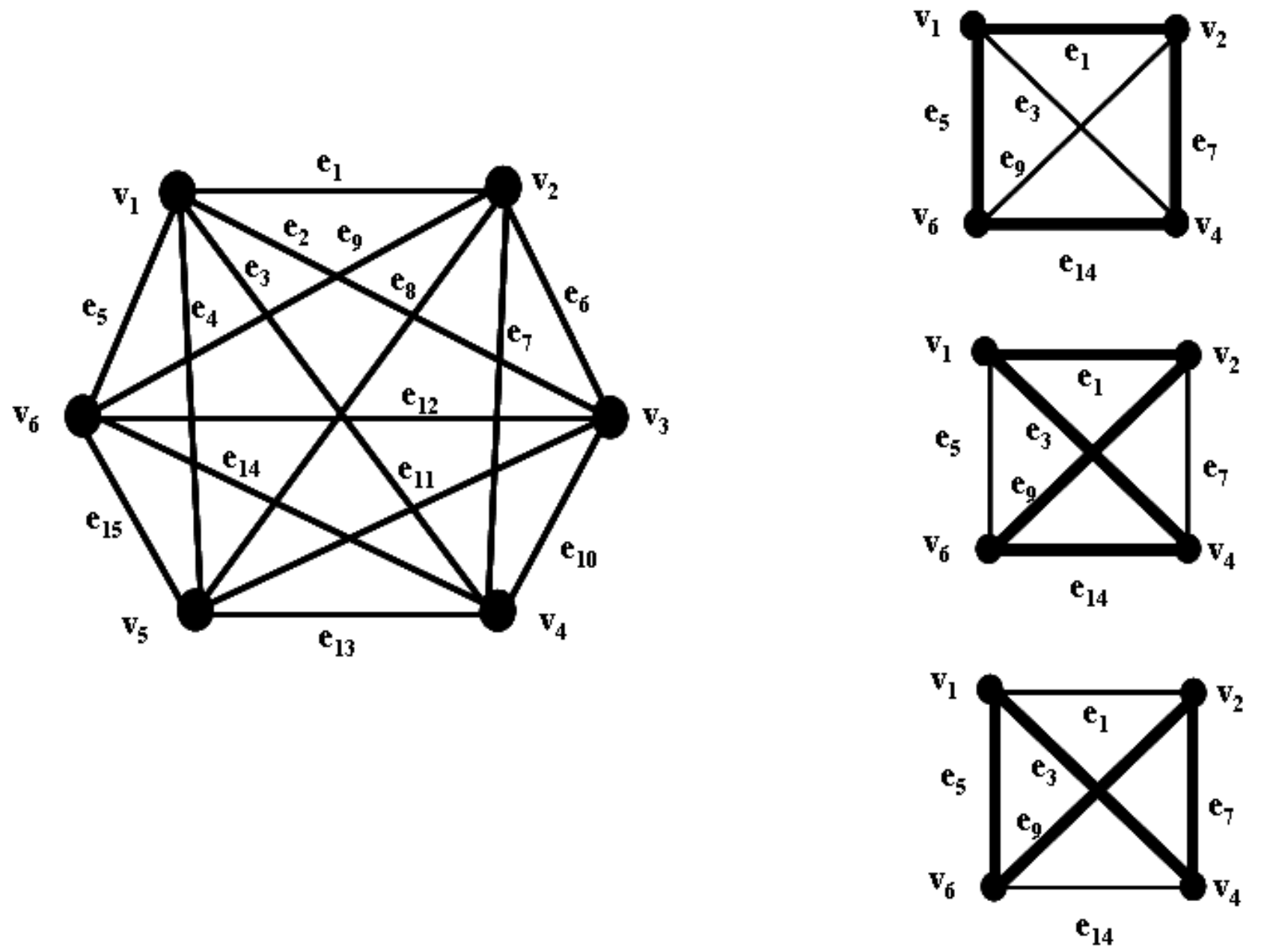

Рис. 3.11. Циклы для подмножества вершин {$v_1$,$v_2$,$v_4$,$v_6$}.

$z_{13}$ = {$e_1$,$e_5$,$e_7$,$e_{14}$} = {$v_1$,$v_2$,$v_4$,$v_6$} = 6+14+11+11 = 44;
$z_{14}$ = {$e_1$,$e_3$,$e_9$,$e_{14}$} = {$v_1$,$v_2$,$v_4$,$v_6$} = 6+8+10+7 = 31;
$z_{15}$ = {$e_3$,$e_5$,$e_7$,$e_9$} = {$v_1$,$v_2$,$v_4$,$v_6$} = 8+14+11+10 = 43;

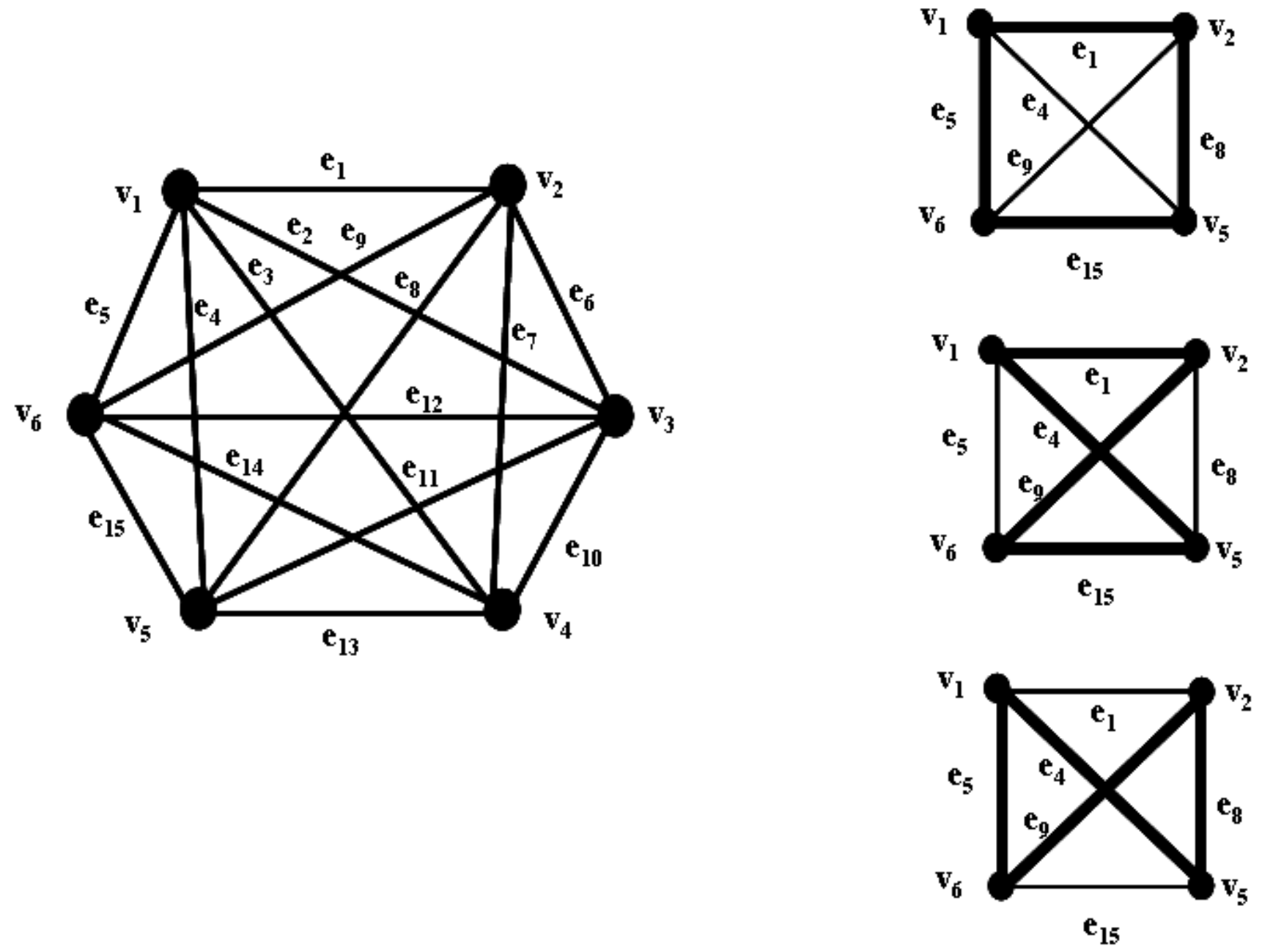

Рис. 3.12. Циклы для подмножества вершин {$v_1$,$v_2$,$v_5$,$v_6$}.

$z_{16}$ = {$e_1$,$e_5$,$e_8$,$e_{15}$} = {$v_1$,$v_2$,$v_5$,$v_6$} = 6+14+7+7 = 34;
$z_{17}$ = {$e_1$,$e_4$,$e_9$,$e_{15}$} = {$v_1$,$v_2$,$v_5$,$v_6$} = 6+7+10+7 = 30;
$z_{18}$ = {$e_4$,$e_5$,$e_8$,$e_9$} = {$v_1$,$v_2$,$v_5$,$v_6$} = 7+14+7+10 = 38;

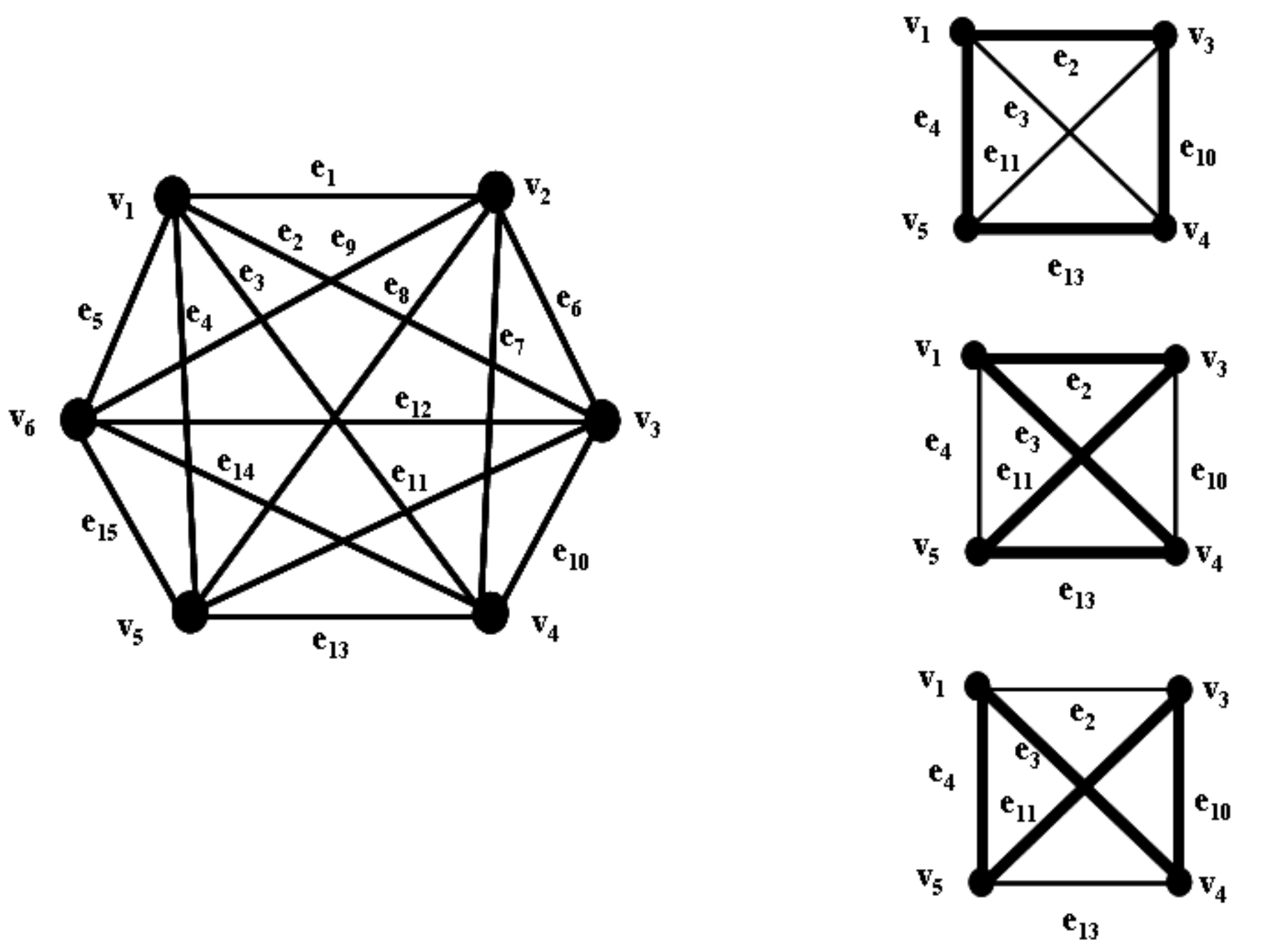

Рис. 3.13. Циклы для подмножества вершин $\{v_1,v_3,v_4,v_5\}$.

$z_{19} = \{e_2,e_4,e_{10},e_{13}\} = \{v_1,v_3,v_4,v_5\} = 4+7+4+5 = 20;$
$z_{20} = \{e_2,e_3,e_{11},e_{13}\} = \{v_1,v_3,v_4,v_5\} = 4+8+3+5 = 20;$
$z_{21} = \{e_3,e_4,e_{10},e_{11}\} = \{v_1,v_3,v_4,v_5\} = 8+7+4+3 = 22;$

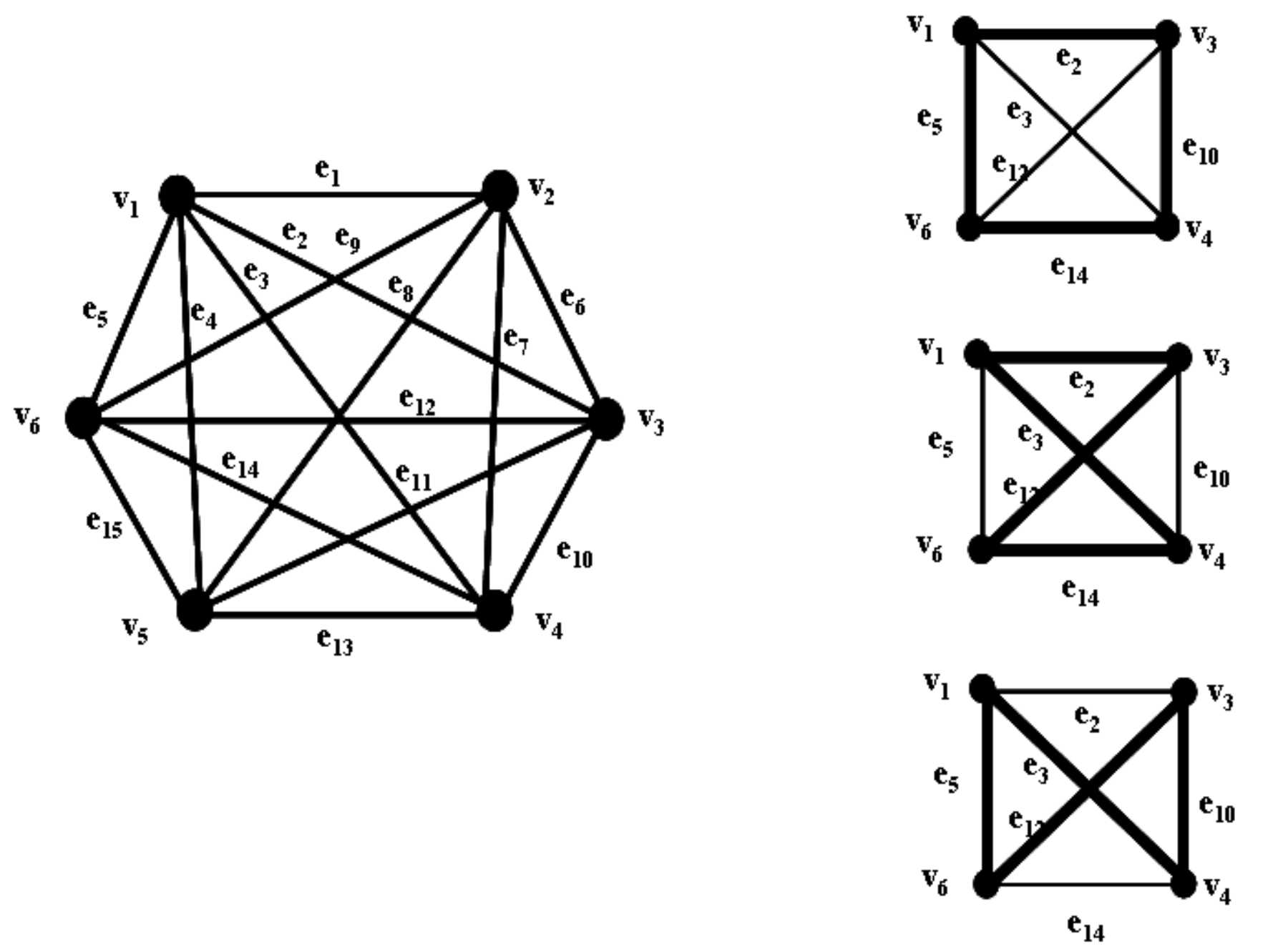

Рис. 3.14. Циклы для подмножества вершин $\{v_1,v_3,v_4,v_6\}$.

$z_{22} = \{e_2,e_5,e_{10},e_{14}\} = \{v_1,v_3,v_4,v_6\} = 4+8+4+11 = 26;$
$z_{23} = \{e_2,e_3,e_{12},e_{14}\} = \{v_1,v_3,v_4,v_6\} = 4+8+10+11 = 33;$
$z_{24} = \{e_3,e_5,e_{10},e_{12}\} = \{v_1,v_3,v_4,v_6\} = 8+14+4+10 = 36;$

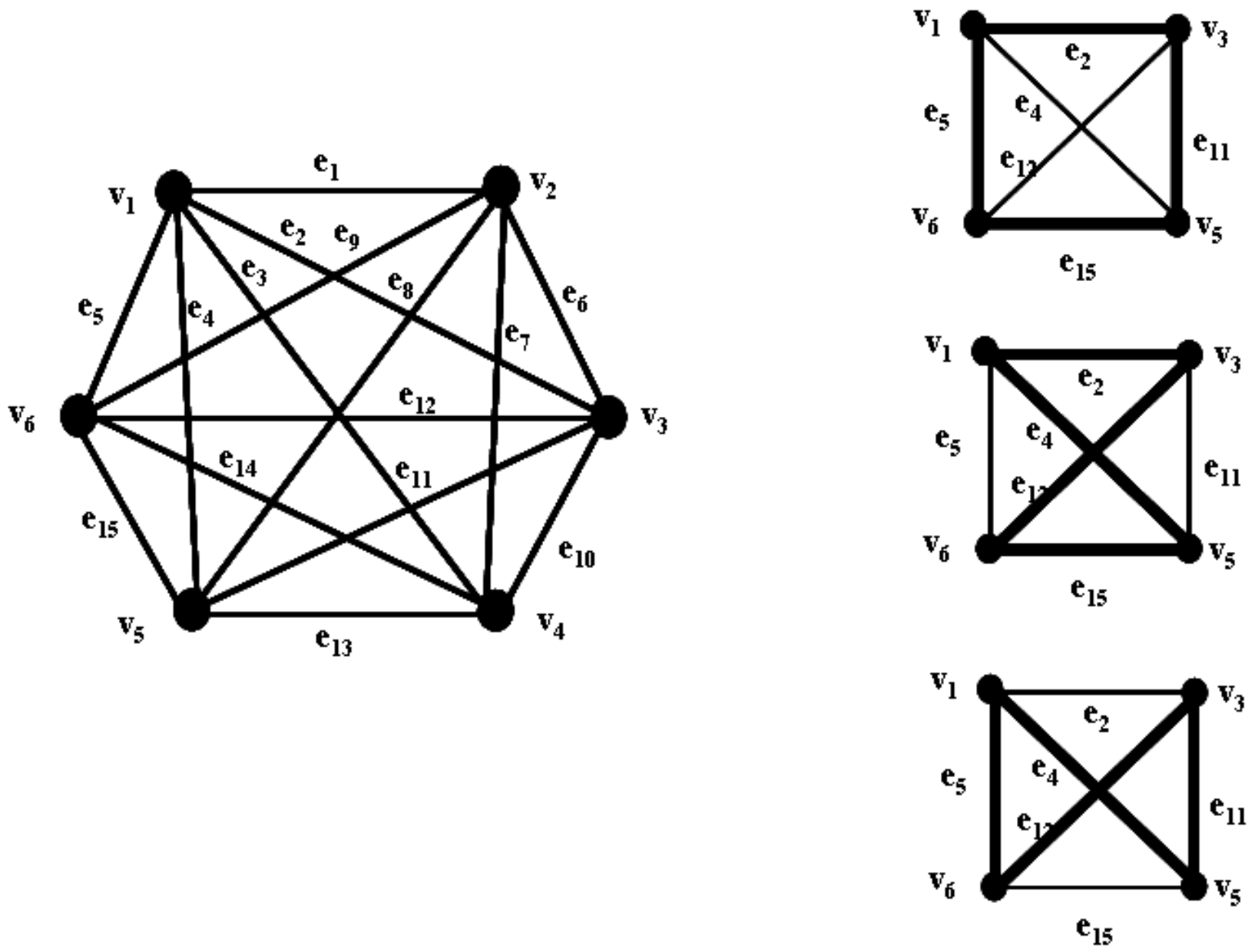

Рис. 3.15. Циклы для подмножества вершин $\{v_1,v_3,v_5,v_6\}$.

$z_{25} = \{e_2,e_5,e_{11},e_{15}\} = \{v_1,v_3,v_5,v_6\} = 4+14+3+7 = 28;$
$z_{26} = \{e_2,e_4,e_{12},e_{15}\} = \{v_1,v_3,v_5,v_6\} = 4+7+10+7 = 28;$
$z_{27} = \{e_4,e_5,e_{11},e_{12}\} = \{v_1,v_3,v_5,v_6\} = 7+14+3+10 = 34;$

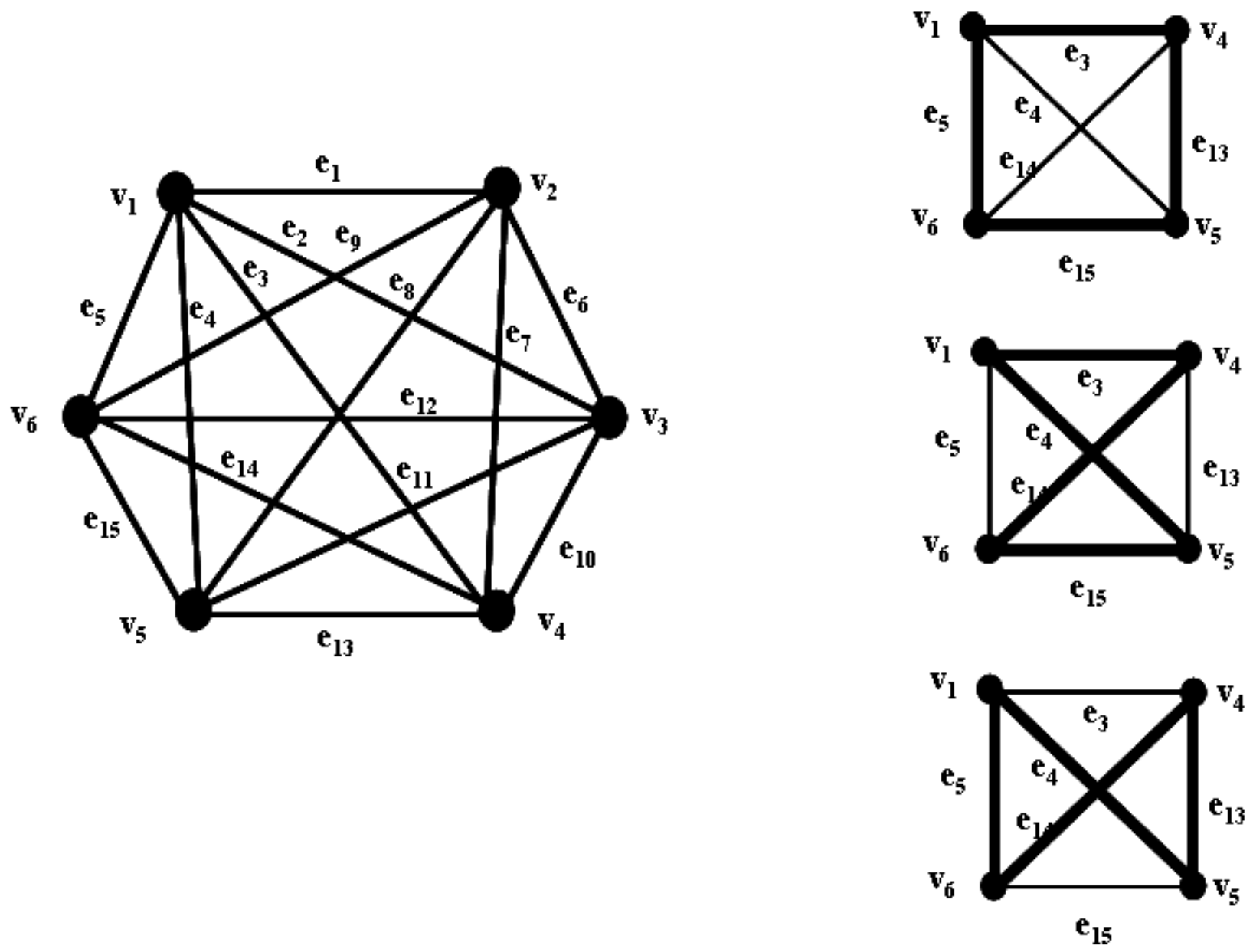

Рис. 3.16. Циклы для подмножества вершин $\{v_1,v_4,v_5,v_6\}$.

$z_{28} = \{e_3,e_5,e_{13},e_{15}\} = \{v_1,v_4,v_5,v_6\} = 8+14+5+7 = 34;$
$z_{29} = \{e_3,e_4,e_{14},e_{15}\} = \{v_1,v_4,v_5,v_6\} = 8+7+11+7 = 33;$
$z_{30} = \{e_4,e_5,e_{13},e_{15}\} = \{v_1,v_4,v_5,v_6\} = 7+14+11+7 = 39;$

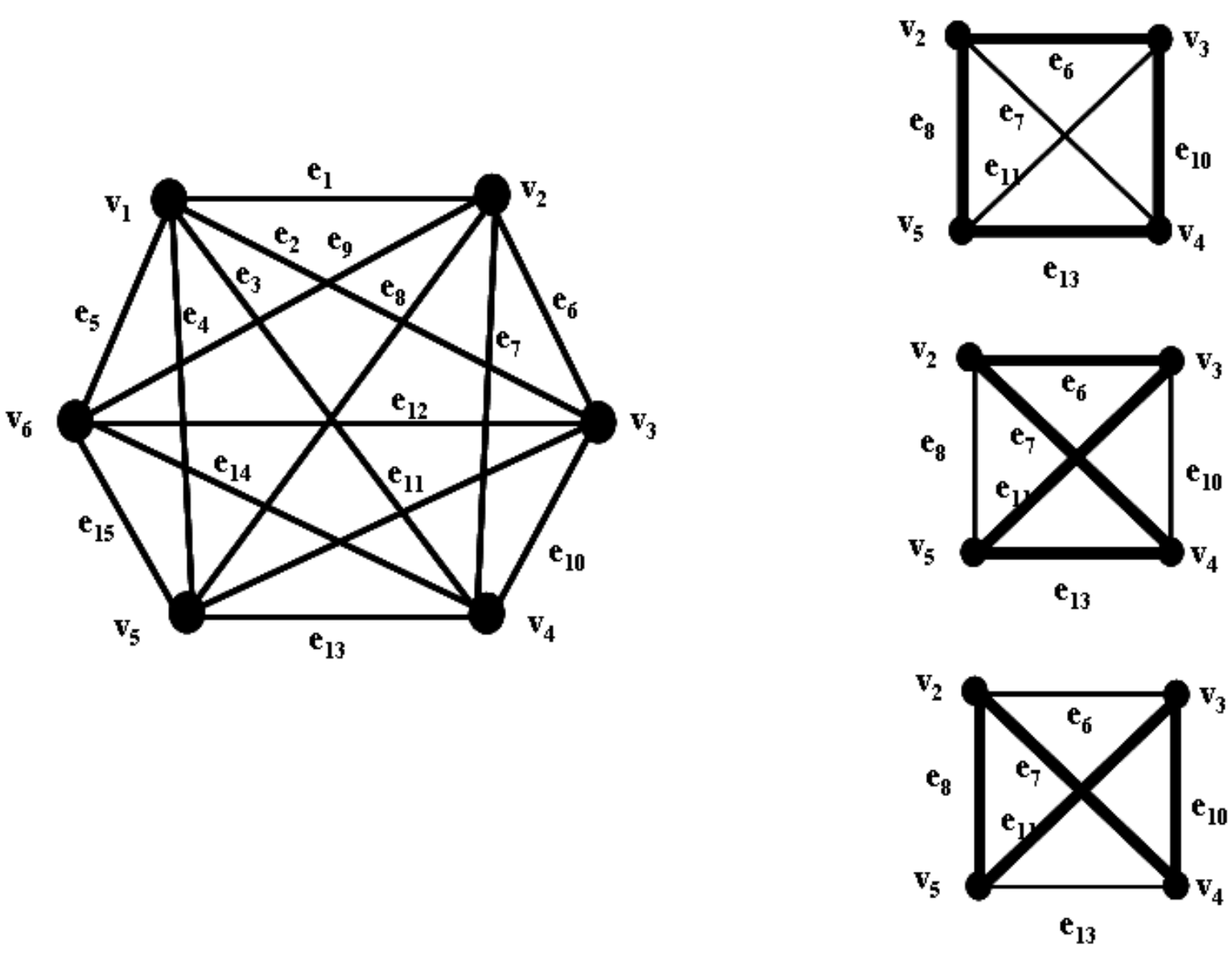

Рис. 3.17. Циклы для подмножества вершин {$v_2$,$v_3$,$v_4$,$v_5$}.

$z_{31}$ = {$e_6$,$e_8$,$e_{10}$,$e_{13}$} = {$v_2$,$v_3$,$v_4$,$v_5$} = 7+7+4+5 = 23;
$z_{32}$ = {$e_6$,$e_7$,$e_{11}$,$e_{13}$} = {$v_2$,$v_3$,$v_4$,$v_5$} = 7+11+3+5 = 26;
$z_{33}$ = {$e_7$,$e_8$,$e_{10}$,$e_{11}$} = {$v_2$,$v_3$,$v_4$,$v_5$} = 11+7+4+3 = 28;

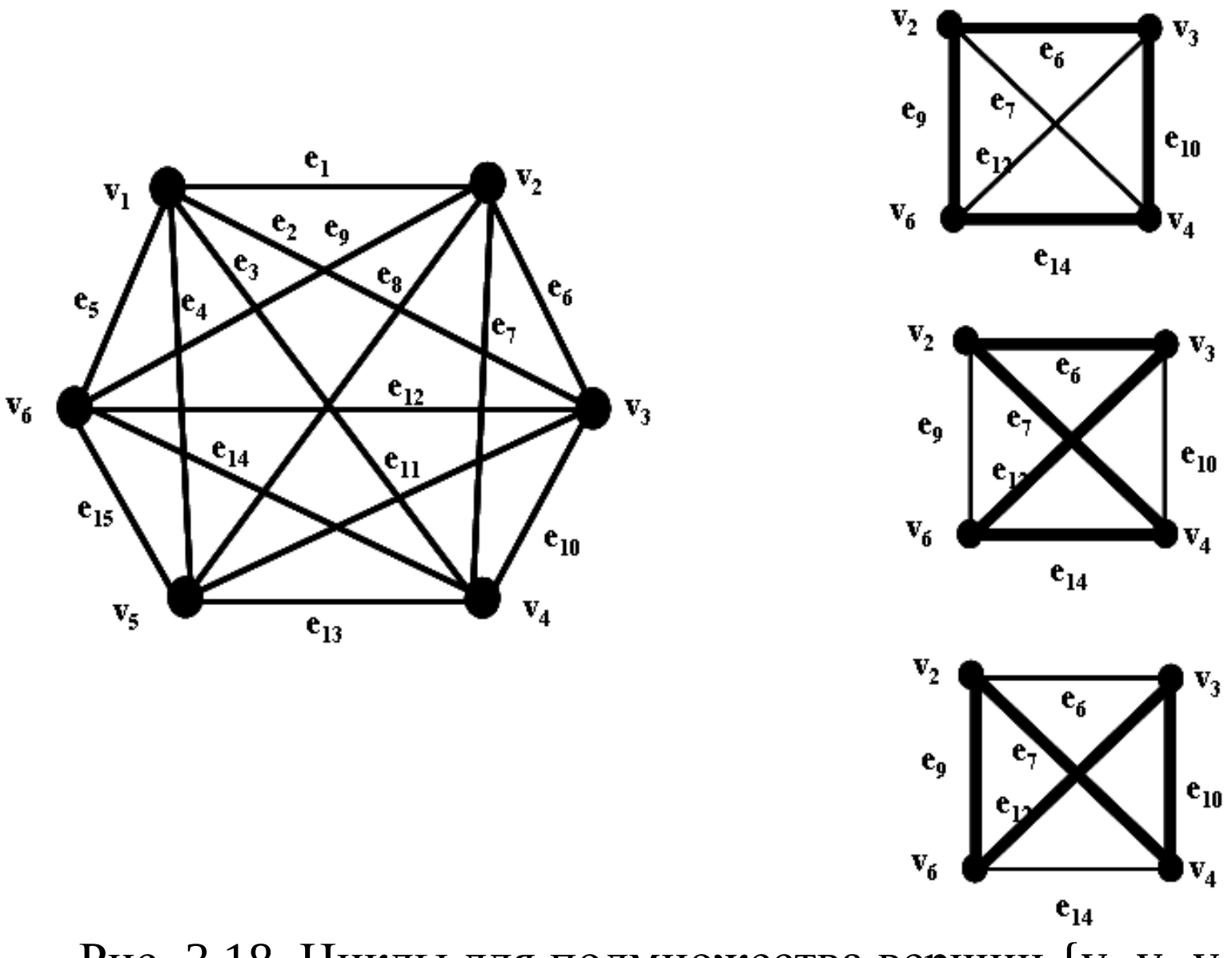

Рис. 3.18. Циклы для подмножества вершин {$v_2$,$v_3$,$v_4$,$v_6$}.

$z_{34}$ = {$e_6$,$e_9$,$e_{10}$,$e_{14}$} = {$v_2$,$v_3$,$v_4$,$v_6$} = 7+10+4+11 = 33;
$z_{35}$ = {$e_6$,$e_7$,$e_{12}$,$e_{14}$} = {$v_2$,$v_3$,$v_4$,$v_6$} = 7+11+10+11 = 39;
$z_{36}$ = {$e_7$,$e_9$,$e_{10}$,$e_{12}$} = {$v_2$,$v_3$,$v_4$,$v_6$} = 11+10+4+10 = 35;

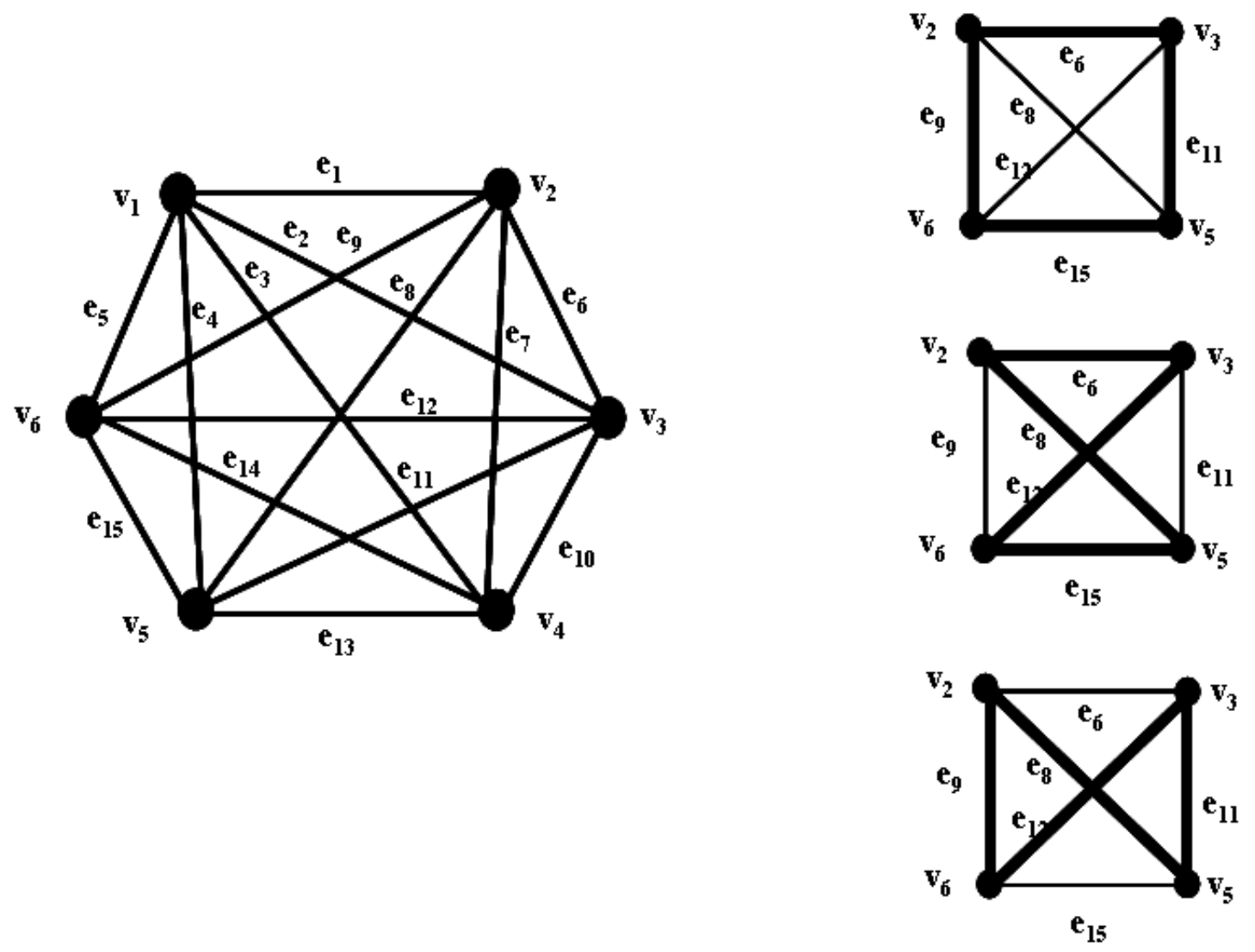

Рис. 3.19. Циклы для подмножества вершин $\{v_2,v_3,v_5,v_6\}$.

$z_{37} = \{e_6,e_9,e_{11},e_{15}\} = \{v_2,v_3,v_5,v_6\} = 7+10+3+7 = 27;$
$z_{38} = \{e_6,e_8,e_{12},e_{15}\} = \{v_2,v_3,v_5,v_6\} = 7+7+10+7 = 31;$
$z_{39} = \{e_8,e_9,e_{11},e_{12}\} = \{v_2,v_3,v_5,v_6\} = 7+10+3+10 = 30;$

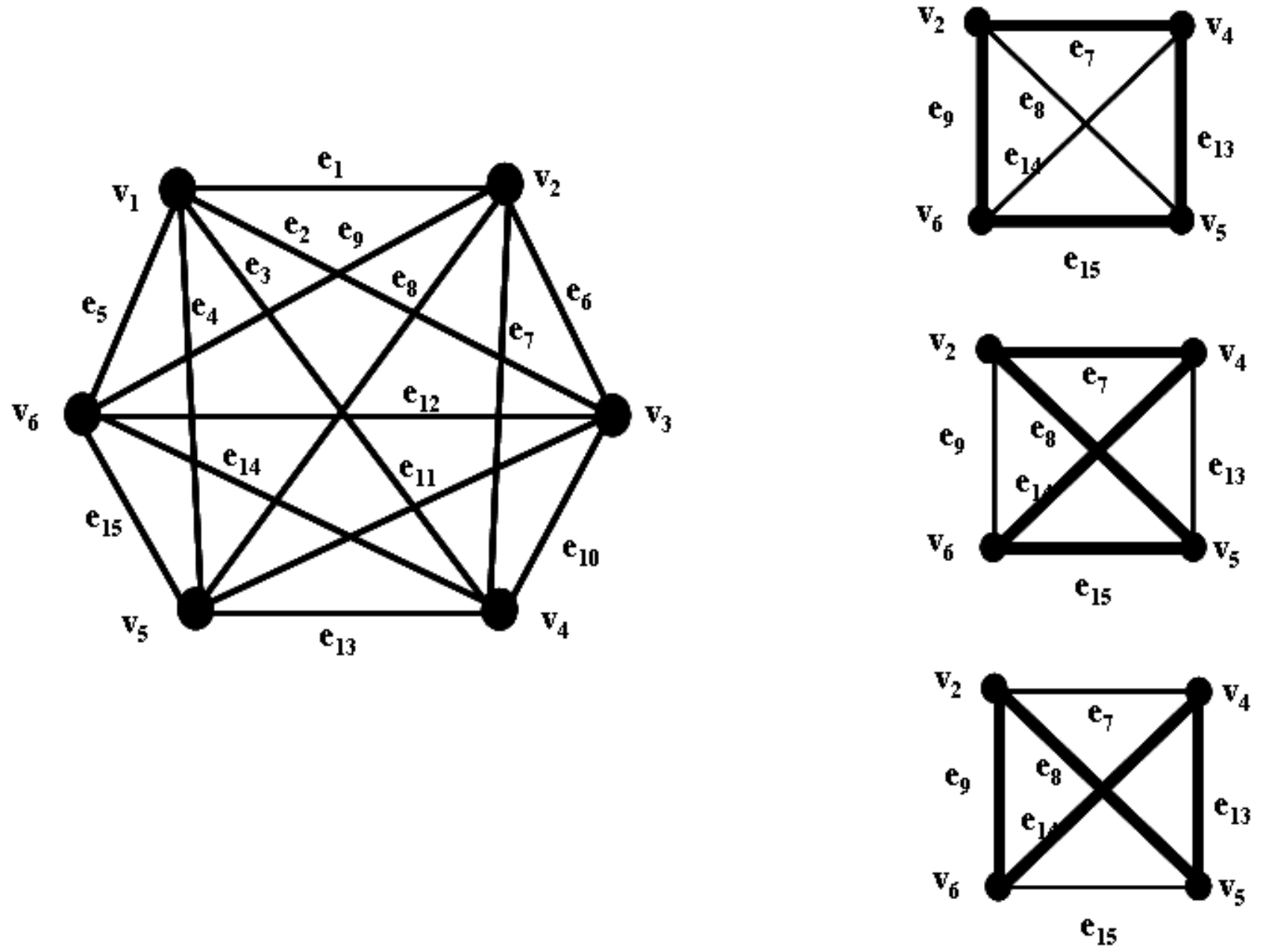

Рис. 3.20. Циклы для подмножества вершин $\{v_2,v_4,v_5,v_6\}$.

$z_{40} = \{e_7,e_9,e_{13},e_{15}\} = \{v_2,v_4,v_5,v_6\} = 11+10+5+7 = 33;$
$z_{41} = \{e_7,e_8,e_{14},e_{15}\} = \{v_2,v_4,v_5,v_6\} = 11+7+11+7 = 36;$
$z_{42} = \{e_8,e_9,e_{13},e_{14}\} = \{v_2,v_4,v_5,v_6\} = 7+10+5+11 = 33;$

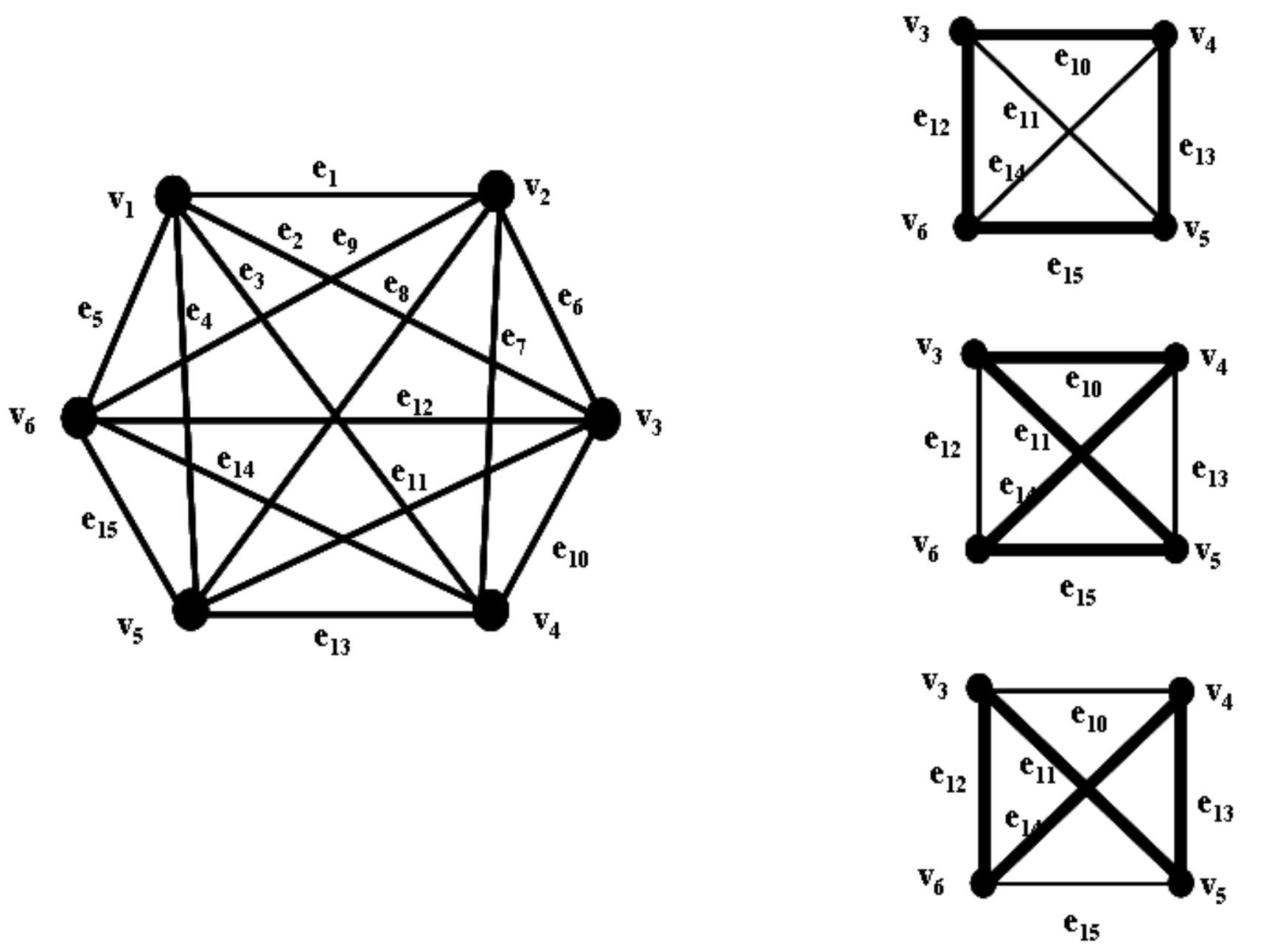

Рис. 3.21. Циклы для подмножества вершин $\{v_3,v_4,v_5,v_6\}$.

$z_{43} = \{e_{10},e_{12},e_{13},e_{15}\} = \{v_3,v_4,v_5,v_6\} = 4+10+5+7 = 26$;
$z_{44} = \{e_{10},e_{11},e_{14},e_{15}\} = \{v_3,v_4,v_5,v_6\} = 4+3+11+7 = 25$;
$z_{45} = \{e_{11},e_{12},e_{13},e_{14}\} = \{v_3,v_4,v_5,v_6\} = 3 +10+5+11 = 29$.

Из каждой тройки циклов выберем цикл длиной четыре с минимальным весом:

$z_1 = \{e_1,e_4,e_6,e_{10}\} = \{v_1,v_2,v_3,v_4\} = 6+7+7+4 = 24$;
$z_5 = \{e_1,e_2,e_8,e_{11}\} = \{v_1,v_2,v_3,v_5\} = 6+4+7+3 = 20$;
$z_7 = \{e_1,e_5,e_6,e_{12}\} = \{v_1,v_2,v_3,v_6\} = 6+4+7+10 = 27$;
$z_{11} = \{e_1,e_3,e_8,e_{13}\} = \{v_1,v_2,v_4,v_5\} = 6+7+8+5 = 26$;
$z_{14} = \{e_1,e_3,e_9,e_{14}\} = \{v_1,v_2,v_4,v_6\} = 6+8+10+7 = 31$;
$z_{17} = \{e_1,e_4,e_9,e_{15}\} = \{v_1,v_2,v_5,v_6\} = 6+7+10+7 = 30$;
$z_{19} = \{e_2,e_4,e_{10},e_{13}\} = \{v_1,v_3,v_4,v_5\} = 4+7+4+5 = 20$;
$z_{20} = \{e_2,e_3,e_{11},e_{13}\} = \{v_1,v_3,v_4,v_5\} = 4+8+3+5 = 20$;
$z_{22} = \{e_2,e_5,e_{10},e_{14}\} = \{v_1,v_3,v_4,v_6\} = 4+8+4+11 = 26$;
$z_{25} = \{e_2,e_5,e_{11},e_{15}\} = \{v_1,v_3,v_5,v_6\} = 4+14+3+7 = 28$;
$z_{26} = \{e_2,e_4,e_{12},e_{15}\} = \{v_1,v_3,v_5,v_6\} = 4+7+10+7 = 28$;
$z_{29} = \{e_3,e_4,e_{14},e_{15}\} = \{v_1,v_4,v_5,v_6\} = 8+7+11+7 = 33$;
$z_{31} = \{e_6,e_8,e_{10},e_{13}\} = \{v_2,v_3,v_4,v_5\} = 7+7+4+5 = 23$;
$z_{34} = \{e_6,e_9,e_{10},e_{14}\} = \{v_2,v_3,v_4,v_6\} = 7+10+4+11 = 33$;
$z_{37} = \{e_6,e_9,e_{11},e_{15}\} = \{v_2,v_3,v_5,v_6\} = 7+10+3+7 = 27$;
$z_{40} = \{e_7,e_9,e_{13},e_{15}\} = \{v_2,v_4,v_5,v_6\} = 11+10+5+7 = 33$;
$z_{42} = \{e_8,e_9,e_{13},e_{14}\} = \{v_2,v_4,v_5,v_6\} = 7+10+5+11 = 33$;
$z_{44} = \{e_{10},e_{11},e_{14},e_{15}\} = \{v_3,v_4,v_5,v_6\} = 4+3+11+7 = 25$.

Выбираем циклы минимальной стоимости:

$z_5 = \{e_1,e_2,e_8,e_{11}\} = \{v_1,v_2,v_3,v_5\} = 6+4+7+3 = 20$;
$z_{19} = \{e_2,e_4,e_{10},e_{13}\} = \{v_1,v_3,v_4,v_5\} = 4+7+4+5 = 20$;
$z_{20} = \{e_2,e_3,e_{11},e_{13}\} = \{v_1,v_3,v_4,v_5\} = 4+8+3+5 = 20$;

Ищем соприкасающиеся циклы для цикла $z_{20}$.

$z_{20} = \{e_2,e_3,e_{11},e_{13}\} = \{v_1,v_3,v_4,v_5\} = 4+8+3+5 = 20$;

$c_1 = \{e_1,e_2,e_6\} \leftrightarrow \{v_1,v_2,v_3\}=6+4+7=17;$
$c_2 = \{e_1,e_3,e_7\} \leftrightarrow \{v_1,v_2,v_4\}=6+8+11=25;$
$c_7 = \{e_2,e_5,e_{12}\} \leftrightarrow \{v_1,v_3,v_6\}=4+14+10=28;$
$c_9 = \{e_3,e_5,e_{14}\} \leftrightarrow \{v_1,v_4,v_6\}=8+14+11=32;$
$c_{12} = \{e_6,e_8,e_{11}\} \leftrightarrow \{v_2,v_3,v_5\}=7+7+3=15;$
$c_{14} = \{e_7,e_8,e_{13}\} \leftrightarrow \{v_2,v_4,v_5\}=11+7+5=23;$
$c_{19} = \{e_{11},e_{12},e_{15}\} \leftrightarrow \{v_3,v_5,v_6\}=3+10+7=20;$
$c_{20} = \{e_{13},e_{14},e_{15}\} \leftrightarrow \{v_4,v_5,v_6\}=5+11+7=23.$

Строим циклы длиной пять:

$z_{20} \oplus c_1 = \{e_2,e_3,e_{11},e_{13}\} \oplus \{e_1,e_2,e_6\}=\{ e_1,e_3,e_6,e_{11},e_{13}\} \leftrightarrow$
$\leftrightarrow \{v_1,v_2,v_3,v_4,v_5\}= 8+3+5+6+7=29$:
$z_{20} \oplus c_2 = \{e_2,e_3,e_{11},e_{13}\} \oplus \{e_1,e_3,e_7\}=\{e_1,e_2,e_7,e_{11},e_{13}\} \leftrightarrow$
$\leftrightarrow \{v_1,v_2,v_3,v_4,v_5\}= 4+3+5+6+11=29;$
$z_{20} \oplus c_7 = \{e_2,e_3,e_{11},e_{13}\} \oplus \{e_2,e_5,e_{12}\}=\{e_3,e_5,e_{11},e_{12},e_{13}\} \leftrightarrow$
$\leftrightarrow \quad \{v_1,v_3,v_4,v_5,v_6\}= 8+3+5+14+10 = 40;$
$z_{20} \oplus c_9 = \{e_2,e_3,e_{11},e_{13}\} \oplus \{e_3,e_4,e_{12}\}=\{e_2,e_4,e_{11},e_{12},e_{13}\} \leftrightarrow$
$\leftrightarrow \{v_1,v_3,v_4,v_5,v_6\}= 4+3+5+14+11 =37;$
$z_{20} \oplus c_{12} = \{e_2,e_3,e_{11},e_{13}\} \oplus \{e_6,e_8,e_{11}\}=\{e_2,e_3,e_6,e_8,e_{13}\} \leftrightarrow$
$\leftrightarrow \{v_1,v_2,v_3,v_4,v_5\}= 4+8+5+7+7=31;$
$z_{20} \oplus c_{14} = \{e_2,e_3,e_{11},e_{13}\} \oplus \{e_7,e_8,e_{13}\}=\{e_2,e_3,e_7,e_8,e_{11}\} \leftrightarrow$
$\leftrightarrow \{v_1,v_2,v_3,v_4,v_5\}= 4+8+3+11+7=33;$
$z_{20} \oplus c_{14} = \{e_2,e_3,e_{11},e_{13}\} \oplus \{e_{11},e_{12},e_{15}\}=\{e_2,e_3,e_{12},e_{13},e_{15}\} \leftrightarrow$
$\leftrightarrow \{v_1,v_2,v_3,v_4,v_5\}= 4+8+5+10+7=34;$
$z_{20} \oplus c_{14} = \{e_2,e_3,e_{11},e_{13}\} \oplus \{e_{13},e_{14},e_{15}\}=\{e_2,e_3,e_{11},e_{14},e_{15}\} \leftrightarrow$
$\leftrightarrow \{v_1,v_3,v_4,v_5,v_6\}= 4+8+3+11+7=33;$

Ищем соприкасающиеся циклы для цикла $z_{19}$.
$z_{19} = \{e_2,e_4,e_{10},e_{13}\} = \{v_1,v_3,v_4,v_5\} = 4+7+4+5 = 20;$

$c_1 = \{e_1,e_2,e_6\} \leftrightarrow \{v_1,v_2,v_3\}=6+4+7=17;$
$c_3 = \{e_1,e_4,e_8\} \leftrightarrow \{v_1,v_2,v_5\}=6+7+7=20;$
$c_7 = \{e_2,e_5,e_{12}\} \leftrightarrow \{v_1,v_3,v_6\}=4+14+10=28;$
$c_{10} = \{e_4,e_5,e_{15}\} \leftrightarrow \{v_1,v_5,v_6\}=7+14+5=28;$
$c_{11} = \{e_6,e_7,e_{10}\} \leftrightarrow \{v_2,v_3,v_4\}=7+11+4=22;$
$c_{14} = \{e_7,e_8,e_{13}\} \leftrightarrow \{v_2,v_4,v_5\}=11+7+5=23;$
$c_{18} = \{e_{10},e_{12},e_{14}\} \leftrightarrow \{v_3,v_4,v_6\}=4+10+11=24;$
$c_{20} = \{e_{13},e_{14},e_{15}\} \leftrightarrow \{v_4,v_5,v_6\}=5+11+7=23.$

Строим циклы длиной пять:

$z_{19} \oplus c_1 = \{e_2,e_4,e_{10},e_{13}\} \oplus \{e_1,e_2,e_6\}=\{e_1,e_4,e_6,e_{10},e_{13}\} \leftrightarrow$
$\leftrightarrow \{v_1,v_2,v_3,v_4,v_5\}=7+4+5+6+7 =29;$
$z_{19} \oplus c_3 = \{e_2,e_4,e_{10},e_{13}\} \oplus \{e_1,e_4,e_8\}=\{e_1,e_2,e_8,e_{10},e_{13}\} \leftrightarrow$
$\leftrightarrow \{v_1,v_2,v_3,v_4,v_5\}=4+4+5+6+7 =26;$
$z_{19} \oplus c_7 = \{e_2,e_4,e_{10},e_{13}\} \oplus \{e_2,e_5,e_{12}\}=\{e_4,e_5,e_{10},e_{12},e_{13}\} \leftrightarrow$
$\leftrightarrow \{v_1,v_3,v_4,v_5,v_6\}=7+4+5+14+10 =40;$
$z_{19} \oplus c_{10} = \{e_2,e_4,e_{10},e_{13}\} \oplus \{e_4,e_5,e_{15}\}=\{e_2,e_5,e_{10},e_{13},e_{15}\} \leftrightarrow$
$\leftrightarrow \{v_1,v_3,v_4,v_5,v_6\}= 4+4+5+14+5 =32;$
$z_{19} \oplus c_{11} = \{e_2,e_4,e_{10},e_{13}\} \oplus \{e_6,e_7,e_{10}\}=\{e_2,e_4,e_6,e_7,e_{13}\} \leftrightarrow$

$\leftrightarrow \{v_1,v_2,v_3,v_4,v_5\}= 4+7+5+7+11=34;$

$z_{19}\oplus\ c_{14} = \{e_2,e_4,e_{10},e_{13}\}\oplus\{e_7,e_8,e_{13}\}=\{e_2,e_4,e_7,e_8,e_{10}\}\leftrightarrow$

$\leftrightarrow \{v_1,v_2,v_3,v_4,v_5\}= 4+7+4+11+7 =33;$

$z_{19}\oplus\ c_{18} = \{e_2,e_4,e_{10},e_{13}\}\oplus\{e_{10},e_{12},e_{14}\}=\{e_2,e_4,e_{12},e_{13},e_{14}\}\leftrightarrow$

$\leftrightarrow \{v_1,v_3,v_4,v_5,v_6\}= 7+4+5+10+11 =37;$

$z_{19}\oplus\ c_{20} = \{e_2,e_4,e_{10},e_{13}\}\oplus\{e_{13},e_{14},e_{15}\}=\{e_2,e_4,e_{10},e_{14},e_{15}\}\leftrightarrow$

$\leftrightarrow \{v_1,v_3,v_4,v_5,v_6\}= 4+7+4+11+7 =33;$

Ищем соприкасающиеся циклы для цикла $z_5$.

$z_5 = \{e_1,e_2,e_8,e_{11}\} = \{v_1,v_2,v_3,v_5\} = 6+4+7+3 = 20;$

$c_2 = \{e_1,e_3,e_7\}\leftrightarrow\{v_1,v_2,v_4\}=6+8+11=25;$

$c_4 = \{e_1,e_5,e_9\}\leftrightarrow\{v_1,v_2,v_6\}=6+14+10=30;$

$c_5 = \{e_2,e_3,e_{10}\}\leftrightarrow\{v_1,v_3,v_4\}=4+8+4=16;$

$c_7 = \{e_2,e_5,e_{12}\}\leftrightarrow\{v_1,v_3,v_6\}=4+14+10=28;$

$c_{14} = \{e_7,e_8,e_{13}\}\leftrightarrow\{v_2,v_4,v_5\}=11+7+5=23;$

$c_{16} = \{e_8,e_9,e_{15}\}\leftrightarrow\{v_2,v_5,v_6\}=7+10+7=24;$

$c_{17} = \{e_{10},e_{11},e_{13}\}\leftrightarrow\{v_3,v_4,v_5\}=4+3+5=12;$

$c_{19} = \{e_{11},e_{12},e_{15}\}\leftrightarrow\{v_3,v_5,v_6\}=3+10+7=20;$

Строим циклы длиной пять:

$z_5 = \{e_1,e_2,e_8,e_{11}\} = \{v_1,v_2,v_3,v_5\} = 6+4+7+3 = 20;$

$z_5\oplus\ c_2 = \{e_1,e_2,e_8,e_{11}\}\oplus\{e_1,e_3,e_7\}=\{e_2,e_3,e_7,e_8,e_{11}\}\leftrightarrow$

$\leftrightarrow \{v_1,v_2,v_3,v_4,v_5\}= 4+7+3+8+11=33;$

$z_5\oplus\ c_4 = \{e_1,e_2,e_8,e_{11}\}\oplus\{e_1,e_5,e_9\}=\{e_2,e_5,e_8,e_9,e_{11}\}\leftrightarrow$

$\leftrightarrow \{v_1,v_2,v_3,v_5,v_6\}= 4+7+3+14+10=38;$

$z_5\oplus\ c_5 = \{e_1,e_2,e_8,e_{11}\}\oplus\{e_2,e_3,e_{10}\}=\{e_1,e_3,e_8,e_{10},e_{11}\}\leftrightarrow$

$\leftrightarrow \{v_1,v_2,v_3,v_4,v_5\}= 6+7+3+8+4=28;$

$z_5\oplus\ c_7 = \{e_1,e_2,e_8,e_{11}\}\oplus\{e_2,e_5,e_{12}\}=\{e_1,e_5,e_8,e_{11},e_{12}\}\leftrightarrow$

$\leftrightarrow \{v_1,v_2,v_3,v_4,v_5\}= 6+7+3+14+10=40;$

$z_5\oplus\ c_{14} = \{e_1,e_2,e_8,e_{11}\}\oplus\{e_7,e_8,e_{13}\}=\{e_1,e_2,e_7,e_{11},e_{13}\}\leftrightarrow$

$\leftrightarrow \{v_1,v_2,v_3,v_4,v_5\}= 6+4+3+11+5=29;$

$z_5\oplus\ c_{16} = \{e_1,e_2,e_8,e_{11}\}\oplus\{e_8,e_9,e_{15}\}=\{e_1,e_2,e_9,e_{11},e_{15}\}\leftrightarrow$

$\leftrightarrow \{v_1,v_2,v_3,v_5,v_6\}= 6+4+3+10+7=30;$

$z_5\oplus\ c_{17} = \{e_1,e_2,e_8,e_{11}\}\oplus\{e_{10},e_{11},e_{13}\}=\{e_1,e_2,e_8,e_{10},e_{13}\}\leftrightarrow$

$\leftrightarrow \{v_1,v_2,v_3,v_4,v_5\}= 6+4+7+4+5=26;$

$z_5\oplus\ c_{19} = \{e_1,e_2,e_8,e_{11}\}\oplus\{e_{11},e_{12},e_{15}\}=\{e_1,e_2,e_8,e_{12},e_{15}\}\leftrightarrow$

$\leftrightarrow \{v_1,v_2,v_3,v_5,v_6\}= 6+4+7+10+7=34;$

Выбираем циклы длиной пять с минимальной стоимостью. Это два дубль-цикла.

$z_{19}\oplus\ c_3 = \{e_1,e_2,e_8,e_{10},e_{13}\}\leftrightarrow\{v_1,v_2,v_3,v_4,v_5\}=4+4+5+6+7 = 26;$

$z_5\oplus\ c_{17} = \{e_1,e_2,e_8,e_{10},e_{13}\}\leftrightarrow \{v_1,v_2,v_3,v_4,v_5\}= 6+4+7+4+5 = 26;$

Ищем соприкасающиеся циклы для цикла $z_{19}\oplus\ c_3$.

$z_{19}\oplus\ c_3 = \{e_1,e_2,e_8,e_{10},e_{13}\}\leftrightarrow\{v_1,v_2,v_3,v_4,v_5\}=4+4+5+6+7 = 26;$

$c_4 = \{e_1,e_5,e_9\}\leftrightarrow\{v_1,v_2,v_6\}=6+14+10=30;$

$c_7 = \{e_2,e_5,e_{12}\}\leftrightarrow\{v_1,v_3,v_6\}=4+14+10=28;$

$c_{16} = \{e_8,e_9,e_{15}\}\leftrightarrow\{v_2,v_5,v_6\}=7+10+7=24;$

$c_{18} = \{e_{10},e_{12},e_{14}\}\leftrightarrow\{v_3,v_4,v_6\}=4+10+11=24;$

$c_{20} = \{e_{13},e_{14},e_{15}\} \leftrightarrow \{v_4,v_5,v_6\}=5+11+7=23.$

Строим циклы длиной шесть:

$z_{19} \oplus c_3 \oplus c_4 = \{e_1,e_2,e_8,e_{10},e_{13}\} \oplus \{e_1,e_5,e_9\}=\{e_2,e_5,e_8,e_9,e_{10},e_{13}\} \leftrightarrow$
$\leftrightarrow \{v_1,v_2,v_3,v_4,v_5,v_6\}= 4+4+5+7+14+10 = 44;$
$z_{19} \oplus c_3 \oplus c_7 = \{e_1,e_2,e_8,e_{10},e_{13}\} \oplus \{e_2,e_5,e_{12}\}=\{e_1,e_5,e_8,e_{10},e_{12},e_{13}\} \leftrightarrow$
$\leftrightarrow \{v_1,v_2,v_3,v_4,v_5,v_6\}=4+5+6+7+14+10 = 46;$
$z_{19} \oplus c_3 \oplus c_{16} = \{e_1,e_2,e_8,e_{10},e_{13}\} \oplus \{e_8,e_9,e_{15}\}=\{e_1,e_2,e_9,e_{10},e_{13},e_{15}\} \leftrightarrow$
$\leftrightarrow \{v_1,v_2,v_3,v_4,v_5,v_6\}=4+4+5+6+10+7 = 36;$
$z_{19} \oplus c_3 \oplus c_{18} = \{e_1,e_2,e_8,e_{10},e_{13}\} \oplus \{e_{10},e_{12},e_{14}\}=\{e_1,e_2,e_8,e_{12},e_{13},e_{14}\} \leftrightarrow$
$\leftrightarrow \{v_1,v_2,v_3,v_4,v_5,v_6\}= 4+5+6+7+10+11 = 43;$
$z_{19} \oplus c_3 \oplus c_{18} = \{e_1,e_2,e_8,e_{10},e_{13}\} \oplus \{e_{13},e_{14},e_{15}\}=\{e_1,e_2,e_8,e_{10},e_{14},e_{15}\} \leftrightarrow$
$\leftrightarrow \{v_1,v_2,v_3,v_4,v_5,v_6\}=4+4+6+7+11+7 = 39;$

Выбираем оптимальный маршрут из шести ребер.

$z_{19} \oplus c_3 \oplus c_{16} = \{e_1,e_2,e_8,e_{10},e_{13}\} \oplus \{e_8,e_9,e_{15}\}=\{e_1,e_2,e_9,e_{10},e_{13},e_{15}\} \leftrightarrow$
$\leftrightarrow \{v_1,v_2,v_3,v_4,v_5,v_6\}=4+4+5+6+10+7 = 36;$

### 3.6. Вычислительная сложность алгоритма решения задачи коммивояжёра

Для формирования циклов длиной четыре, согласно формулы (3.2), нужно затратить

$$k_4 = 3\frac{n(n-1)(n-2)(n-3)}{4!} \approx \frac{n^4}{8}. \qquad (3.11)$$

Для выбора соприкасающихся изометрических циклов с целью построения оптимального маршрута длиной пять нужно перебрать все множество изометрических циклов, количество которых определяется формулой (3.1)

$$k_c = \frac{n(n-1)(n-2)}{6} \approx \frac{n^3}{6}. \qquad (3.12)$$

Вычислительную сложность построения маршрута из пяти рёбер можно определить как $f(n)=\frac{n^4}{8}+\frac{n^3}{6}$.

Вычислительная сложность построения маршрута из шести рёбер определяется как $f(n)= \frac{n^4}{8}+\frac{n^3}{6}+\frac{n^3}{6}= \frac{n^4}{8}+\frac{n^3}{3}$.

Вычислительную сложность построения маршрута из количества рёбер равного количеству вершин можно определить как

$$f(n)=\frac{n^4}{8}+\frac{(n-4)n^3}{6}=\frac{3n^4+4(n^4-4n^3)}{24}=\frac{1}{24}(7n^4-16n^3). \qquad (3.13)$$

Исходя из последнего выражения вычислительную сложность алгоритма можно записать как о($n^4$).

Таким образом, вычислительная сложность предсталенного алгоритма решения задачи коммивояжёра носит полиномиальный характер.

### 3.7. Алгоритм вычисления оптимального пути в симметричной задаче коммивояжёра

Опишем алгоритм построения гамильтонового цикла в полном графе.

**Инициализация.** Выделяем в полном графе G(E,V) систему изометрических циклов $C_\tau$. Вводим матрицу расстояний между вершинами. Записываем каждый цикл $c_i$ в виде подмножества рёбер $E_i$ и подмножества вершин $V_i$. Каждому изометрическому циклу ставим в соответствие его вес как аддитивную сумму его ребер.

**Шаг 1. [Построение подмножества циклов длиной четыре].** Определяем для каждой четвёрки вершин оптимальный маршрут. Строим подмножество оптимальных маршрутов.

**Шаг 2. [Выбор оптимального маршрута длиной четыре].** Осуществляем выбор из подмножества циклов оптимального маршрута длиной четыре.

**Шаг 3. [Выбор соприкасающихся изометрических циклов].** Для выбранного простого цикла определяем подмножество соприкасающихся циклов.

**Шаг 4. [Суммирование циклов].** Осуществляем кольцевое суммирование предыдущего простого цикла и соприкасающихся циклов.

**Шаг 5**. [**Выбор оптимального маршрута**]. Выбираем из суммарного подмножества циклов оптимальный маршрут и добавляем новую вершину в подмножество вершин простого цикла.

**Шаг 6. [Проверка подмножества вершин].** Если в подмножество вершин построенного простого цикла включены все вершины графа, то конец работы алгоритма. Иначе идти на шаг 3.

В результате получаем оптимальный маршрут в заданном взвешенном графе.

**Замечание**: для минимизации объёма памяти в процессе расчёта можно хранить только простые циклы с минимальными параметрами.

## Заключение

В работе представлен полиномиальный алгоритм для решения симметричной задачи коммивояжёра с вычислительной сложностью $o(n^4)$.

Для демонстрации метода приведены примеры расчёта оптимального маршрута для полных графов $K_4$, $K_5$ и $K_6$ без применения средств вычислительной техники.

## Список используемой литературы